\documentclass[12pt]{article}

\usepackage{amsmath,amsthm,amssymb}
\allowdisplaybreaks
\usepackage[colorlinks,
linkcolor=red,
anchorcolor=blue,
citecolor=magenta
]{hyperref}

\theoremstyle{plain}
\newtheorem{theorem}{Theorem}[section]
\newtheorem{lemma}{Lemma}[section]
\newtheorem{proposition}{Proposition}[section]
\newtheorem{corollary}[theorem]{Corollary}

\theoremstyle{definition}
\newtheorem{definition}{Definition}[section]
\newtheorem{remark}{Remark}[section]

\numberwithin{equation}{section}

\newcommand{\vs}{\vspace}

\renewcommand{\div}{\operatorname{div}}
\begin{document}

	\title{Global Well-Posedness of Two-Dimensional Compressible Navier--Stokes Equations with Vacuum under a Scaling-Invariant Smallness Condition}
	
	\author{
		Boqiang L\"{u},~~ Xin Song\thanks{Email: lvbq86@ncu.edu.cn \& lvbq86@163.com (B.~L\"{u}),
			songxin2357@163.com (X.~Song).}\\[4pt]
		{\small School of Mathematics and Computer Sciences,}\\
		{\small \& Institute of Mathematics and Interdisciplinary Sciences,}\\
		{\small Nanchang University, Nanchang 330031, P.~R.~China}
	}

	\date{}

	\maketitle \vs{-.7cm}
	
	\begin{abstract}
		In this paper, we consider the Cauchy problem for the two-dimensional
		isentropic compressible Navier--Stokes equations with vacuum as far-field density.
		For the initial data that satisfy suitable regularity conditions,
		we establish the global well-posedness of strong solutions, provided that the scaling-invariant initial
		quantity \(
		\bar{\rho}
		\big(
		\|\rho_0\|_{L^1}
		+
		\bar{\rho}^{-2\gamma} \|\rho_0\|_{L^1}^{-1}
		+
		\bar{\rho}^{\frac{\alpha(\gamma+1)}{2}} \|\rho_0 |x|^\alpha\|_{L^1}
		\big)
		(1+\bar{\rho} E_0)^{3/2}
		\big(
		\|\nabla u_0\|_{L^2}^2 + \|\rho_0^\gamma\|_{L^2}^2
		\big)
		\) is sufficiently small. Moreover, if the initial
		data possess additional higher order regularity and satisfy the corresponding
		compatibility condition, the global strong solution becomes a  classical
		one.
		The initial density is allowed to contain vacuum, and vacuum at infinity is included.
		In particular, compactly supported initial density is admissible. The result provides
		a scaling-invariant smallness framework for the global well-posedness theory of
		two-dimensional compressible barotropic flows with far-field vacuum.
	\end{abstract}


	\section{ Introduction}
	We consider the Navier--Stokes equations
	\begin{equation}\label{1.1}
		\left\{\begin{array}{l}
			\rho_t+\operatorname{div}(\rho u)=0, \\[2mm]
			(\rho u)_t+\operatorname{div}(\rho u \otimes u)
			-\mu \Delta u-(\mu+\lambda) \nabla \operatorname{div} u+\nabla P(\rho)=0
		\end{array}\right.
	\end{equation}
	for viscous compressible barotropic flows in two spatial dimensions. Here
	\(t\ge0\) and \(x\in\mathbb R^2\) are the time and space variables, while
	\(\rho=\rho(x,t)\), \(u=(u^1,u^2)(x,t)\), and
	\begin{equation}\label{1.2.1}
		P(\rho)=R\rho^\gamma \qquad (R>0,\ \gamma>1)
	\end{equation}
	denote the density, velocity and pressure, respectively. Without loss of generality,
	we take \(R=1\). The constant viscosity coefficients satisfy the physical condition
	\begin{equation}\label{1.3.1}
		\mu>0,\qquad \mu+\lambda\ge0 .
	\end{equation}
	We study the Cauchy problem for \eqref{1.1} with
	\((\rho,u)\) vanishing at infinity (in some weak sense) with given initial data
	\(\rho_0\) and \(u_0\), as
	\begin{equation}\label{1.4.1}
		\rho(x,0)=\rho_0(x),\qquad \rho u(x,0)=\rho_0u_0(x).
	\end{equation}
	
	There is a vast literature on the global-in-time existence and large-time behavior of
	solutions to \eqref{1.1}. In the absence of vacuum, Matsumura--Nishida
	\cite{MatsumuraNishida1980} established global classical solutions for small
	perturbations of non-vacuum equilibrium states. Global solutions with discontinuous
	initial data and small oscillations around a non-vacuum state were studied by Hoff
	\cite{Hoff1995}. In the framework of critical Besov spaces, Danchin
	\cite{Danchin2000} proved the global well-posedness of the compressible
	Navier--Stokes system for small initial velocity and initial density close to a
	positive constant. For finite-energy initial data allowing vacuum, the fundamental
	global weak-solution theory was
	established by Lions \cite{Lions1998} for the barotropic pressure
	\(P(\rho)=\rho^\gamma\) with
	\(\gamma\ge 3n/(n+2)\), and was further developed by Feireisl--Novotn\'y--%
	Petzeltov\'a \cite{Feireisl2001}, who lowered the condition to
	\(\gamma>n/2\). In particular, in two space dimensions this includes all
	\(\gamma>1\). However, the
	regularity and uniqueness of such weak solutions remain largely open, especially when
	the density may vanish.
	
	When vacuum is present, the analysis becomes much more subtle. The momentum equation
	loses its uniform parabolic character in vacuum regions, and even the local theory
	requires a careful formulation of the initial data and compatibility conditions. The
	local strong solution theory for nonnegative initial density was developed in
	\cite{ChoeKim2003,ChoChoeKim2004,ChoKim2006}. For the Cauchy problem with far-field
	vacuum, there is an additional difficulty: one cannot use a positive lower bound of
	the density to convert the energy bounds into unweighted estimates for the velocity,
	and this loss is particularly relevant for compactly supported density.
	Xin \cite{Xin1998} shows that, for the full compressible
	Navier--Stokes system without heat conduction, compactly supported initial density
	leads to finite-time blow-up of smooth finite-energy solutions in any space
	dimension; the corresponding compact-support and finite-energy assumptions were
	later weakened by Xin--Yan \cite{XinYan2013} for a large class of initial data
	containing vacuum. This makes global regularity in the presence of vacuum and compactly
	supported density a delicate issue.
	
	A breakthrough for the three-dimensional isentropic system was obtained by Huang--Li--
	Xin \cite{HuangLiXin2012}, who proved global well-posedness of classical
	solutions with small initial energy but possibly large oscillations and vacuum. Their
	result showed that interior vacuum and far-field vacuum may coexist with global
	regularity, provided the solution is sought in a suitable homogeneous framework and the
	initial energy is small. For the full heat-conductive compressible Navier--Stokes
	system, related global large-data results with vacuum at infinity were obtained by Wen--%
	Zhu \cite{WenZhu2013}.
	
	For the two-dimensional Cauchy problem with vacuum as far-field density, the
	local-in-time strong solution theory in the weighted framework used below follows
	from \cite{Liang2014 2D}. Starting from this local theory, Li--Xin
	\cite{Li-xin-2019} developed a global strong and classical solution theory under a
	small initial energy assumption. Their solution class and global a priori estimates
	are formulated in a homogeneous Sobolev framework. One of the key observations in
	\cite{Li-xin-2019} is that the two-dimensional far-field vacuum problem cannot be
	treated simply by repeating the three-dimensional argument of
	\cite{HuangLiXin2012}, since that argument uses the a priori \(L^6\)-bound on the
	velocity. By contrast, for the two-dimensional problem with positive far-field
	density, the \(L^p\)-norm of \(u\) can be bounded by
	\(\|\sqrt\rho u\|_{L^2}\) and \(\|\nabla u\|_{L^2}\), whereas this argument is no
	longer available in the same way when the far-field density is vacuum. Li--Xin
	overcome this difficulty by combining weighted estimates for the density,
	time-weighted energy estimates, decay rates of the pressure and the velocity gradient; in
	particular, compactly supported initial density is allowed. In the present paper, we
	follow the decay mechanism of \cite{Li-xin-2019}, but the lower-order global
	estimates are not organized through their homogeneous Sobolev framework. Instead, they
	are arranged by the weighted density information, the normalization of \(\bar\rho\),
	and the scaling-invariant initial quantity introduced below. For comparison, in
	two-dimensional models with density-dependent viscosities, Vaigant--Kazhikhov
	\cite{VaigantKazhikhov1995} obtained a celebrated global existence result for large
	initial data, and Jiu--Wang--Xin \cite{JiuWangXin2016} extended this direction to
	allow vacuum and large data. The viscosity structure in these papers is different
	from the constant-viscosity system considered here.
	
	The small-energy condition in \cite{Li-xin-2019}, however, is not formulated in a
	scaling-invariant way. Scaling is a natural guide for
	compressible Navier--Stokes equations. For the barotropic system, Lei--Xin
	\cite{lx19b} formulated the scaling law
	\[
	\rho_\ell(x,t)=\ell^{1/\gamma}
	\rho\bigl(\ell^{\frac{\gamma+1}{2\gamma}}x,\ell t\bigr),
	\qquad
	u_\ell(x,t)=\ell^{\frac{\gamma-1}{2\gamma}}
	u\bigl(\ell^{\frac{\gamma+1}{2\gamma}}x,\ell t\bigr),
	\qquad \ell>0,
	\]
	which leaves the compressible Navier--Stokes equations invariant. In \cite{lx19b},
	this scaling is used to assign dimensions to the space-time variables and the
	unknowns, to exclude certain type-I singularities, and, when \(\gamma>1\), to
	construct an explicit self-similar type-II blow-up solution.
	For smallness assumptions, the same scaling gives a useful test: a non-invariant
	quantity can be changed by rescaling the same initial data, and hence its size may
	depend on the chosen normalization rather than only on the intrinsic size of the data.
	For the heat-conductive compressible Navier--Stokes equations with vacuum, J. Li
	\cite{Li2020ARMA} established global strong solutions under the smallness of a quantity
	invariant under the natural scaling of the full system, under the additional viscosity
	condition \(2\mu>\lambda\); notably, the total mass may be finite or infinite. More
	recently, Wen \cite{Wen2025} considered the three-dimensional compressible
	Navier--Stokes equations with far-field vacuum and proved global well-posedness for the
	isentropic flow under the smallness of the scaling-invariant quantity
	\[
	\bar\rho^3
	E_0
	\left(
	\|\nabla u_0\|_{L^2}^2+\|\rho_0^\gamma\|_{L^2}^2
	\right)
	\left(1+\bar\rho^{3+\gamma}E_0^2\right),
	\]
	together with the corresponding result for the Navier--Stokes--Fourier system with a
	different scaling-invariant quantity. The new dissipative estimate in \cite{Wen2025}
	is used to handle the supercritical convective term without imposing extra restrictions
	on the viscosity coefficients. This scaling-invariant viewpoint has also appeared in
	recent preprints. Xu--Zhong \cite{XuZhong2026coefficient} studied the
	three-dimensional compressible Navier--Stokes equations with far-field vacuum under
	a coefficient-coupled scaling-invariant smallness condition. They established global
	well-posedness and decay rates for strong solutions, with a smallness threshold
	independent of the initial data and the known system parameters; their result also
	does not require the usual compatibility condition on the initial data. For domains
	with boundaries, Lin Xu--Zhong \cite{XuZhong2026halfspace} considered the
	three-dimensional isentropic system in the upper half-space with slip boundary
	conditions. They proved global existence and uniqueness of strong solutions with
	vacuum and large oscillations under a scaling-invariant smallness assumption. The
	natural scaling preserves both the half-space geometry and the slip boundary
	conditions, so this work gives a half-space counterpart of the whole-space theory.
	Both results are adjacent to the present problem, but they concern three-dimensional
	equations (and, in the latter case, a boundary-value problem), whereas here we study
	the two-dimensional whole-space Cauchy problem with weighted far-field vacuum.
	For related MHD systems with vacuum, Li--Xu--Zhang
	\cite{LiXuZhang2013} established a global classical solution theory with small energy
	and large oscillations. Liu--Zhong \cite{LiuZhong2022} treated the
	three-dimensional full compressible MHD system with far-field vacuum under a smallness
	condition and the additional restriction \(3\mu>\lambda\). Xu--Zhong
	\cite{XuZhong2026} recently obtained a global strong solution result for the
	three-dimensional non-isentropic compressible MHD equations under a scaling-invariant
	smallness condition, removing the artificial restriction \(3\mu>\lambda\) imposed in
	\cite{LiuZhong2022}. The precise invariant quantities in these results differ
	according to the scaling laws and the main a priori estimates of the corresponding
	models.
	
	The aim of the present paper is to establish a two-dimensional scaling-invariant
	smallness result in the weighted far-field vacuum setting of Li--Xin
	\cite{Li-xin-2019}. The scaling recalled in \cite{lx19b} is used here not for the
	study of singularity rates, but to organize the initial smallness condition and the
	global a priori estimates. The problem considered here is also different from the
	heat-conductive full system in \cite{Li2020ARMA} and from the three-dimensional
	far-field vacuum theory in \cite{Wen2025}. In two dimensions, the estimates must
	retain the weighted density structure and the expansion control of the mass
	distribution from \cite{Li-xin-2019}. At the same time, the lower-order global estimates are not
	organized in the homogeneous Sobolev framework used in \cite{Li-xin-2019}, but through
	the normalization of \(\bar\rho\) and a scaling-invariant initial quantity. The
	quantity \(I_0\) below is designed to be adapted simultaneously to the
	two-dimensional scaling, the far-field vacuum, and the weighted density estimates.
	
	More precisely, we introduce below the quantity \(I_0\), where the factor \(\Xi_0\)
	contains the mass \(M_0\), its reciprocal in the scaled form, and the weighted moment
	\(\|\rho_0|x|^\alpha\|_{L^1}\). The definition of \(I_0\) is chosen so that it is
	invariant under the natural scaling of \eqref{1.1} in \(\mathbb R^2\). The smallness
	threshold for \(I_0\) depends only on \(\mu,\lambda,\gamma,q,\alpha\), and is
	independent of the initial data, time, and the upper bound \(\bar\rho\) of the initial
	density. A key step of the proof is to use the scaling to normalize the density upper
	bound, obtain estimates in the normalized variables, and then return to the original
	variables. In this way, the constants appearing in the continuation argument are
	organized through the scaling-invariant quantity rather than through a non-invariant
	small-energy condition.
	
	The proof is inspired by the estimates in \cite{Li-xin-2019}, but the lower-order global estimates are reorganized to respect the scaling structure.
	We use weighted Sobolev inequalities to control the velocity in the far-field vacuum setting, effective viscous flux and vorticity estimates
	to recover derivatives of \(u\), time-weighted estimates to derive the large-time decay
	needed in two dimensions, and a Zlotnik-type argument to obtain a uniform upper bound
	for the density. The scaling-invariant smallness of \(I_0\) closes these estimates
	globally in time. A basic difference from the three-dimensional theories in
	\cite{HuangLiXin2012,Wen2025} is that the homogeneous Sobolev embedding no longer
	gives an a priori \(L^6\)-bound for the velocity in the present two-dimensional
	far-field vacuum problem. The mass expansion estimate and the weighted density moment
	are therefore used to recover the required velocity control, while the decay of \(P\)
	and \(\nabla u\) is obtained from time-weighted estimates in normalized scaling variables;
	these arguments are related to those in Li--Xin \cite{Li-xin-2019}. The resulting theorem
	allows the density to contain vacuum both inside the domain and at infinity; in
	particular, compactly supported initial density is admissible. If the initial data have
	additional higher order regularity and satisfy the compatibility condition, the global
	strong solution obtained here becomes a global classical solution.
	
	Before stating the main results, we fix the notation used throughout the paper. For
	\(R>0\), let
	\[
	B_R:=\{x\in\mathbb R^2:\ |x|<R\}.
	\]
	For \(1\le q\le\infty\) and any positive integer \(m\), we write
	\(L^q=L^q(\mathbb R^2)\), \(W^{m,q}=W^{m,q}(\mathbb R^2)\), and
	\(H^m=W^{m,2}\). The notation \(\|\cdot\|_q\) denotes the \(L^q\)-norm on
	\(\mathbb R^2\). For homogeneous Sobolev spaces, set
	\[
	D^{k,r}:=
	\left\{
	v\in L^1_{\operatorname{loc}}(\mathbb R^2):\
	\|\nabla^k v\|_{L^r}<\infty
	\right\},
	\qquad
	D^k:=D^{k,2}.
	\]
	We also use
	\[
	\int f\,dx:=\int_{\mathbb R^2}f\,dx .
	\]
	For vector-valued functions, the same notation is used componentwise. We now
	recall the notion of strong solutions.
	
	\begin{definition}\label{def1.1}
		If all derivatives involved in \eqref{1.1} for \((\rho,u)\) are regular
		distributions, and equations \eqref{1.1} hold almost everywhere in
		\(\mathbb R^2\times(0,T)\), then \((\rho,u)\) is called a strong solution to
		\eqref{1.1}.
	\end{definition}
	
	Assume that the initial mass is finite and positive, namely
	\begin{equation}\label{mass-M0}
		M_0:=\int_{\mathbb R^2}\rho_0\,dx\in(0,+\infty).
	\end{equation}
	Then there exists a constant \(N_0\ge1\) such that
	\begin{equation}\label{1.6.1}
		\int_{B_{N_0}}\rho_0\,dx\ge\frac12 M_0>0.
	\end{equation}
	Set
	\[
	E_0:=
	\int_{\mathbb R^2}
	\left(
	\frac12\rho_0|u_0|^2+\frac{1}{\gamma-1}\rho_0^\gamma
	\right)\,dx,
	\]
	\[
	\Xi_0
	:=
	\bar\rho^\gamma M_0+(\bar\rho^\gamma M_0)^{-1}
	+\bar\rho^{\gamma+\frac{\alpha(\gamma+1)}2}
	\|\rho_0|x|^\alpha\|_{L^1}.
	\]
	We define the two-dimensional scaling-invariant initial quantity by
	\begin{equation}\label{definition-I0}
		\begin{aligned}
			I_0
			:=\,
			\Xi_0(1+\bar\rho E_0)^{3/2}
			\bar\rho^{1-\gamma}
			\left(
			\|\nabla u_0\|_{L^2}^2
			+
			\|\rho_0^\gamma\|_{L^2}^2
			\right).
		\end{aligned}
	\end{equation}

	The main results of this paper read as follows:
	
	\begin{theorem}\label{th1.1}
		Let $\gamma>1$, $q>2$, $\alpha>1$, $\bar{\rho}\geq 1 $  be fixed constants.
		Assume that the initial data satisfy
		\begin{equation}\label{assumption}
			\begin{gathered}
				0 \leq \rho_0 \leq \bar{\rho}, \quad
				\bar{x}^\alpha \rho_0 \in L^1(\mathbb{R}^2) \cap H^1(\mathbb{R}^2) \cap W^{1,q}(\mathbb{R}^2), \\
				u_0 \in D^1(\mathbb{R}^2) , \quad \sqrt{\rho_0} u_0 \in L^2(\mathbb{R}^2),
			\end{gathered}
		\end{equation}
		where $\bar{x} :=\left(e+|x|^2\right)^{1 / 2} \log ^2\left(e+|x|^2\right)$.
		In addition, there exists a positive constant $\varepsilon_0$ depending only on $\mu, \lambda, \gamma, q, \alpha$, which is independent of the initial data, time, and
		\(\bar\rho\)
		such that if
		\begin{equation}\label{initial-scale-smallness}
			I_0\le \varepsilon_0,
		\end{equation}
		then the above Cauchy problem admits a unique global strong solution
		$(\rho,u)$ on $\mathbb R^2\times(0,\infty)$ satisfying, for every $T>0$,
		\begin{equation}\label{1.8.1}
			0\le \rho(x,t)\le 2\bar\rho
			\qquad \forall \ (x,t)\in \mathbb{R}^2\times[0,T],
		\end{equation}
		\begin{equation}\label{1.9.1}
			\left\{\begin{array}{l}
				\rho \in C\left([0, T] ; L^1 \cap H^1 \cap W^{1, q}\right) \\
				\bar{x}^\alpha  \rho \in L^{\infty}\left(0, T ; L^1 \cap H^1 \cap W^{1, q}\right) \\
				\sqrt{\rho} u, \nabla u, \bar{x}^{-1} u, \sqrt{t} \sqrt{\rho} u_t \in L^{\infty}\left(0, T ; L^2\right) \\
				\nabla u \in L^2\left(0, T ; H^1\right) \cap L^{(q+1) / q}\left(0, T ; W^{1, q}\right) \\
				\sqrt{t} \nabla u \in L^2\left(0, T ; W^{1, q}\right) \\
				\sqrt{\rho} u_t, \sqrt{t} \nabla u_t, \sqrt{t} \bar{x}^{-1} u_t \in L^2\left(\mathbb{R}^2 \times(0, T)\right)
			\end{array}\right.
		\end{equation}
		Moreover, there exists a constant $N_1>0$, independent of $T$, such that
		\begin{equation}\label{1.10.1}
			\inf_{0\le t\le T}
			\int_{B_{N_1(1+t)}}\rho(x,t)\,dx
			\ge \frac{1}{4} M_0.
		\end{equation}
	\end{theorem}
	
	\begin{remark}\label{rem:wen-wu}
		 Upon completing the present work, the authors discussed this topic privately with Professor Huanyao Wen in a friendly exchange. It turns out that, based on another interesting  scaling-invariant smallness condition as follows
		\[
		\hat\rho E_0+\hat\rho E_0\bigl(1+(\hat\rho E_0)^3\bigr)
		(\hat\rho^\gamma M_0)^{-\frac{6}{5}}R_*^{\frac{4}{5}}
		\Bigl[(\hat\rho E_0)^2
		+\bigl[R_*^4+(\hat\rho^\gamma M_0)^2\bigr]
		(\hat\rho^{1-\gamma}H_0)^2\Bigr]\le\varepsilon_* ,
		\]
		where
		\(H_0=\|\nabla u_0\|_2^2+\|\rho_0^\gamma\|_2^2\) and
		\[
		R_*:=1+\frac{\hat\rho^{\frac{\gamma+1}{2}}}{M_0}
		\int_{\mathbb R^2}\rho_0|x|\,dx
		+\sqrt{\frac{2\hat\rho^{1-\gamma}E_0}{M_0}},
		\]
		Wen-Wu \cite{WenWu2026} adopted distinct analytical approaches and independently  established the global existence and uniqueness of strong solutions to system \eqref{1.1}. It should be noted that neither of these two conditions implies the other.
	\end{remark}

	\begin{remark}
		By Lemma \ref{lem6.1} in the appendix, the quantity
		\[
		I_0
		=
		\Xi_0(1+\bar\rho E_0)^{3/2}
		\bar\rho^{1-\gamma}
		\left(
		\|\nabla u_0\|_{L^2}^2+\|\rho_0^\gamma\|_{L^2}^2
		\right)
		\]
		is invariant under the natural scaling of the two-dimensional isentropic
		compressible Navier--Stokes equations. This is the reason why the smallness
		condition in Theorem \ref{th1.1} is independent of the particular choice of the
		density upper bound \(\bar\rho\).
	\end{remark}
	
	\begin{remark}
		Lemmas \ref{lem6.2}--\ref{lem6.3} further clarify why the particular
		scaling-invariant smallness condition in Theorem \ref{th1.1} is reasonable. If
		the power of \(\Xi_0\) is chosen too large, Lemma \ref{lem6.2} shows that the
		corresponding scaling-invariant quantity admits a positive lower bound
		independent of the initial data and \(\bar\rho\), so that no arbitrarily small
		universal threshold can be imposed. In contrast, Lemma \ref{lem6.3} shows that
		the present choice \(I_0\), with the first power \(\Xi_0\), does not have such
		an obstruction and can be made arbitrarily small along a suitable sequence of
		initial data.
	\end{remark}
	
	\begin{remark}
		For the isentropic compressible Navier--Stokes equations in \(\mathbb R^N\)
		\((N\ge 3)\), the natural scaling leads to the invariant quantity
		\[
		I_0^{(N)}
		:=
		\bar\rho^{\,N+(N-3)\gamma}
		E_0(\|\nabla u_0\|_{L^2}^2+\|\rho_0^\gamma\|_{L^2}^2)
		\Bigl(1+\bar\rho^{\,N+(N-2)\gamma}E_0^2\Bigr).
		\]
		In particular, when \(N=3\), this becomes
		\[
		I_0^{(3)}
		=
		\bar\rho^3
		E_0(\|\nabla u_0\|_{L^2}^2+\|\rho_0^\gamma\|_{L^2}^2)
		\Bigl(1+\bar\rho^{3+\gamma}E_0^2\Bigr).
		\]
		This is the scaling-invariant quantity used by Wen \cite{Wen2025} to extend the
		local solution to a global-in-time solution in the three-dimensional case.
	\end{remark}

	\begin{remark}
		Li--Xin \cite{Li-xin-2019} introduced the small-energy condition \(E_0\le\varepsilon_*\). The present work uses the scaling-invariant smallness condition \(I_0\le\varepsilon_0\). Lemma \ref{lem6.4} shows that these conditions are not equivalent. More precisely, \(E_0\le\varepsilon_*\) does not imply \(I_0\le\varepsilon_0\). Conversely, \(I_0\le\varepsilon_0\) does not imply \(E_0\le\varepsilon_*\).
	\end{remark}
	
	\begin{theorem}\label{th1.2}
		In addition to the assumptions in Theorem \ref{th1.1}, assume further that $(\rho_0, u_0)$ satisfy
		\begin{equation}\label{1.12.1}
			\nabla^2 \rho_0,\ \nabla^2 P(\rho_0) \in L^2 \cap L^q,\quad
			\bar{x}^{\delta_0} \nabla^2 \rho_0,\ \bar{x}^{\delta_0} \nabla^2 P(\rho_0),\ \nabla^2 u_0 \in L^2,
		\end{equation}
		for some constant $\delta_0 \in (0,1)$, and the following compatibility condition:
		\begin{equation}\label{1.13.1}
			-\mu \Delta u_0 - (\mu+\lambda) \nabla \mathrm{div}\, u_0 + \nabla P(\rho_0) = \rho_0^{1/2} g,
		\end{equation}
		with some $g \in L^2$. Then, in addition to \eqref{1.8.1}-\eqref{1.10.1}, the strong solution $(\rho, u)$ obtained by Theorem \ref{th1.1} satisfies for any $0 < T < \infty$,
		\begin{equation}\label{1.14.1}
			\begin{cases}
				\nabla^2 \rho,\ \nabla^2 P(\rho) \in C([0,T]; L^2 \cap L^q), \\
				\bar{x}^{\delta_0} \nabla^2 \rho,\ \bar{x}^{\delta_0} \nabla^2 P(\rho),\ \nabla^2 u \in L^\infty(0,T; L^2), \\
				\sqrt{\rho} u_t,\ \sqrt{t} \nabla u_t,\ \sqrt{t} \bar{x}^{-1} u_t,\ t \sqrt{\rho} u_{tt},\ t \nabla^2 u_t \in L^\infty(0,T; L^2), \\
				t \nabla^3 u \in L^\infty(0,T; L^2 \cap L^q), \\
				\nabla u_t,\ \bar{x}^{-1} u_t,\ t \nabla u_{tt},\ t \bar{x}^{-1} u_{tt} \in L^2(0,T; L^2), \\
				t \nabla^2 (\rho u) \in L^\infty(0,T; L^{(q+2)/2}).
			\end{cases}
		\end{equation}
	\end{theorem}
	
	\begin{remark}
		The solution obtained in Theorem \ref{th1.2} becomes a classical one for positive
		time. Although the smallness assumption is imposed on the scaling-invariant
		quantity \(I_0\), the initial oscillations of the density and velocity may still
		be large. In particular, both interior vacuum and far-field vacuum are allowed,
		and there is no requirement on the size of the vacuum set. Therefore, compactly
		supported initial density is included in our framework, provided the weighted
		regularity and compatibility assumptions in Theorems \ref{th1.1}--\ref{th1.2}
		are satisfied.
	\end{remark}
	
	\begin{remark}
		It is worth noting that Theorems \ref{th1.1}--\ref{th1.2} are closely related
		to, but different from, the two-dimensional theory of Li--Xin
		\cite{Li-xin-2019}. In \cite{Li-xin-2019}, the smallness is imposed on the
		initial energy after the density upper bound and the mass distribution parameters
		are fixed. Here the smallness is imposed on the scaling-invariant quantity
		\(I_0\). The factor \(\Xi_0\) records, in scaling-invariant form, the total mass,
		the reciprocal of the scaled mass, and the weighted moment of the initial
		density. These quantities are needed in the far-field vacuum case because there
		is no positive lower bound of the density at infinity, and the weighted kinetic
		energy term alone does not give an unweighted \(L^p\)-control of the velocity.
		Another difference from \cite{Li-xin-2019} is that the lower-order global
		argument is not organized through the homogeneous Sobolev framework used there.
		Instead, the weighted density information is combined with the normalization of
		\(\bar\rho\) and the scaling-invariant quantity \(I_0\). The key technical step
		behind this formulation is the normalization of the density upper bound. We use
		the exact scaling of the system and take
		\(\ell=\bar\rho^{-\gamma}\), so that
		\[
		\widetilde\rho(y,\tau)=\bar\rho^{-1}
		\rho\bigl(\bar\rho^{-\frac{\gamma+1}{2}}y,\bar\rho^{-\gamma}\tau\bigr),
		\qquad
		\widetilde u(y,\tau)=\bar\rho^{-\frac{\gamma-1}{2}}
		u\bigl(\bar\rho^{-\frac{\gamma+1}{2}}y,\bar\rho^{-\gamma}\tau\bigr).
		\]
		Then \(0\le\widetilde\rho_0\le1\),
		\(\widetilde M_0=\bar\rho^\gamma M_0\),
		\(\widetilde E_0=\bar\rho E_0\), and
		\(\Xi_0=\widetilde M_0+\widetilde M_0^{-1}
		+\|\widetilde\rho_0|y|^\alpha\|_{L^1}.\)
		Moreover,
		\[
		\|\nabla\widetilde u_0\|_{L^2}^2+\|\widetilde P_0\|_{L^2}^2
		=
		\bar\rho^{1-\gamma}
		(\|\nabla u_0\|_{L^2}^2+\|P_0\|_{L^2}^2).
		\]
		Hence \(I_0\) is exactly the smallness quantity obtained after the density
		upper bound has been normalized to one. All lower-order estimates in Section 3
		are proved first in these normalized variables and then scaled back, which is
		what makes the final threshold independent of \(\bar\rho\).
	\end{remark}
	
	We now make some comments on the analysis of this paper. By the local existence
	theory \cite{Liang2014 2D}, it is enough to establish global a priori estimates
	for smooth solutions. The central point is Proposition \ref{p3.1}, which is a
	bootstrap statement improving simultaneously the density upper bound
	\(\rho\le2\bar\rho\) to \(\rho\le\frac32\bar\rho\) and the invariant
	lower-order bound to a strictly smaller one. Once this proposition is proved, the
	standard continuation argument gives the global strong solution; the additional
	regularity in Theorem \ref{th1.2} is then obtained from the higher order
	estimates and the compatibility condition.
	
	The proof of Proposition \ref{p3.1} contains two linked mechanisms. The first one
	is the scale-invariant lower-order energy closure. After the above normalization,
	we may work with \(\bar\rho=1\). Multiplying the momentum equation by the
	material derivative \(\dot u\), using the pressure equation
	\(P_t+u\cdot\nabla P+\gamma P\operatorname{div}u=0\), and exploiting the
	effective viscous flux
	\(F=(2\mu+\lambda)\operatorname{div}u-P\), we obtain the short-time estimates for
	\(\nabla u\), \(P\), \(\sqrt\rho\dot u\), and \(P\) in \(L^3\). The pressure
	terms cannot simply be treated as lower-order errors; they are closed by a
	complementary pressure estimate, which is why the quantity
	\(\|\nabla u\|_{L^2}^2+\|P\|_{L^2}^2\) appears as a single block in \(I_0\).
	The time-weighted estimate with
	\(\sigma(t)=\min\{1,\bar\rho^\gamma t\}\) then gives control of
	\(\sigma\nabla\dot u\). These estimates form the invariant analogue of the basic
	energy bounds in \cite{Li-xin-2019}, but here each term has to be compatible with
	the scaling identities for \(\bar\rho\), \(M_0\), \(E_0\), and the weighted
	moment.
	
	The second mechanism is the uniform upper bound of the density. This is the part
	where the two-dimensional nature of the problem is most visible. The mass
	expansion estimate gives a positive amount of mass in balls whose radius grows at
	most linearly in time. Together with the weighted moment bound, it allows us to
	use the weighted Sobolev inequality to control the velocity, despite the lack of
	a direct energy bound for \(u\). The density equation is then rewritten along
	particle paths in terms of the effective viscous flux and treated by Zlotnik's
	inequality. For short time, the \(L^1_tL^\infty_x\)-norm of \(F\) follows from
	the weighted Sobolev inequality and the material derivative estimate. For large
	time, one needs the decay of \(P\), \(\nabla u\), \(\sqrt\rho\dot u\), and
	\(\nabla\dot u\). In particular, Lemma \ref{lem:large-time-L2-decay} yields \(\|\nabla u(t)\|_{L^2}+\|P(t)\|_{L^2}+\|\sqrt\rho\dot u(t)\|_{L^2}\le C t^{-1/2}\) for \(t\ge\bar\rho^{-\gamma}\), where \(C\) depends on the initial data and fixed parameters but not on time; \(\nabla\dot u\) is controlled by the corresponding time-weighted integral estimate. This large-time decay argument follows the core idea of
	Li--Xin \cite{Li-xin-2019}, but it is carried out after the
	\(\bar\rho\)-normalization and is then scaled back to the original variables.
	This is what makes the density estimate close under the invariant smallness of
	\(I_0\). After the uniform density bound is obtained, the invariant lower-order
	bound is improved by the closed energy estimates, and Proposition \ref{p3.1}
	follows. Finally, the bounds for the gradients of the density and velocity are
	derived through the effective flux-vorticity decomposition, a Beale--Kato--Majda
	type logarithmic inequality, and a logarithmic Gronwall argument. The higher order
	estimates needed for the classical solution are then proved in the same spirit as
	\cite{Li-xin-2019} and the local theory \cite{Liang2014 2D}. The appendix is
	used to verify the scaling invariance of \(I_0\) and to justify the choice of the
	power of \(\Xi_0\) in the smallness condition, in particular by showing the
	obstruction caused by higher powers and the absence of such an obstruction for
	the present \(I_0\).

	The paper is organized as follows. In Section 2, we collect several fundamental facts and inequalities that are essential for the later analysis. Sections 3 and 4 are devoted to establishing the required a priori estimates on classical solutions, which are used to extend the local solution to the global-in-time solution. The main results including Theorems \ref{th1.1}--\ref{th1.2} are proved in Section 5. Finally, in Section 6, we verify the scaling invariance of \(I_0\) and discuss the reasonableness of the power of \(\Xi_0\) appearing in the scaling-invariant smallness condition.

	\section{Preliminaries}
	\quad \quad In this section, we review several fundamental results and classical inequalities that will play a key role in the subsequent analysis.
	
	We start by recalling the local existence theory for strong and classical solutions, whose detailed proof can be found in the literature \cite{Liang2014 2D}.
	
	\begin{lemma}\label{lem2.1}
		Suppose that the initial data $(\rho_0, u_0)$ satisfy:
		\begin{equation}\label{2.1.1}
			\rho_0 \geq 0,\quad
			\bar{x}^\alpha \rho_0 \in L^1 \cap H^1 \cap W^{1,q},\quad
			u_0 \in D^1,\quad
			\rho_0^{1/2} u_0 \in L^2.
		\end{equation}
		Then there exists a sufficiently small time $T > 0$ such that problem \eqref{1.1}--\eqref{1.4.1} on $\mathbb{R}^2 \times (0, T)$ admits a unique strong solution $(\rho, u)$, which also satisfies estimates \eqref{1.9.1} and \eqref{1.10.1}. Furthermore, if $(\rho_0, u_0)$ additionally meets conditions \eqref{1.12.1} and \eqref{1.13.1}, then the solution $(\rho, u)$ also verifies the regularity results \eqref{1.14.1}.
	\end{lemma}
	
	Next, we introduce the well-known Gagliardo--Nirenberg inequality (see, e.g., \cite{Nirenberg1959 PDE}), which will be employed repeatedly in what follows.
	
	\begin{lemma}[Gagliardo--Nirenberg]\label{lem2.2}
		For any $p \in [2, \infty)$, $q \in (1, \infty)$, and $r \in (2, \infty)$, there exists a  constant $C > 0$ depending on $p, q, r$ such that for all $f \in H^1(\mathbb{R}^2)$ and $g \in L^q(\mathbb{R}^2) \cap D^{1, r}(\mathbb{R}^2)$, the following inequalities hold:
		\begin{align}
			\|f\|_{L^p(\mathbb{R}^2)}^p &\leq C \|f\|_{L^2(\mathbb{R}^2)}^2 \|\nabla f\|_{L^2(\mathbb{R}^2)}^{p-2}, \label{2.2.1} \\
			\|g\|_{C(\overline{\mathbb{R}^2})} &\leq C \|g\|_{L^q(\mathbb{R}^2)}^{q(r-2)/(2r+q(r-2))} \|\nabla g\|_{L^r(\mathbb{R}^2)}^{2r/(2r+q(r-2))}. \label{2.3.1}
		\end{align}
	\end{lemma}
	
	We then present the weighted $L^p$ estimates for functions in the Hilbert space $D^1(\mathbb{R}^2)$, which are established in \cite[Theorem B.1]{Lions1996}.
	
	\begin{lemma}\label{lem2.3}
		For any $m \in [2, \infty)$ and $\theta \in (1 + m/2, \infty)$, there exists a positive constant $C$ such that for all $v \in D^1(\mathbb{R}^2)$,
		\begin{equation}\label{2.4.1}
			\left( \int_{\mathbb{R}^2} \frac{|v|^m}{e + |x|^2} \left( \log(e + |x|^2) \right)^{-\theta} dx \right)^{1/m}
			\leq C \|v\|_{L^2(B_1)} + C \|\nabla v\|_{L^2(\mathbb{R}^2)}.
		\end{equation}
	\end{lemma}
	
	Combining Lemma \ref{lem2.3} with the classical Poincar\'e inequality leads to further useful estimates.
	
	\begin{lemma}\label{lem2.4.1}
		Assume that \(\rho \in L^\infty(\mathbb R^2)\) satisfies
		\begin{equation}\label{2.5}
			0\le \rho(x)\le m_1,\qquad
			\int_{B_{N_*}}\rho(x)\,dx\ge m_2,\qquad
			\rho\bar x^a\in L^1(\mathbb R^2),
		\end{equation}
		where \(m_1,N_*\ge1\), and \(m_2,a>0\). Then, for any
		\(r\in[2,\infty)\) and
		\[
		v\in
		\left\{
		v\in D^1(\mathbb R^2)
		\ \big|\
		\rho^{1/2}v\in L^2(\mathbb R^2)
		\right\},
		\]
		there exists a constant \(C=C(a,r)\) such that
		\begin{equation}\label{2.6}
			\begin{aligned}
				\left(\int_{\mathbb R^2}\rho |v|^r\,dx\right)^{1/r}
				\le{}&
				C(a,r)
				m_1^{\frac{a}{r(4+a)}}
				\left(1+\|\rho\bar x^a\|_{L^1}\right)^{\frac{4}{r(4+a)}}\\
				&\times
				N_*
				\left(1+\sqrt{\frac{m_1}{m_2}}\right)
				\left(
				\|\rho^{1/2}v\|_{L^2}
				+
				\|\nabla v\|_{L^2}
				\right).
			\end{aligned}
		\end{equation}
	\end{lemma}
	
	\begin{proof}
		For \(f\in L^1(B_{N_*})\), denote
		\[
		f_{B_{N_*}}
		:=
		\frac1{|B_{N_*}|}
		\int_{B_{N_*}}f(x)\,dx.
		\]
		We first estimate the average \(v_{B_{N_*}}\). Since
		\[
		v_{B_{N_*}}
		=
		\frac{1}{\int_{B_{N_*}}\rho\,dx}
		\int_{B_{N_*}}\rho v\,dx
		+
		\frac{1}{\int_{B_{N_*}}\rho\,dx}
		\int_{B_{N_*}}\rho(v_{B_{N_*}}-v)\,dx,
		\]
		we have
		\[
		\begin{aligned}
			|v_{B_{N_*}}|
			&\le
			\frac{1}{\int_{B_{N_*}}\rho\,dx}
			\left|\int_{B_{N_*}}\rho v\,dx\right|
			+
			\frac{1}{\int_{B_{N_*}}\rho\,dx}
			\int_{B_{N_*}}\rho|v-v_{B_{N_*}}|\,dx\\
			&\le
			\frac{1}{\left(\int_{B_{N_*}}\rho\,dx\right)^{1/2}}
			\|\rho^{1/2}v\|_{L^2(B_{N_*})}+
			\frac{1}{\left(\int_{B_{N_*}}\rho\,dx\right)^{1/2}}
			\|\rho^{1/2}(v-v_{B_{N_*}})\|_{L^2(B_{N_*})}.
		\end{aligned}
		\]
		Using \(0\le\rho\le m_1\) and \(\int_{B_{N_*}}\rho\,dx\ge m_2\), we get
		\[
		\|\rho^{1/2}(v-v_{B_{N_*}})\|_{L^2(B_{N_*})}
		\le
		m_1^{1/2}
		\|v-v_{B_{N_*}}\|_{L^2(B_{N_*})}.
		\]
		By the Poincar\'e inequality \cite[(7.45)]{GilbargTrudinger1983} on \(B_{N_*}\), we have
		\[
		\|v-v_{B_{N_*}}\|_{L^2(B_{N_*})}
		\le
		4N_*\|\nabla v\|_{L^2(B_{N_*})},
		\]
		which implies that
		\begin{equation}\label{2.7-new}
			|v_{B_{N_*}}|
			\le
			\frac1{\sqrt{m_2}}\|\rho^{1/2}v\|_{L^2}
			+
			CN_*\sqrt{\frac{m_1}{m_2}}\|\nabla v\|_{L^2}.
		\end{equation}
		
		We next estimate the local \(L^2\)-norm on the fixed ball \(B_1\), not on
		\(B_{N_*}\). Since \(N_*\ge1\), one has \(B_1\subset B_{N_*}\). Hence
		\[
		\begin{aligned}
			\|v\|_{L^2(B_1)}
			&\le
			\|v-v_{B_{N_*}}\|_{L^2(B_1)}
			+
			|B_1|^{1/2}|v_{B_{N_*}}|\\
			&\le
			\|v-v_{B_{N_*}}\|_{L^2(B_{N_*})}
			+
			C|v_{B_{N_*}}|\\
			&\le
			CN_*\|\nabla v\|_{L^2}
			+
			C|v_{B_{N_*}}|.
		\end{aligned}
		\]
		Using \eqref{2.7-new}, we obtain
		\begin{equation}\label{local-L2-B1}
			\|v\|_{L^2(B_1)}
			\le
			C\frac1{\sqrt{m_2}}\|\rho^{1/2}v\|_{L^2}
			+
			CN_*
			\left(1+\sqrt{\frac{m_1}{m_2}}\right)
			\|\nabla v\|_{L^2}.
		\end{equation}
		
		Set
		\[
		\sigma:=\frac{4}{4+a}\in(0,1).
		\]
		By H\"older's inequality, we obtain
		\[
		\begin{aligned}
			\int_{\mathbb R^2}\rho|v|^r\,dx
			&=
			\int_{\mathbb R^2}
			(\rho\bar x^a)^\sigma
			\left(|v|^r\bar x^{-a\sigma}\right)
			\rho^{1-\sigma}\,dx\\
			&\le
			\|\rho\bar x^a\|_{L^1}^{\sigma}
			\left\|v\bar x^{-\frac{a\sigma}{r}}\right\|_{L^{\frac{4r}{a\sigma}}}^{r}
			\|\rho\|_{L^\infty}^{1-\sigma}\\
			&\le
			\|\rho\bar x^a\|_{L^1}^{\sigma}
			m_1^{1-\sigma}
			\left\|v\bar x^{-\frac{a\sigma}{r}}\right\|_{L^{\frac{4r}{a\sigma}}}^{r}.
		\end{aligned}
		\]
		Thus
		\begin{equation}\label{rho-v-r-before-weighted-sobolev}
			\left(\int_{\mathbb R^2}\rho|v|^r\,dx\right)^{1/r}
			\le
			m_1^{\frac{1-\sigma}{r}}
			\|\rho\bar x^a\|_{L^1}^{\frac{\sigma}{r}}
			\left\|v\bar x^{-\frac{a\sigma}{r}}\right\|_{L^{\frac{4r}{a\sigma}}}.
		\end{equation}
		
		By the weighted Sobolev inequality \eqref{2.4.1} applied with the fixed ball \(B_1\), we have
		\[
		\left\|v\bar x^{-\frac{a\sigma}{r}}\right\|_{L^{\frac{4r}{a\sigma}}}
		\le
		C(a,r)
		\left(
		\|v\|_{L^2(B_1)}
		+
		\|\nabla v\|_{L^2(\mathbb R^2)}
		\right).
		\]
		Substituting \eqref{local-L2-B1} into the previous estimate, we get
		\[
		\begin{aligned}
			\left\|v\bar x^{-\frac{a\sigma}{r}}\right\|_{L^{\frac{4r}{a\sigma}}}
			\le{}&
			C(a,r)
			\frac1{\sqrt{m_2}}\|\rho^{1/2}v\|_{L^2}+
			C(a,r)
			N_*
			\left(1+\sqrt{\frac{m_1}{m_2}}\right)
			\|\nabla v\|_{L^2}.
		\end{aligned}
		\]
		Since \(N_*\ge1\) and \(m_1\ge1\), one has
		\[
		\frac1{\sqrt{m_2}}
		\le
		N_*
		\left(1+\sqrt{\frac{m_1}{m_2}}\right).
		\]
		Therefore, we have
		\[
		\left\|v\bar x^{-\frac{a\sigma}{r}}\right\|_{L^{\frac{4r}{a\sigma}}}
		\le
		C(a,r)
		N_*
		\left(1+\sqrt{\frac{m_1}{m_2}}\right)
		\left(
		\|\rho^{1/2}v\|_{L^2}
		+
		\|\nabla v\|_{L^2}
		\right).
		\]
		Combining this with \eqref{rho-v-r-before-weighted-sobolev}, and using
		\[
		\frac{1-\sigma}{r}
		=
		\frac{a}{r(4+a)},
		\qquad
		\frac{\sigma}{r}
		=
		\frac{4}{r(4+a)},
		\]
		we obtain
		\[
		\begin{aligned}
			\left(\int_{\mathbb R^2}\rho|v|^r\,dx\right)^{1/r}
			\le{}&
			C(a,r)
			m_1^{\frac{a}{r(4+a)}}
			\|\rho\bar x^a\|_{L^1}^{\frac{4}{r(4+a)}}\times
			N_*
			\left(1+\sqrt{\frac{m_1}{m_2}}\right)
			\left(
			\|\rho^{1/2}v\|_{L^2}
			+
			\|\nabla v\|_{L^2}
			\right),
		\end{aligned}
		\]
		which gives \eqref{2.6}. The proof is complete.
	\end{proof}

	We now turn to establishing a uniform-in-time upper bound for the density $\rho$, for which we require the classical Zlotnik inequality (see, e.g., \cite{Zlotnik2000}).
	
	\begin{lemma}
		Let the function $y$ satisfy the differential equation
		\begin{equation*}
			y'(t) = g(y(t)) + b'(t) \quad \text{on } [0, T], \quad y(0) = y^0,
		\end{equation*}
		where $g \in C(\mathbb{R})$ and $y, b \in W^{1,1}(0, T)$. Assume further that $g(\infty) = -\infty$ and that $b$ satisfies
		\begin{equation}
			b(t_2) - b(t_1) \leq N_0 + N_1(t_2 - t_1)
		\end{equation}
		for all $0 \leq t_1 < t_2 < T$, with constants $N_0, N_1 \geq 0$. If there exists a constant $\bar{\zeta}$ such that
		\begin{equation}
			g(\zeta) \leq -N_1 \quad \text{for all } \zeta \geq \bar{\zeta},
		\end{equation}
		then the solution $y$ admits the uniform upper bound
		\begin{equation*}
			y(t) \leq \max\left\{y^0, \bar{\zeta}\right\} + N_0 < \infty \quad \text{on } [0, T].
		\end{equation*}
	\end{lemma}
	
	Subsequently, we will also employ a Beale--Kato--Majda type inequality, originally established in \cite{Beale-Kato-Majda-1984,Kato1986Remarks} for the case $\mathrm{div}u \equiv 0$. This estimate will play a key role in controlling the norms $\|\nabla u\|_{L^\infty}$ and $\|\nabla \rho\|_{L^2 \cap L^q}$ for $q > 2$.
	\begin{lemma}
		For any exponent $q \in (2, \infty)$, there exists a constant $C = C(q) > 0$ such that for all vector fields $u$ with $\nabla u \in L^2(\mathbb{R}^2) \cap D^{1,q}(\mathbb{R}^2)$, the following estimate holds:
		\begin{equation}\label{2.9.1}
			\begin{split}
				\|\nabla u\|_{L^\infty(\mathbb{R}^2)}
				&\leq C\left( \|\mathrm{div}\,u\|_{L^\infty(\mathbb{R}^2)} + \|\omega\|_{L^\infty(\mathbb{R}^2)} \right)
				\log\left(e + \|\nabla^2 u\|_{L^q(\mathbb{R}^2)} \right) \\
				& + C\|\nabla u\|_{L^2(\mathbb{R}^2)} + C.
			\end{split}
		\end{equation}
	\end{lemma}

	Finally we introduce the notation $\nabla^\perp := (-\partial_2, \partial_1)$ and denote the material derivative of a function $f$ by $\dot{f} := f_t + u \cdot \nabla f$.
	From the momentum equations in \eqref{1.1}, we derive the following elliptic system:
	\begin{equation}
		\Delta F = \mathrm{div}(\rho \dot{u}), \quad \mu \Delta \omega = \nabla^\perp \cdot (\rho \dot{u}),
		\label{eq:2.9}
	\end{equation}
	where
	\begin{equation}\label{2.11.1}
		\omega := \partial_1 u^2 - \partial_2 u^1, \quad
		F := (2\mu + \lambda)\mathrm{div}\,u - P.
	\end{equation}
	This system will be used to establish the following elementary estimates, which rely on \eqref{2.1.1} and standard $L^p$-theory for elliptic equations.
	
	\begin{lemma}[\cite{Li-xin-2019}]\label{lem2.9}
		Let $(\rho, u)$ be a smooth solution to the compressible Navier--Stokes system \eqref{1.1}. Then for any $p \geq 2$, there exists a positive constant $C = C(p, \mu, \lambda)$ such that
		\begin{align}
			&\|\nabla F\|_{L^p(\mathbb{R}^2)} + \|\nabla \omega\|_{L^p(\mathbb{R}^2)} \leq C \|\rho \dot{u}\|_{L^p(\mathbb{R}^2)}, \label{2.12.1} \\
			&\|F\|_{L^p(\mathbb{R}^2)} + \|\omega\|_{L^p(\mathbb{R}^2)} \leq C \|\rho \dot{u}\|_{L^2(\mathbb{R}^2)}^{1-2/p} \left( \|\nabla u\|_{L^2(\mathbb{R}^2)} + \|P\|_{L^2(\mathbb{R}^2)} \right)^{2/p}, \label{2.13.1} \\
			&\|\nabla u\|_{L^p(\mathbb{R}^2)} \leq C \|\rho \dot{u}\|_{L^2(\mathbb{R}^2)}^{1-2/p} \left( \|\nabla u\|_{L^2(\mathbb{R}^2)} + \|P\|_{L^2(\mathbb{R}^2)} \right)^{2/p} + C \|P\|_{L^p(\mathbb{R}^2)}. \label{2.14.1}
		\end{align}
	\end{lemma}
	
	\begin{remark}
		The restriction \(p\geq 2\) in \eqref{2.12.1} can be relaxed. Indeed, by \eqref{eq:2.9}, both \(\nabla F\) and \(\nabla \omega\) are obtained from \(\rho \dot{u}\) through Calderon--Zygmund singular integral operators. Hence the standard \(L^p\)-boundedness of singular integrals (see, for instance, \cite{CalderonZygmund1952}) yields
		\[
		\|\nabla F\|_{L^p(\mathbb{R}^2)}+\|\nabla \omega\|_{L^p(\mathbb{R}^2)}
		\leq C\|\rho \dot{u}\|_{L^p(\mathbb{R}^2)},\quad 1<p<\infty.
		\]
	\end{remark}

	\section{Lower Order A Priori Estimates}
	\quad \quad In this section, we focus on deriving the essential a priori estimates for smooth solutions to the Cauchy problem \eqref{1.1}--\eqref{1.4.1} in the whole space $\mathbb{R}^2$. The local existence of strong solutions is provided by Lemma \ref{lem2.1}. For the conditional a priori estimates, let $(\rho, u)$ be a smooth solution to \eqref{1.1}--\eqref{1.4.1} on $\mathbb{R}^2 \times(0, T]$ for some $T>0$, with smooth initial data $(\rho_0, u_0)$ verifying the conditions in \eqref{assumption}. Short-time intervals are restricted to the interval of existence; long-time estimates apply only when $T>\bar{\rho}^{-\gamma}$.
	
	We first record the scaling normalization that will be used throughout this
	section. Recall that the two-dimensional isentropic compressible
	Navier--Stokes system is invariant under the scaling
	\[
	\rho_\ell(x,t)=\ell^{1/\gamma}
	\rho\bigl(\ell^{\frac{\gamma+1}{2\gamma}}x,\ell t\bigr),
	\qquad
	u_\ell(x,t)=\ell^{\frac{\gamma-1}{2\gamma}}
	u\bigl(\ell^{\frac{\gamma+1}{2\gamma}}x,\ell t\bigr).
	\]
	Taking
	\[
	\ell=\bar\rho^{-\gamma},
	\]
	we define, for \(0\le \tau\le \bar\rho^\gamma T\),
	\[
	\widetilde\rho(y,\tau)
	:=
	\bar\rho^{-1}
	\rho\bigl(\bar\rho^{-\frac{\gamma+1}{2}}y,
	\bar\rho^{-\gamma}\tau\bigr), \quad
	\widetilde u(y,\tau)
	:=
	\bar\rho^{-\frac{\gamma-1}{2}}
	u\bigl(\bar\rho^{-\frac{\gamma+1}{2}}y,
	\bar\rho^{-\gamma}\tau\bigr).
	\]
	Then it follows from  Theorem 1.1 in \cite{lx19b} that \((\widetilde\rho,\widetilde u)\) solves the same system on
	\(\mathbb R^2\times [0,\bar\rho^\gamma T]\). Moreover, if
	\(0\le \rho\le 2\bar\rho\), then
	\(0\le \widetilde\rho\le 2\),
	and if \(0\le \rho_0\le \bar\rho\), then
	\(0\le \widetilde\rho_0\le 1\).
	
	Let
	\[
	x=\bar\rho^{-\frac{\gamma+1}{2}}y,
	\qquad
	t=\bar\rho^{-\gamma}\tau .
	\]
	Then
	the pressure, material derivative and effective viscous flux satisfy
	\[
	\widetilde P(y,\tau)
	:=
	\widetilde\rho^\gamma(y,\tau)
	=
	\bar\rho^{-\gamma}P(x,t),
	\]
	\[
	\dot{\widetilde u}(y,\tau)
	:=
	\widetilde u_\tau+\widetilde u\cdot\nabla_y\widetilde u
	=
	\bar\rho^{-\frac{3\gamma-1}{2}}\dot u(x,t),
	\]
	and
	\[
	\widetilde F(y,\tau)
	:=
	(2\mu+\lambda)\operatorname{div}_y\widetilde u-\widetilde P
	=
	\bar\rho^{-\gamma}F(x,t).
	\]
	Furthermore,
	\[
	\nabla_y\widetilde u(y,\tau)
	=
	\bar\rho^{-\gamma}\nabla_x u(x,t),
	\qquad
	\nabla_y\dot{\widetilde u}(y,\tau)
	=
	\bar\rho^{-2\gamma}\nabla_x\dot u(x,t).
	\]
	
	The initial mass and energy transform as
	\[
	\widetilde M_0
	:=
	\int_{\mathbb R^2}\widetilde\rho_0\,dy
	=
	\bar\rho^\gamma M_0,
	\qquad
	\widetilde E_0
	:=
	\int_{\mathbb R^2}
	\left(
	\frac12\widetilde\rho_0|\widetilde u_0|^2
	+\frac{1}{\gamma-1}\widetilde\rho_0^\gamma
	\right)dy
	=
	\bar\rho E_0 .
	\]
	For every \(m\ge 0\),
	\[
	\int_{\mathbb R^2}\widetilde\rho_0(y)|y|^m\,dy
	=
	\bar\rho^{\gamma+\frac{m(\gamma+1)}2}
	\int_{\mathbb R^2}\rho_0(x)|x|^m\,dx .
	\]
	In particular,
	\[
	\Xi_0
	:=
	\bar\rho^\gamma M_0
	+(\bar\rho^\gamma M_0)^{-1}
	+\bar\rho^{\gamma+\frac{\alpha(\gamma+1)}2}
	\|\rho_0|x|^\alpha\|_{L^1}
	=
	\widetilde M_0+\widetilde M_0^{-1}
	+\|\widetilde\rho_0|y|^\alpha\|_{L^1}.
	\]
	
	The basic \(L^2\)-quantities obey
	\begin{equation}\label{scale-large-L2}
		\|\nabla_y\widetilde u(\tau)\|_{L^2_y}^2
		+\|\widetilde P(\tau)\|_{L^2_y}^2
		=
		\bar\rho^{1-\gamma}
		\left(
		\|\nabla_x u(t)\|_{L^2_x}^2
		+\|P(t)\|_{L^2_x}^2
		\right),
	\end{equation}
	and
	\[
	\int_{\mathbb R^2}\widetilde\rho|\dot{\widetilde u}|^2\,dy
	=
	\bar\rho^{1-2\gamma}
	\int_{\mathbb R^2}\rho|\dot u|^2\,dx,
	\quad
	\int_{\mathbb R^2}|\nabla_y\dot{\widetilde u}|^2\,dy
	=
	\bar\rho^{1-3\gamma}
	\int_{\mathbb R^2}|\nabla_x\dot u|^2\,dx.
	\]
	In addition,
	\begin{equation}\label{scale-large-P3}
		\|\widetilde P(\tau)\|_{L^3_y}^3
		=
		\bar\rho^{1-2\gamma}
		\|P(t)\|_{L^3_x}^3.
	\end{equation}
	Consequently, for \(0\le t\le T\), it follows that
	\[
	\int_0^{\bar\rho^\gamma t}
	\int_{\mathbb R^2}\widetilde\rho|\dot{\widetilde u}|^2\,dy\,ds
	=
	\bar\rho^{1-\gamma}
	\int_0^t
	\int_{\mathbb R^2}\rho|\dot u|^2\,dx\,dr, \quad
	\int_0^{\bar\rho^\gamma t}
	\|\widetilde P(s)\|_{L^3_y}^3\,ds
	=
	\bar\rho^{1-\gamma}
	\int_0^t
	\|P(r)\|_{L^3_x}^3\,dr.
	\]
	
	We shall also use the time weight
	\[
	\widetilde\sigma(\tau):=\min\{1,\tau\},
	\qquad
	\sigma(t):=\min\{1,\bar\rho^\gamma t\}.
	\]
	Then
	\[
	\widetilde\sigma(\tau)=\sigma(t)
	\qquad
	\text{whenever } \tau=\bar\rho^\gamma t.
	\]
	In particular,
	\[
	\int_0^{\bar\rho^\gamma t}
	\widetilde\sigma(s)^2
	\int_{\mathbb R^2}|\nabla_y\dot{\widetilde u}|^2\,dy\,ds
	=
	\bar\rho^{1-2\gamma}
	\int_0^t
	\sigma(r)^2
	\int_{\mathbb R^2}|\nabla_x\dot u|^2\,dx\,dr.
	\]
	The preceding scaling identities allow us to reduce the estimates below to the
	normalized case
	$
	\bar\rho=1,
	$
	and then recover the corresponding estimates in the original variables by
	scaling back. We shall use this convention throughout the rest of this section
	without further comment.

	The following proposition plays a crucial role in deriving the global existence of strong and classical solutions to the two-dimensional problem \eqref{1.1}--\eqref{1.4.1}.
	\begin{proposition}\label{p3.1}
		Under the hypotheses of Theorem \ref{th1.1}, in addition assume that the initial-value problem has a strong solution $(\rho,u)$ on $[0,T]$, where $T$ is a positive constant. Then there exists a small positive constant $\varepsilon_1\in(0,1]$, independent of the initial data, $\bar\rho$, and $T$, such that if
		\begin{equation}\label{3.1.3}
			\left\{
			\begin{aligned}
				&\rho \le 2\bar\rho,\\
				&
				\Xi_0(1+\bar\rho E_0)^{3/2}
				\bar\rho^{1-\gamma}
				\sup_{0\le t\le T}
				\left(
				\|\nabla u(t)\|_{L^2}^2+\|P(t)\|_{L^2}^2
				\right)
				\le 2\varepsilon
			\end{aligned}
			\right.
		\end{equation}
		and $\varepsilon\in (c_0\varepsilon_0,\varepsilon_1]$ for a positive constant $c_0$, depending only on $\mu,\lambda,\gamma,\alpha$, then
		\begin{equation}\label{3.2}
			\left\{
			\begin{aligned}
				&\rho \le \frac{3}{2}\bar\rho,\\
				&
				\Xi_0(1+\bar\rho E_0)^{3/2}
				\bar\rho^{1-\gamma}
				\sup_{0\le t\le T}
				\left(
				\|\nabla u(t)\|_{L^2}^2+\|P(t)\|_{L^2}^2
				\right)
				\le \frac{3}{2}\varepsilon .
			\end{aligned}
			\right.
		\end{equation}
	\end{proposition}
	
	We first remark that the proof of Proposition \ref{p3.1} will be deferred to the end of this section.
	In what follows, we adopt the standard convention that $C$ stands for a generic positive constant
	which may depend on $\mu,\lambda,\gamma,\alpha$.
	Moreover, we write $C(\alpha)$ to explicitly indicate that the constant $C$ depends on~$\alpha$.
	
	Our analysis begins with the following standard energy identity for the system, which describes the basic energy balance.
	
	\begin{lemma}\label{lem3.1}
		Let $(\rho, u)$ solve \eqref{1.1}--\eqref{1.4.1} smoothly on $\mathbb{R}^2 \times (0, T]$. Then the following energy identity holds, for all $t \in [0, T]$,
		\begin{equation}\label{3.1.1}
			\int_{\mathbb{R}^2} \left( \frac{1}{2} \rho |u|^2 + \frac{1}{\gamma - 1} \rho^\gamma \right) dx
			+ \int_0^t \int_{\mathbb{R}^2} \left[ \mu |\nabla u|^2 + (\mu + \lambda) |\mathrm{div}\,u|^2 \right] dx ds
			=  E_0.
		\end{equation}
	\end{lemma}
	\begin{proof}
		This energy equality is a standard result, which follows directly from the smoothness of the solution.
		Specifically, multiplying the momentum equation \eqref{1.1}$_2$ by $u$, integrating by parts over $\mathbb{R}^2$
		and combining with the continuity equation \eqref{1.1}$_1$ we obtain the desired identity \eqref{3.1.1}.
	\end{proof}
	
	Next, we give a standard energy estimate for the system $(\rho, u)$,
	together with some preliminary $L^2$-estimate for $\nabla u$ and $\rho\dot u$.
	
	\begin{lemma}\label{lem3.2}
		Let \((\rho,u)\) solve \eqref{1.1}--\eqref{1.4.1} smoothly on
		\(\mathbb R^2\times(0,T]\) satisfying \eqref{3.1.3} and the assumptions of
		Theorem \ref{th1.1}.  Then, for any \(t\in[0,T]\), there holds
		\begin{equation}\label{eq:short-time-u}
			\begin{aligned}
				&\|\nabla u(t)\|_{L^2}^2
				+\int_0^t\int_{\mathbb R^2}\rho|\dot u|^2\,dxds \\
				\le{}&
				C_{1}\Bigl(\|\nabla u_0\|_{L^2}^2+\|P_0\|_{L^2}^2\Bigr)
				+C_{1}\|P(t)\|_{L^2}^2
				+C_{1}\int_0^t\|P(s)\|_{L^3}^3\,ds,
			\end{aligned}
		\end{equation}
		where \(\dot u=u_t+u\cdot\nabla u\), \(P_0=\rho_0^\gamma\), and \(C_{1}\) depends only on
		\(\mu,\lambda,\gamma\).
	\end{lemma}
	
	\begin{proof}
		We first normalize the upper bound of the density. Observe that, for any \(t\in[0,T]\) and \(\tau=\bar\rho^\gamma t\),
		\begin{equation}\label{scale-material-dissipation-exact}
			\begin{aligned}
				&\int_0^\tau\int_{\mathbb R^2}
				\widetilde\rho(y,\sigma)
				|\dot{\widetilde u}(y,\sigma)|^2\,dy\,d\sigma \\
				={}&
				\int_0^{\bar\rho^\gamma t}\int_{\mathbb R^2}
				\bar\rho^{-1}
				\rho\!\left(\bar\rho^{-\frac{\gamma+1}{2}}y,\bar\rho^{-\gamma}\sigma\right)
				\bar\rho^{-(3\gamma-1)}
				\left|
				\dot u\!\left(\bar\rho^{-\frac{\gamma+1}{2}}y,\bar\rho^{-\gamma}\sigma\right)
				\right|^2
				\,dy\,d\sigma \\
				={}&
				\bar\rho^{1-\gamma}
				\int_0^t\int_{\mathbb R^2}
				\rho(x,s)|\dot u(x,s)|^2\,dx\,ds.
			\end{aligned}
		\end{equation}
		Similarly,
		\begin{equation}\label{scale-pressure-dissipation-exact}
			\begin{aligned}
				\int_0^\tau\|\widetilde P(\sigma)\|_{L^3_y}^3\,d\sigma
				&=
				\int_0^{\bar\rho^\gamma t}
				\int_{\mathbb R^2}
				\bar\rho^{-3\gamma}
				P^3\!\left(\bar\rho^{-\frac{\gamma+1}{2}}y,\bar\rho^{-\gamma}\sigma\right)
				\,dy\,d\sigma \\
				&=
				\bar\rho^{1-\gamma}
				\int_0^t\|P(s)\|_{L^3_x}^3\,ds.
			\end{aligned}
		\end{equation}
		
		Thus it is enough to prove the following normalized estimate: for every
		\(\tau\in[0,\bar\rho^\gamma T]\),
		\begin{equation}\label{eq:short-time-u-normalized}
			\begin{aligned}
				&\|\nabla\widetilde u(\tau)\|_{L^2_y}^2
				+
				\int_0^\tau\int_{\mathbb R^2}
				\widetilde\rho|\dot{\widetilde u}|^2\,dy\,d\sigma \\
				\le{}&
				C\left(
				\|\nabla\widetilde u_0\|_{L^2_y}^2
				+\|\widetilde P_0\|_{L^2_y}^2
				\right)
				+
				C\|\widetilde P(\tau)\|_{L^2_y}^2
				+
				C\int_0^\tau\|\widetilde P(\sigma)\|_{L^3_y}^3\,d\sigma.
			\end{aligned}
		\end{equation}
		Indeed, once \eqref{eq:short-time-u-normalized} is proved, substituting
		\eqref{scale-material-dissipation-exact}, and \eqref{scale-pressure-dissipation-exact}
		into \eqref{eq:short-time-u-normalized}, and then multiplying both sides by
		\(\bar\rho^{\gamma-1}\), gives exactly \eqref{eq:short-time-u}.
		
		For simplicity, in the proof of \eqref{eq:short-time-u-normalized}, we drop the tildes
		and write \(t\) for the normalized time variable. Thus, in the following part of the
		proof, we work under
		\[
		\bar\rho=1,\qquad 0\le\rho\le2.
		\]
		In the normalized variables, the bootstrap bound gives
		\begin{equation}\label{3.2.0-normal}
			\Xi_0(1+E_0)^{3/2}
			\sup_{0\le s\le  T}
			\left(
			\|\nabla u(s)\|_{L^2}^2+\|P(s)\|_{L^2}^2
			\right)
			\le 2\varepsilon .
		\end{equation}
		
		Recall that the momentum equation \eqref{1.1}$_{2}$ can be written as
		\begin{equation}\label{3.2.2}
			\rho \dot u-\mu\Delta u-(\mu+\lambda)\nabla\div u+\nabla P=0.
		\end{equation}
		Multiplying \eqref{3.2.2} by \(\dot u\) and integrating over \(\mathbb R^2\), we obtain
		\begin{equation}\label{3.2.3}
			\int_{\mathbb R^2}\rho|\dot u|^2\,dx
			=
			I_1+I_2+I_3,
		\end{equation}
		where
		\[
		I_1:=\mu\int_{\mathbb R^2}\Delta u\cdot\dot u\,dx,
		\qquad
		I_2:=(\mu+\lambda)\int_{\mathbb R^2}\nabla\div u\cdot\dot u\,dx,
		\]
		and
		\[
		I_3:=-\int_{\mathbb R^2}\nabla P\cdot\dot u\,dx.
		\]
		
		For \(I_1\), since \(\dot u=u_t+u\cdot\nabla u\), integration by parts gives
		\[
		I_1
		=
		-\mu\int_{\mathbb R^2}\nabla u:\nabla u_t\,dx
		-\mu\int_{\mathbb R^2}\partial_i u^j\,\partial_i(u\cdot\nabla u^j)\,dx.
		\]
		Thus
		\[
		-\mu\int_{\mathbb R^2}\nabla u:\nabla u_t\,dx
		=
		-\frac{\mu}{2}\frac{d}{dt}\int_{\mathbb R^2}|\nabla u|^2\,dx.
		\]
		Moreover,
		\[
		\partial_i(u\cdot\nabla u^j)
		=
		\partial_i u^k\,\partial_k u^j+u^k\partial_{ik}u^j.
		\]
		Hence
		\[
		\left|
		-\mu\int_{\mathbb R^2}\partial_i u^j\,\partial_i(u\cdot\nabla u^j)\,dx
		\right|
		\le C\int_{\mathbb R^2}|\nabla u|^3\,dx,
		\]
		which implies that
		\begin{equation}\label{3.2.4}
			I_1
			\le
			-\frac{\mu}{2}\frac{d}{dt}\|\nabla u\|_{L^2}^2
			+
			C\int_{\mathbb R^2}|\nabla u|^3\,dx.
		\end{equation}
		
		For \(I_2\), we have
		\[
		I_2
		=-(\mu+\lambda)\int_{\mathbb R^2}\div u\,\div \dot u\,dx.
		\]
		Since
		\[
		\div\dot u=(\div u)_t+\div(u\cdot\nabla u),
		\]
		we obtain
		\[
		I_2
		=
		-\frac{\mu+\lambda}{2}\frac{d}{dt}\|\div u\|_{L^2}^2
		-(\mu+\lambda)\int_{\mathbb R^2}\div u\,\div(u\cdot\nabla u)\,dx.
		\]
		Using
		\[
		\div(u\cdot\nabla u)
		=
		u\cdot\nabla(\div u)+\partial_i u^j\partial_j u^i,
		\]
		and integrating by parts in the transport term, we get
		\[
		\left|
		\int_{\mathbb R^2}\div u\,\div(u\cdot\nabla u)\,dx
		\right|
		\le C\int_{\mathbb R^2}|\nabla u|^3\,dx.
		\]
		Thus
		\begin{equation}\label{3.2.5}
			I_2
			\le
			-\frac{\mu+\lambda}{2}\frac{d}{dt}\|\div u\|_{L^2}^2
			+
			C\int_{\mathbb R^2}|\nabla u|^3\,dx.
		\end{equation}
		
		For \(I_3\), integrating by parts gives
		\[
		I_3=\int_{\mathbb R^2}P\,\div\dot u\,dx.
		\]
		Using again
		\[
		\div\dot u=(\div u)_t+\div(u\cdot\nabla u),
		\]
		we have
		\[
		I_3
		=
		\int_{\mathbb R^2}P(\div u)_t\,dx
		+
		\int_{\mathbb R^2}P\,\div(u\cdot\nabla u)\,dx.
		\]
		Since
		\[
		\int_{\mathbb R^2}P(\div u)_t\,dx
		=
		\frac{d}{dt}\int_{\mathbb R^2}P\div u\,dx
		-
		\int_{\mathbb R^2}P_t\div u\,dx.
		\]
		and
		\[
		P_t+u\cdot\nabla P+\gamma P\div u=0,
		\]
		we get
		\[
		-\int_{\mathbb R^2}P_t\div u\,dx
		=
		\int_{\mathbb R^2}(u\cdot\nabla P)\div u\,dx
		+
		\gamma\int_{\mathbb R^2}P(\div u)^2\,dx.
		\]
		On the other hand,
		\[
		\int_{\mathbb R^2}P\,\div(u\cdot\nabla u)\,dx
		=
		\int_{\mathbb R^2}P\,u\cdot\nabla(\div u)\,dx
		+
		\int_{\mathbb R^2}P\,\partial_i u^j\partial_j u^i\,dx.
		\]
		The two transport terms satisfy
		\[
		\int_{\mathbb R^2}(u\cdot\nabla P)\div u\,dx
		+
		\int_{\mathbb R^2}P\,u\cdot\nabla(\div u)\,dx
		=
		-\int_{\mathbb R^2}P(\div u)^2\,dx.
		\]
		Therefore,
		\[
		I_3
		=
		\frac{d}{dt}\int_{\mathbb R^2}P\div u\,dx
		+
		(\gamma-1)\int_{\mathbb R^2}P(\div u)^2\,dx
		+
		\int_{\mathbb R^2}P\,\partial_i u^j\partial_j u^i\,dx,
		\]
		and hence
		\begin{equation}\label{3.2.6}
			I_3
			\le
			\frac{d}{dt}\int_{\mathbb R^2}P\div u\,dx
			+
			C\int_{\mathbb R^2}P|\nabla u|^2\,dx.
		\end{equation}
		
		Substituting \eqref{3.2.4}, \eqref{3.2.5}, and \eqref{3.2.6} into \eqref{3.2.3}, we obtain
		\begin{equation}\label{3.2.7}
			\begin{aligned}
				&\int_{\mathbb R^2}\rho|\dot u|^2\,dx
				+
				\frac{d}{dt}
				\int_{\mathbb R^2}
				\left(
				\frac{\mu}{2}|\nabla u|^2
				+\frac{\mu+\lambda}{2}(\div u)^2
				-P\div u
				\right)\,dx \\
				\le{}&
				C\int_{\mathbb R^2}|\nabla u|^3\,dx
				+
				C\int_{\mathbb R^2}P|\nabla u|^2\,dx.
			\end{aligned}
		\end{equation}
		
		Recall that
		\[
		F:=(2\mu+\lambda)\div u-P,
		\qquad
		\omega:=\partial_1u^{2}-\partial_2u^{1}.
		\]
		By the $L^{p}$ estimate,
		\[
		\|\nabla u\|_{L^3}
		\le
		C\Bigl(
		\|F\|_{L^3}+\|\omega\|_{L^3}+\|P\|_{L^3}
		\Bigr).
		\]
		We estimate \(F\) and \(\omega\) in \(L^3\). The interpolation inequality
		\[
		\|f\|_{L^3}^3
		\le
		C\|f\|_{L^2}\|\nabla f\|_{L^{4/3}}^2
		\]
		gives
		\[
		\|F\|_{L^3}^3+\|\omega\|_{L^3}^3
		\le
		C\|F\|_{L^2}\|\nabla F\|_{L^{4/3}}^2
		+
		C\|\omega\|_{L^2}\|\nabla\omega\|_{L^{4/3}}^2.
		\]
		By Lemma \ref{lem2.9} with \(p=2\), and by the \(L^p\)-extension stated after
		Lemma \ref{lem2.9} with \(p=4/3\),
		\[
		\|F\|_{L^2}+\|\omega\|_{L^2}
		\le
		C\left(\|\nabla u\|_{L^2}+\|P\|_{L^2}\right),
		\]
		and
		\[
		\|\nabla F\|_{L^{4/3}}+\|\nabla\omega\|_{L^{4/3}}
		\le
		C\|\rho\dot u\|_{L^{4/3}}.
		\]
		Since \(0\le\rho\le2\),
		\[
		\begin{aligned}
			\|\rho\dot u\|_{L^{4/3}}^{4/3}
			&=
			\int_{\mathbb R^2}
			(\rho|\dot u|^2)^{2/3}\rho^{2/3}\,dx \\
			&\le
			\left(
			\int_{\mathbb R^2}\rho|\dot u|^2\,dx
			\right)^{2/3}
			\left(
			\int_{\mathbb R^2}\rho^2\,dx
			\right)^{1/3}  \\
			&\le
			C
			\left(
			\int_{\mathbb R^2}\rho|\dot u|^2\,dx
			\right)^{2/3}
			\left(
			\int_{\mathbb R^2}\rho\,dx
			\right)^{1/3}.
		\end{aligned}
		\]
		Thus
		\[
		\|\rho\dot u\|_{L^{4/3}}^2
		\le
		C
		\left(
		\int_{\mathbb R^2}\rho\,dx
		\right)^{1/2}
		\int_{\mathbb R^2}\rho|\dot u|^2\,dx.
		\]
		Combining the above estimates, we obtain
		\[
		\|\nabla u\|_{L^3}^3
		\le
		C
		\left(
		\int_{\mathbb R^2}\rho\,dx
		\right)^{1/2}
		\Bigl(
		\|\nabla u\|_{L^2}+\|P\|_{L^2}
		\Bigr)
		\int_{\mathbb R^2}\rho|\dot u|^2\,dx
		+
		C\|P\|_{L^3}^3.
		\]
		Observe that
		\[
		\int_{\mathbb R^2}\rho\,dx=M_0.
		\]
		Moreover, by \eqref{3.2.0-normal} and \(\Xi_0\ge M_0\), it follows that
		\[
		\begin{aligned}
			&\left(
			\int_{\mathbb R^2}\rho\,dx
			\right)^{1/2}
			\Bigl(
			\|\nabla u(t)\|_{L^2}+\|P(t)\|_{L^2}
			\Bigr) \\
			\le&
			C
			M_0^{1/2}
			\Bigl(
			\|\nabla u(t)\|_{L^2}^2+\|P(t)\|_{L^2}^2
			\Bigr)^{1/2}  \\
			\le&
			C
			\Xi_0^{1/2}
			\frac{(2\varepsilon)^{1/2}}{\Xi_0^{1/2}(1+E_0)^{3/4}}
			\le
			C\varepsilon^{1/2}.
		\end{aligned}
		\]
		Therefore,
		\begin{equation}\label{grad-u-L3-small}
			\|\nabla u\|_{L^3}^3
			\le
			C\varepsilon^{1/2}
			\int_{\mathbb R^2}\rho|\dot u|^2\,dx
			+
			C\|P\|_{L^3}^3.
		\end{equation}
		By H\"older's inequality and Young's inequality,
		\[
		\int_{\mathbb R^2}P|\nabla u|^2\,dx
		\le
		\|P\|_{L^3}\|\nabla u\|_{L^3}^2
		\le
		\delta\|\nabla u\|_{L^3}^3
		+
		C_\delta\|P\|_{L^3}^3.
		\]
		Using \eqref{grad-u-L3-small}, we get
		\begin{equation}\label{P-grad-u-small}
			\int_{\mathbb R^2}P|\nabla u|^2\,dx
			\le
			C\varepsilon^{1/2}
			\int_{\mathbb R^2}\rho|\dot u|^2\,dx
			+
			C\|P\|_{L^3}^3.
		\end{equation}
		Substituting \eqref{grad-u-L3-small} and \eqref{P-grad-u-small} into
		\eqref{3.2.7}, we obtain
		\begin{equation}\label{short-time-diff-pre}
			\begin{aligned}
				&\int_{\mathbb R^2}\rho|\dot u|^2\,dx
				+
				\frac{d}{dt}
				\int_{\mathbb R^2}
				\left(
				\frac{\mu}{2}|\nabla u|^2
				+\frac{\mu+\lambda}{2}(\div u)^2
				-P\div u
				\right)\,dx \\
				\le&
				C\varepsilon^{1/2}
				\int_{\mathbb R^2}\rho|\dot u|^2\,dx
				+
				C\|P\|_{L^3}^3.
			\end{aligned}
		\end{equation}
		Choosing \(\varepsilon_1>0\) sufficiently small such that
		\(C\varepsilon_1^{1/2}\le1/2\) for \(0<\varepsilon\le\varepsilon_1\),
		\eqref{short-time-diff-pre} gives
		\begin{equation}\label{short-time-diff-final}
			\frac12\int_{\mathbb R^2}\rho|\dot u|^2\,dx
			+
			\frac{d}{dt}
			\int_{\mathbb R^2}
			\left(
			\frac{\mu}{2}|\nabla u|^2
			+\frac{\mu+\lambda}{2}(\div u)^2
			-P\div u
			\right)\,dx
			\le
			C\|P\|_{L^3}^3.
		\end{equation}
		Integrating \eqref{short-time-diff-final} over \([0,t]\), and using
		\[
		\left|\int_{\mathbb R^2}P\div u\,dx\right|
		\le
		\frac{\mu}{4}\|\nabla u\|_{L^2}^2+C\|P\|_{L^2}^2
		\]
		at both \(t\) and \(0\), gives
		\[
		\begin{aligned}
			&\|\nabla u(t)\|_{L^2}^2
			+\int_0^t\int_{\mathbb R^2}\rho|\dot u|^2\,dxds\\
			\le&
			C\left(\|\nabla u_0\|_{L^2}^2+\|P_0\|_{L^2}^2\right)
			+C\|P(t)\|_{L^2}^2
			+C\int_0^t\|P(s)\|_{L^3}^3\,ds.
		\end{aligned}
		\]
		This proves the desired estimate in the normalized case, which yields \eqref{eq:short-time-u-normalized}.
	\end{proof}

	The estimate in Lemma \ref{lem3.2} is not yet closed, since its right-hand side still
	contains the pressure terms
	\[
	\|P(t)\|_{L^2}^2+\int_0^t\|P(s)\|_{L^3}^3\,ds .
	\]
	We next derive a complementary estimate for these pressure terms.
	\begin{lemma}\label{lem3.3.1}
		Let \((\rho,u)\) solve \eqref{1.1}--\eqref{1.4.1} smoothly on
		\(\mathbb R^2\times(0,T]\) satisfying \eqref{3.1.3} and the assumptions of
		Theorem \ref{th1.1}. Let \(C_1\ge1\) be the constant in Lemma $\ref{lem3.2}$. Then, for any
		\(t\in[0,T]\), there holds
		\begin{equation}\label{3.3.1}
			\begin{aligned}
				&\|P(t)\|_{L^2}^2
				+\frac{2\gamma-1}{2(2\mu+\lambda)}
				\int_0^t\|P(s)\|_{L^3}^3\,ds \\
				\le{}&
				\|P_0\|_{L^2}^2
				+
				\frac{1}{
					16C_1\left(1+\frac{2(2\mu+\lambda)}{2\gamma-1}\right)}
				\int_0^t\int_{\mathbb R^2}\rho|\dot u|^2\,dxds,
			\end{aligned}
		\end{equation}
		where $C_{1}$ depends only on
		\(\mu,\lambda,\gamma\).
	\end{lemma}
	
	\begin{proof}
		As in the proof of Lemma \ref{lem3.2}, after applying the scaling
		\[
		\widetilde\rho(y,\tau)
		=
		\bar\rho^{-1}\rho\!\left(\bar\rho^{-\frac{\gamma+1}{2}}y,\bar\rho^{-\gamma}\tau\right),
		\qquad
		\widetilde u(y,\tau)
		=
		\bar\rho^{-\frac{\gamma-1}{2}}
		u\!\left(\bar\rho^{-\frac{\gamma+1}{2}}y,\bar\rho^{-\gamma}\tau\right),
		\]
		it suffices to prove the estimate in the normalized variables. Thus, in the following,
		we work under
		\[
		\bar\rho=1,\qquad 0\le \rho\le2.
		\]
		In the normalized variables, \eqref{3.1.3} gives
		\[
		\Xi_0(1+E_0)^{3/2}
		\sup_{0\le s\le T}
		\left(
		\|\nabla u(s)\|_{L^2}^2+\|P(s)\|_{L^2}^2
		\right)
		\le 2\varepsilon.
		\]
		
		Since \(P=\rho^\gamma\), the continuity equation gives
		\begin{equation}\label{3.3.2}
			P_t+u\cdot\nabla P+\gamma P\div u=0.
		\end{equation}
		Multiplying \eqref{3.3.2} by \(2P\) and integrating over \(\mathbb R^2\), we obtain
		\[
		\frac{d}{dt}\|P\|_{L^2}^2
		+
		(2\gamma-1)\int_{\mathbb R^2}P^2\div u\,dx
		=0.
		\]
		Note that
		\[
		\div u=\frac{F+P}{2\mu+\lambda}.
		\]
		Hence
		\begin{equation}\label{3.3.3}
			\frac{d}{dt}\|P\|_{L^2}^2
			+
			\frac{2\gamma-1}{2\mu+\lambda}\|P\|_{L^3}^3
			=
			-\frac{2\gamma-1}{2\mu+\lambda}
			\int_{\mathbb R^2}FP^2\,dx.
		\end{equation}
		
		We estimate the last term. By H\"older's inequality,
		\[
		\left|\int_{\mathbb R^2}FP^2\,dx\right|
		\le
		\|F\|_{L^3}\|P\|_{L^3}^2.
		\]
		Using the interpolation inequality
		\[
		\|F\|_{L^3}^3
		\le
		C\|F\|_{L^2}\|\nabla F\|_{L^{4/3}}^2,
		\]
		Lemma \ref{lem2.9} with \(p=2\), and the \(L^p\)-extension in the remark after
		Lemma \ref{lem2.9} with \(p=4/3\), we get
		\[
		\|F\|_{L^3}
		\le
		C M_0^{1/6}
		\left(\|\nabla u\|_{L^2}^2+\|P\|_{L^2}^2\right)^{1/6}
		\|\sqrt\rho\,\dot u\|_{L^2}^{2/3}.
		\]
		Therefore,
		\[
		\left|\int_{\mathbb R^2}FP^2\,dx\right|
		\le
		C M_0^{1/6}
		\left(\|\nabla u\|_{L^2}^2+\|P\|_{L^2}^2\right)^{1/6}
		\|\sqrt\rho\,\dot u\|_{L^2}^{2/3}
		\|P\|_{L^3}^2.
		\]
		By Young's inequality with the fixed coefficient
		\[
		\frac{2\mu+\lambda}{
			16C_1(2\gamma-1)\left(1+\frac{2(2\mu+\lambda)}{2\gamma-1}\right)},
		\]
		we get
		\[
		\begin{aligned}
			\left|\int_{\mathbb R^2}FP^2\,dx\right|
			\le{}&
			\frac{2\mu+\lambda}{
				16C_1(2\gamma-1)\left(1+\frac{2(2\mu+\lambda)}{2\gamma-1}\right)}
			\|\sqrt\rho\,\dot u\|_{L^2}^2\\
			&+
			C
			M_0^{1/4}
			\left(\|\nabla u\|_{L^2}^2+\|P\|_{L^2}^2\right)^{1/4}
			\|P\|_{L^3}^3,
		\end{aligned}
		\]
		where \(C\) depends only on \(\mu,\lambda,\gamma\). Using
		\[
		M_0^{1/4}
		\left(\|\nabla u\|_{L^2}^2+\|P\|_{L^2}^2\right)^{1/4}
		\le
		C\varepsilon^{1/4},
		\]
		which follows from  \(\Xi_0\ge M_0\), and choosing
		\(\varepsilon>0\) sufficiently small such that
		\[
		C\varepsilon^{1/4}\le \frac12,
		\]
		we obtain
		\begin{equation}\label{3.3.4}
			\begin{aligned}
				\left|\int_{\mathbb R^2}FP^2\,dx\right|
				\le{}&
				\frac{2\mu+\lambda}{
					16C_1(2\gamma-1)\left(1+\frac{2(2\mu+\lambda)}{2\gamma-1}\right)}
				\int_{\mathbb R^2}\rho|\dot u|^2\,dx+
				\frac12\|P\|_{L^3}^3.
			\end{aligned}
		\end{equation}
		Substituting \eqref{3.3.4} into \eqref{3.3.3}, we infer
		\[
		\frac{d}{dt}\|P\|_{L^2}^2
		+
		\frac{2\gamma-1}{2(2\mu+\lambda)}
		\|P\|_{L^3}^3
		\le
		\frac{1}{
			16C_1\left(1+\frac{2(2\mu+\lambda)}{2\gamma-1}\right)}
		\int_{\mathbb R^2}\rho|\dot u|^2\,dx.
		\]
		Integrating over \([0,t]\) gives \eqref{3.3.1} in the normalized variables. Scaling back to
		the original variables proves \eqref{3.3.1}. The proof is complete.
	\end{proof}
	
	Combining  Lemma \ref{lem3.2} with
	Lemma \ref{lem3.3.1}, we have the following corollary.
	
	\begin{corollary}\label{cor3.1.1}
		Let \((\rho,u)\) solve \eqref{1.1}--\eqref{1.4.1} smoothly on
		\(\mathbb R^2\times(0,T]\) satisfying \eqref{3.1.3} and the assumptions of
		Theorem \ref{th1.1}. Then, for any \(t\in[0,T]\),
		\begin{equation}\label{3.4.1}
			\begin{aligned}
				&\sup_{0\le s\le t}
				\left(
				\|\nabla u(s)\|_{L^2}^2+\|P(s)\|_{L^2}^2
				\right)
				+\int_0^t\int_{\mathbb R^2}\rho|\dot u|^2\,dxds
				+\int_0^t\|P(s)\|_{L^3}^3\,ds \\
				\le{}&
				C\left(
				\|\nabla u_0\|_{L^2}^2+\|P_0\|_{L^2}^2
				\right),
			\end{aligned}
		\end{equation}
		where \(C\) depends only on \(\mu,\lambda,\gamma\).
	\end{corollary}
	
	\begin{proof}
		Let \(C_1\ge1\) be the constant in Lemma \ref{lem3.2}. By Lemma \ref{lem3.2}, for any
		\(s\in[0,t]\), one has
		\begin{equation}\label{3.4.2}
			\begin{aligned}
				&\|\nabla u(s)\|_{L^2}^2
				+\int_0^s\int_{\mathbb R^2}\rho|\dot u|^2\,dxdr \\
				\le{}&
				C_1\left(
				\|\nabla u_0\|_{L^2}^2+\|P_0\|_{L^2}^2
				\right)
				+C_1\|P(s)\|_{L^2}^2
				+C_1\int_0^s\|P(r)\|_{L^3}^3\,dr.
			\end{aligned}
		\end{equation}
		On the other hand, by Lemma \ref{lem3.3.1},
		\begin{equation}\label{3.4.3}
			\begin{aligned}
				&\|P(s)\|_{L^2}^2
				+\frac{2\gamma-1}{2(2\mu+\lambda)}
				\int_0^s\|P(r)\|_{L^3}^3\,dr \\
				\le{}&
				\|P_0\|_{L^2}^2
				+
				\frac{1}{
					16C_1\left(1+\frac{2(2\mu+\lambda)}{2\gamma-1}\right)}
				\int_0^s\int_{\mathbb R^2}\rho|\dot u|^2\,dxdr.
			\end{aligned}
		\end{equation}
		Multiplying \eqref{3.4.3} by $4C_1\left(1+\frac{2(2\mu+\lambda)}{2\gamma-1}\right)$ and adding the result to \eqref{3.4.2}, we obtain
		\[
		\begin{aligned}
			&\|\nabla u(s)\|_{L^2}^2
			+\frac34\int_0^s\int_{\mathbb R^2}\rho|\dot u|^2\,dxdr \\
			&+
			\left[
			4C_1\left(1+\frac{2(2\mu+\lambda)}{2\gamma-1}\right)-C_1
			\right]\|P(s)\|_{L^2}^2\\
			&+
			\left[
			4C_1\left(1+\frac{2(2\mu+\lambda)}{2\gamma-1}\right)
			\frac{2\gamma-1}{2(2\mu+\lambda)}
			-C_1
			\right]\int_0^s\|P(r)\|_{L^3}^3\,dr\\
			\le&
			C_1\left(
			\|\nabla u_0\|_{L^2}^2+\|P_0\|_{L^2}^2
			\right)
			+
			4C_1\left(1+\frac{2(2\mu+\lambda)}{2\gamma-1}\right)
			\|P_0\|_{L^2}^2.
		\end{aligned}
		\]
		Since
		\[
		4C_1\left(1+\frac{2(2\mu+\lambda)}{2\gamma-1}\right)-C_1\ge C_1
		\]
		and
		\[
		4C_1\left(1+\frac{2(2\mu+\lambda)}{2\gamma-1}\right)
		\frac{2\gamma-1}{2(2\mu+\lambda)}
		-C_1
		=
		4C_1\left(1+\frac{2\gamma-1}{2(2\mu+\lambda)}\right)-C_1
		\ge C_1,
		\]
		we arrive at
		\[
		\begin{aligned}
			&\|\nabla u(s)\|_{L^2}^2+\|P(s)\|_{L^2}^2
			+\int_0^s\int_{\mathbb R^2}\rho|\dot u|^2\,dxdr
			+\int_0^s\|P(r)\|_{L^3}^3\,dr\\
			&\le
			C\left(
			\|\nabla u_0\|_{L^2}^2+\|P_0\|_{L^2}^2
			\right),
		\end{aligned}
		\]
		where \(C\) depends only on \(\mu,\lambda,\gamma\). Taking the supremum over
		\(0\le s\le t\) yields \eqref{3.4.1}. The proof is complete.
	\end{proof}
	
	We next establish a weighted estimate for the material derivative, which shows the decay in time.
	\begin{lemma}\label{lem3.5}
		Let \((\rho,u)\) solve \eqref{1.1}--\eqref{1.4.1} smoothly on
		\(\mathbb R^2\times(0,T]\) satisfying \eqref{3.1.3} and the assumptions of
		Theorem \ref{th1.1}. Define
		\[
		\sigma(t):=\min\{1,\bar\rho^\gamma t\}.
		\]
		Then, for any \(t\in[0,T]\), there holds
		\begin{equation}\label{3.5.1}
			\begin{aligned}
				&\sup_{0\le s\le t}
				\sigma(s)^2
				\int_{\mathbb R^2}\rho|\dot u|^2\,dx
				+
				\int_0^t
				\sigma(s)^2
				\int_{\mathbb R^2}|\nabla\dot u|^2\,dxds\\
				&\qquad\le
				C\bar\rho^\gamma
				\left(
				\|\nabla u_0\|_{L^2}^2+\|P_0\|_{L^2}^2
				\right),
			\end{aligned}
		\end{equation}
		where \(\dot u=u_t+u\cdot\nabla u\), \(P_0=\rho_0^\gamma\), and \(C\) depends only on
		\(\mu,\lambda,\gamma\).
	\end{lemma}
	
	\begin{proof}
		We first normalize the upper bound of the density. Clearly, it is enough to prove that, for every \(\tau\in[0,\bar\rho^\gamma T]\),
		\begin{equation}\label{3.5.5}
			\begin{aligned}
				&\sup_{0\le s\le \tau}
				\widetilde\sigma(s)^2
				\int_{\mathbb R^2}\widetilde\rho|\dot{\widetilde u}|^2\,dy
				+
				\int_0^\tau
				\widetilde\sigma(s)^2
				\int_{\mathbb R^2}|\nabla_y\dot{\widetilde u}|^2\,dyds\\
				&\qquad\le
				C\left(
				\|\nabla_y\widetilde u_0\|_{L^2_y}^2
				+\|\widetilde P_0\|_{L^2_y}^2
				\right).
			\end{aligned}
		\end{equation}
		Indeed, rewriting  \eqref{3.5.5} in the original variables by the above scaling
		relations and multiplying
		both sides by \(\bar\rho^{2\gamma-1}\), gives exactly \eqref{3.5.1}.
		
		We now prove \eqref{3.5.5}.  The normalized momentum equation \eqref{1.1}$_2$
		can be rewritten as
		\begin{equation}\label{3.5.6}
			\widetilde\rho\,\dot{\widetilde u}
			-\mu\Delta_y\widetilde u
			-(\mu+\lambda)\nabla_y\div_y\widetilde u
			+\nabla_y\widetilde P=0.
		\end{equation}
		Applying \(\partial_s+\widetilde u\cdot\nabla_y\) to \eqref{3.5.6}, multiplying the
		resulting equation by \(\widetilde\sigma(s)^2\dot{\widetilde u}\), and integrating by parts
		over \(\mathbb R^2\), we obtain
		\begin{equation}\label{3.5.7}
			\begin{aligned}
				&\frac{d}{ds}
				\left[
				\widetilde\sigma(s)^2
				\int_{\mathbb R^2}\widetilde\rho|\dot{\widetilde u}|^2\,dy
				\right]
				+
				\mu\widetilde\sigma(s)^2
				\int_{\mathbb R^2}|\nabla_y\dot{\widetilde u}|^2\,dy\\
				\le{}&
				C\widetilde\sigma(s)\widetilde\sigma'(s)
				\int_{\mathbb R^2}\widetilde\rho|\dot{\widetilde u}|^2\,dy
				+
				C\widetilde\sigma(s)^2
				\left(
				\|\nabla_y\widetilde u\|_{L^4_y}^4
				+
				\|\widetilde P\|_{L^4_y}^4
				\right),
			\end{aligned}
		\end{equation}
		where \(C>0\) depends only on \(\mu,\lambda,\gamma\). Observe that
		\[
		0\le \widetilde\sigma\le 1,\qquad 0\le \widetilde\sigma'\le 1.
		\]
		
		By Lemma \ref{lem2.9} with \(p=4\),
		\[
		\|\nabla_y\widetilde u\|_{L^4_y}
		\le
		C\|\widetilde\rho\dot{\widetilde u}\|_{L^2_y}^{1/2}
		\left(
		\|\nabla_y\widetilde u\|_{L^2_y}
		+
		\|\widetilde P\|_{L^2_y}
		\right)^{1/2}
		+
		C\|\widetilde P\|_{L^4_y}.
		\]
		Since \(0\le\widetilde\rho\le2\),
		\[
		\|\widetilde\rho\dot{\widetilde u}\|_{L^2_y}^2
		\le
		C\int_{\mathbb R^2}\widetilde\rho|\dot{\widetilde u}|^2\,dy.
		\]
		Also, since \(0\le\widetilde P=\widetilde\rho^\gamma\le2^\gamma\),
		\[
		\|\widetilde P\|_{L^4_y}^4
		\le
		C\|\widetilde P\|_{L^3_y}^3.
		\]
		Hence we have
		\begin{equation}\label{3.5.8}
			\begin{aligned}
				&\|\nabla_y\widetilde u\|_{L^4_y}^4
				+
				\|\widetilde P\|_{L^4_y}^4
				\le
				C\left(
				\|\nabla_y\widetilde u\|_{L^2_y}^2
				+
				\|\widetilde P\|_{L^2_y}^2
				\right)
				\int_{\mathbb R^2}\widetilde\rho|\dot{\widetilde u}|^2\,dy
				+
				C\|\widetilde P\|_{L^3_y}^3.
			\end{aligned}
		\end{equation}
		
		By the normalized form of Corollary \ref{cor3.1.1}, for every
		\(\tau\in[0,\bar\rho^\gamma T]\), we obtain
		\begin{equation}\label{3.5.9}
			\begin{aligned}
				&\sup_{0\le s\le \tau}
				\left(
				\|\nabla_y\widetilde u(s)\|_{L^2_y}^2
				+
				\|\widetilde P(s)\|_{L^2_y}^2
				\right)
				+
				\int_0^\tau
				\int_{\mathbb R^2}\widetilde\rho|\dot{\widetilde u}|^2\,dyds
				+
				\int_0^\tau
				\|\widetilde P(s)\|_{L^3_y}^3\,ds\\
				\le&
				C\left(
				\|\nabla_y\widetilde u_0\|_{L^2_y}^2
				+
				\|\widetilde P_0\|_{L^2_y}^2
				\right).
			\end{aligned}
		\end{equation}
		Moreover, by the normalized smallness assumption \eqref{3.1.3}, we get
		\[
		\sup_{0\le s\le \bar\rho^\gamma T}
		\left(
		\|\nabla_y\widetilde u(s)\|_{L^2_y}^2
		+
		\|\widetilde P(s)\|_{L^2_y}^2
		\right)
		\le 2\varepsilon.
		\]
		Using \eqref{3.5.8}, \eqref{3.5.9}, and choosing \(\varepsilon>0\) sufficiently small, we obtain
		\begin{equation}\label{3.5.10}
			\begin{aligned}
				&\int_0^\tau
				\widetilde\sigma(s)^2
				\left(
				\|\nabla_y\widetilde u(s)\|_{L^4_y}^4
				+
				\|\widetilde P(s)\|_{L^4_y}^4
				\right)\,ds\le
				C\left(
				\|\nabla_y\widetilde u_0\|_{L^2_y}^2
				+
				\|\widetilde P_0\|_{L^2_y}^2
				\right).
			\end{aligned}
		\end{equation}
		Integrating \eqref{3.5.7} over \([0,\tau]\), and using \(\widetilde\sigma(0)=0\), we get
		\begin{equation}\label{3.5.11}
			\begin{aligned}
				&\widetilde\sigma(\tau)^2
				\int_{\mathbb R^2}\widetilde\rho|\dot{\widetilde u}|^2(y,\tau)\,dy
				+
				\mu\int_0^\tau
				\widetilde\sigma(s)^2
				\int_{\mathbb R^2}|\nabla_y\dot{\widetilde u}|^2\,dyds\\
				\le{}&
				C\int_0^\tau
				\widetilde\sigma(s)\widetilde\sigma'(s)
				\int_{\mathbb R^2}\widetilde\rho|\dot{\widetilde u}|^2\,dyds+
				C\int_0^\tau
				\widetilde\sigma(s)^2
				\left(
				\|\nabla_y\widetilde u(s)\|_{L^4_y}^4
				+
				\|\widetilde P(s)\|_{L^4_y}^4
				\right)\,ds .
			\end{aligned}
		\end{equation}
		By \eqref{3.5.9} and \eqref{3.5.10}, we have
		\begin{equation}\label{3.5.12}
			\begin{aligned}
				&\widetilde\sigma(\tau)^2
				\int_{\mathbb R^2}\widetilde\rho|\dot{\widetilde u}|^2(y,\tau)\,dy
				+
				\int_0^\tau
				\widetilde\sigma(s)^2
				\int_{\mathbb R^2}|\nabla_y\dot{\widetilde u}|^2\,dyds\le
				C\left(
				\|\nabla_y\widetilde u_0\|_{L^2_y}^2
				+
				\|\widetilde P_0\|_{L^2_y}^2
				\right).
			\end{aligned}
		\end{equation}
		Taking the supremum over \(0\le s\le\tau\), we obtain \eqref{3.5.5}.
	\end{proof}

	Next, we derive the rates of decay in time for \(\nabla u\), \(P\), and
	\(\rho^{1/2}\dot u\). These estimates will be used to obtain the uniform-in-time
	upper bound of the density for large time.
	\begin{lemma}\label{lem:large-time-L2-decay}
		Let \((\rho,u)\) solve \eqref{1.1}--\eqref{1.4.1} smoothly on
		\(\mathbb R^2\times(0,T]\) satisfying \eqref{3.1.3} and the assumptions of
		Theorem \ref{th1.1}.
		Then, for any \(t\in[\bar\rho^{-\gamma},T]\), there holds
		\begin{equation}\label{large-time-L2-decay-original}
			\begin{aligned}
				&\sup_{\bar\rho^{-\gamma}\le s\le t}
				\Bigg[
				\bar\rho\,s
				\big(\|\nabla u(s)\|_{L^2}^2+\|P(s)\|_{L^2}^2\big)
				+
				\bar\rho^{1-\gamma} s
				\int_{\mathbb R^2}\rho|\dot u|^2\,dx
				\Bigg] \\
				&+
				\int_{\bar\rho^{-\gamma}}^t
				\bar\rho^{1-\gamma} s
				\int_{\mathbb R^2}|\nabla\dot u|^2\,dxds
				+
				\int_{\bar\rho^{-\gamma}}^t
				\bar\rho\,s
				\|P(s)\|_{L^3}^3\,ds \\
				\le&
				C\Bigg[
				\bar\rho E_0
				+
				\big(\bar\rho^\gamma M_0 + (\bar\rho^\gamma M_0)^{-1}\big)
				\bar\rho^{1-\gamma}
				\big(\|\nabla u_0\|_{L^2}^2+\|P_0\|_{L^2}^2\big)
				\Bigg],
			\end{aligned}
		\end{equation}
		where \(\dot u=u_t+u\cdot\nabla u\), \(P_0=\rho_0^\gamma\), and \(C\) depends only on
		\(\mu,\lambda,\gamma\).
	\end{lemma}
	
	\begin{proof}
		We first normalize the upper bound of the density.
		Thus it suffices to prove that, for any \(1\le \tau\le \bar\rho^\gamma T\),
		\begin{equation}\label{large-time-L2-decay-normalized}
			\begin{aligned}
				&\sup_{1\le s\le \tau}
				\Bigg[
				s
				\left(
				\|\nabla_y\widetilde u(s)\|_{L^2_y}^2
				+
				\|\widetilde P(s)\|_{L^2_y}^2
				\right)
				+
				s\int_{\mathbb R^2}\widetilde\rho|\dot{\widetilde u}|^2\,dy
				\Bigg]\\
				&
				+
				\int_1^\tau
				s\int_{\mathbb R^2}|\nabla_y\dot{\widetilde u}|^2\,dyds
				+
				\int_1^\tau
				s\|\widetilde P(s)\|_{L^3_y}^3\,ds\\
				\le&
				C\Bigg[
				\bar\rho E_0
				+
				\left(
				\bar\rho^\gamma M_0
				+
				(\bar\rho^\gamma M_0)^{-1}
				\right)
				\bar\rho^{1-\gamma}
				\left(
				\|\nabla u_0\|_{L^2}^2+\|P_0\|_{L^2}^2
				\right)
				\Bigg].
			\end{aligned}
		\end{equation}
		
		We now prove \eqref{large-time-L2-decay-normalized}. Consider
		\[
		\begin{aligned}
			\mathcal H(s):=
			\int_{\mathbb R^2}
			\left(
			\frac{\mu}{2}|\nabla_y\widetilde u|^2
			+
			\frac{\mu+\lambda}{2}(\div_y\widetilde u)^2
			-
			\widetilde P\,\div_y\widetilde u
			\right)dy
			+
			C_2\|\widetilde P(s)\|_{L^2_y}^2,
		\end{aligned}
		\]
		where \(C_2>0\) is chosen sufficiently large depending only on
		\(\mu,\lambda,\gamma\). By Cauchy's inequality,
		\[
		\left|
		\int_{\mathbb R^2}\widetilde P\div_y\widetilde u\,dy
		\right|
		\le
		\frac{\mu}{4}\|\nabla_y\widetilde u\|_{L^2_y}^2
		+
		C\|\widetilde P\|_{L^2_y}^2.
		\]
		Thus
		\begin{equation}\label{H-equivalence}
			c
			\left(
			\|\nabla_y\widetilde u(s)\|_{L^2_y}^2
			+
			\|\widetilde P(s)\|_{L^2_y}^2
			\right)
			\le
			\mathcal H(s)
			\le
			C
			\left(
			\|\nabla_y\widetilde u(s)\|_{L^2_y}^2
			+
			\|\widetilde P(s)\|_{L^2_y}^2
			\right).
		\end{equation}
		
		In the normalized variables,
		the velocity estimate in the proof of Lemma \ref{lem3.2} gives
		\[
		\begin{aligned}
			&\frac{d}{ds}
			\int_{\mathbb R^2}
			\left(
			\frac{\mu}{2}|\nabla_y\widetilde u|^2
			+
			\frac{\mu+\lambda}{2}(\div_y\widetilde u)^2
			-
			\widetilde P\div_y\widetilde u
			\right)dy
			+
			\frac12
			\int_{\mathbb R^2}\widetilde\rho|\dot{\widetilde u}|^2\,dy \le
			C_1\|\widetilde P(s)\|_{L^3_y}^3 ,
		\end{aligned}
		\]
		where \(C_1\ge1\) depends only on \(\mu,\lambda,\gamma\).
		
		On the other hand,
		the differential form of the pressure estimate in the proof of Lemma \ref{lem3.3.1} gives
		\[
		\frac{d}{ds}\|\widetilde P(s)\|_{L^2_y}^2
		+
		\frac{2\gamma-1}{2(2\mu+\lambda)}
		\|\widetilde P(s)\|_{L^3_y}^3
		\le
		\frac{1}
		{16C_1
			\left(
			1+\frac{2(2\mu+\lambda)}{2\gamma-1}
			\right)}
		\int_{\mathbb R^2}\widetilde\rho|\dot{\widetilde u}|^2\,dy .
		\]
		Choose
		\[
		C_2
		=
		4C_1
		\left(
		1+\frac{2(2\mu+\lambda)}{2\gamma-1}
		\right).
		\]
		Multiplying the last inequality by \(C_2\) and adding it to the preceding
		velocity estimate, we obtain
		\[
		\begin{aligned}
			\mathcal H'(s)
			&+
			\left[
			\frac12
			-
			\frac{C_2}
			{16C_1
				\left(
				1+\frac{2(2\mu+\lambda)}{2\gamma-1}
				\right)}
			\right]
			\int_{\mathbb R^2}\widetilde\rho|\dot{\widetilde u}|^2\,dy +
			\left[
			\frac{C_2(2\gamma-1)}{2(2\mu+\lambda)}
			-
			C_1
			\right]
			\|\widetilde P(s)\|_{L^3_y}^3
			\le0 .
		\end{aligned}
		\]
		By the choice of \(C_2\), both coefficients in the left-hand side are positive.
		Hence there exists a constant \(c>0\), depending only on
		\(\mu,\lambda,\gamma\), such that
		\begin{equation}\label{H-dissipation}
			\mathcal H'(s)
			+
			c
			\left(
			\int_{\mathbb R^2}\widetilde\rho|\dot{\widetilde u}|^2\,dy
			+
			\|\widetilde P(s)\|_{L^3_y}^3
			\right)
			\le0 .
		\end{equation}
		Multiplying \eqref{H-dissipation} by \(s\), we obtain
		\[
		s\mathcal H'(s)
		+
		cs
		\left(
		\int_{\mathbb R^2}\widetilde\rho|\dot{\widetilde u}|^2\,dy
		+
		\|\widetilde P(s)\|_{L^3_y}^3
		\right)
		\le0.
		\]
		Since
		\[
		s\mathcal H'(s)
		=
		\frac{d}{ds}\bigl(s\mathcal H(s)\bigr)-\mathcal H(s),
		\]
		we get
		\begin{equation}\label{weighted-H-diff}
			\frac{d}{ds}\bigl(s\mathcal H(s)\bigr)
			+
			cs
			\left(
			\int_{\mathbb R^2}\widetilde\rho|\dot{\widetilde u}|^2\,dy
			+
			\|\widetilde P(s)\|_{L^3_y}^3
			\right)
			\le
			\mathcal H(s).
		\end{equation}
		Integrating \eqref{weighted-H-diff} over \([1,\tau]\), where
		\(1\le \tau\le \bar\rho^\gamma T\), yields
		\begin{equation}\label{weighted-H-int}
			\begin{aligned}
				&\tau\mathcal H(\tau)
				+
				c\int_1^\tau
				s
				\left(
				\int_{\mathbb R^2}\widetilde\rho|\dot{\widetilde u}|^2\,dy
				+
				\|\widetilde P(s)\|_{L^3_y}^3
				\right)ds
				\le
				\mathcal H(1)+\int_1^\tau \mathcal H(s)\,ds .
			\end{aligned}
		\end{equation}
		
		We now estimate the right-hand side. By \eqref{H-equivalence} and Corollary
		\ref{cor3.1.1} in the normalized variables,
		\begin{equation}\label{H1-bound}
			\mathcal H(1)
			\le
			C\bar\rho^{1-\gamma}
			\left(
			\|\nabla u_0\|_{L^2}^2+\|P_0\|_{L^2}^2
			\right).
		\end{equation}
		Next, by \eqref{H-equivalence} we have
		\[
		\int_1^\tau\mathcal H(s)\,ds
		\le
		C\int_1^\tau
		\left(
		\|\nabla_y\widetilde u(s)\|_{L^2_y}^2
		+
		\|\widetilde P(s)\|_{L^2_y}^2
		\right)ds.
		\]
		For the first term, \eqref{3.1.1} yields
		\[
		\int_1^\tau
		\|\nabla_y\widetilde u(s)\|_{L^2_y}^2\,ds
		\le
		C\bar\rho E_0.
		\]
		For the pressure term, it follows from momentum equation \eqref{1.1}$_{2}$ that
		\[
		\Delta_y\widetilde F=\div_y(\widetilde\rho\dot{\widetilde u}),
		\]
		where $	\widetilde F:=(2\mu+\lambda)\div_y\widetilde u-\widetilde P$.
		Observe that
		\[
		\widetilde P
		=
		(2\mu+\lambda)\div_y\widetilde u
		+
		(-\Delta_y)^{-1}\div_y(\widetilde\rho\dot{\widetilde u}).
		\]
		Multiplying this identity by \(\widetilde P\) and integrating over \(\mathbb R^2\), we get
		\[
		\|\widetilde P\|_{L^2_y}^2
		\le
		C\|\nabla_y\widetilde u\|_{L^2_y}\|\widetilde P\|_{L^2_y}
		+
		\left|
		\int_{\mathbb R^2}
		(-\Delta_y)^{-1}\div_y(\widetilde\rho\dot{\widetilde u})
		\widetilde P\,dy
		\right|.
		\]
		By the Hardy--Littlewood--Sobolev inequality (see, e.g., \cite{Stein1993}),
		\[
		\left\|
		(-\Delta_y)^{-1}\div_y(\widetilde\rho\dot{\widetilde u})
		\right\|_{L^{4\gamma}_y}
		\le
		C
		\|\widetilde\rho\dot{\widetilde u}\|_{L^{\frac{4\gamma}{2\gamma+1}}_y}.
		\]
		Therefore,
		\[
		\begin{aligned}
			&\left|
			\int_{\mathbb R^2}
			(-\Delta_y)^{-1}\div_y(\widetilde\rho\dot{\widetilde u})
			\widetilde P\,dy
			\right|\le
			C
			\|\widetilde\rho\dot{\widetilde u}\|_{L^{\frac{4\gamma}{2\gamma+1}}_y}
			\|\widetilde P\|_{L^{\frac{4\gamma}{4\gamma-1}}_y}.
		\end{aligned}
		\]
		Since
		\[
		\frac{2\gamma+1}{4\gamma}=\frac12+\frac1{4\gamma},
		\]
		we have
		\[
		\|\widetilde\rho\dot{\widetilde u}\|_{L^{\frac{4\gamma}{2\gamma+1}}_y}
		\le
		\|\widetilde\rho^{1/2}\dot{\widetilde u}\|_{L^2_y}
		\|\widetilde\rho\|_{L^{2\gamma}_y}^{1/2}.
		\]
		Moreover, using \(\widetilde P=\widetilde\rho^\gamma\) and interpolation between
		\(L^1\) and \(L^{2\gamma}\), we get
		\[
		\|\widetilde P\|_{L^{\frac{4\gamma}{4\gamma-1}}_y}
		\le
		\|\widetilde\rho\|_{L^1_y}^{1/2}
		\|\widetilde\rho\|_{L^{2\gamma}_y}^{\gamma-\frac12},
		\]
		which implies that
		\[
		\left|
		\int_{\mathbb R^2}
		(-\Delta_y)^{-1}\div_y(\widetilde\rho\dot{\widetilde u})
		\widetilde P\,dy
		\right|
		\le
		C(\bar\rho^\gamma M_0)^{1/2}
		\|\widetilde\rho^{1/2}\dot{\widetilde u}\|_{L^2_y}
		\|\widetilde P\|_{L^2_y}.
		\]
		Consequently,
		\[
		\|\widetilde P\|_{L^2_y}^2
		\le
		C\|\nabla_y\widetilde u\|_{L^2_y}^2
		+
		C(\bar\rho^\gamma M_0)
		\int_{\mathbb R^2}\widetilde\rho|\dot{\widetilde u}|^2\,dy.
		\]
		It follows that
		\[
		\begin{aligned}
			\int_1^\tau\mathcal H(s)\,ds
			&\le C\int_1^\tau\|\nabla_y\widetilde u(s)\|_{L^2_y}^2\,ds
			+ C(\bar\rho^\gamma M_0)\int_1^\tau \int_{\mathbb R^2}\widetilde\rho|\dot{\widetilde u}|^2\,dyds \\
			&\le C\bar\rho E_0
			+ C(\bar\rho^\gamma M_0)\bar\rho^{1-\gamma}
			\left(\|\nabla u_0\|_{L^2}^2+\|P_0\|_{L^2}^2\right) \\
			&\le C\left[
			\bar\rho E_0
			+
			\left(
			\bar\rho^\gamma M_0
			+
			(\bar\rho^\gamma M_0)^{-1}
			\right)
			\bar\rho^{1-\gamma}
			\left(
			\|\nabla u_0\|_{L^2}^2+\|P_0\|_{L^2}^2
			\right)
			\right].
		\end{aligned}
		\]
		Substituting this and \eqref{H1-bound} into
		\eqref{weighted-H-int}, we obtain
		\begin{equation}\label{large-time-H-decay-final}
			\begin{aligned}
				&\sup_{1\le s\le \tau} s \left(\|\nabla_y\widetilde u(s)\|_{L^2_y}^2 + \|\widetilde P(s)\|_{L^2_y}^2\right)
				+ \int_1^\tau s\left(\int_{\mathbb R^2}\widetilde\rho|\dot{\widetilde u}|^2\,dy + \|\widetilde P(s)\|_{L^3_y}^3\right)ds \\
				\le& C\Bigg[\bar\rho E_0 + \left(\bar\rho^\gamma M_0 + (\bar\rho^\gamma M_0)^{-1}\right)\bar\rho^{1-\gamma}\left(\|\nabla u_0\|_{L^2}^2+\|P_0\|_{L^2}^2\right)\Bigg].
			\end{aligned}
		\end{equation}
		
		It remains to obtain the pointwise decay of
		\(\int\widetilde\rho|\dot{\widetilde u}|^2\).
		Applying \(\partial_s+\widetilde u\cdot\nabla_y\) to \eqref{3.5.6}, multiplying the
		resulting equation by \(\dot{\widetilde u}\) and integrating by parts
		over \(\mathbb R^2\), we obtain
		\begin{equation}\label{dotu-diff}
			\frac{d}{ds}
			\int_{\mathbb R^2}\widetilde\rho|\dot{\widetilde u}|^2\,dy
			+
			c\int_{\mathbb R^2}|\nabla_y\dot{\widetilde u}|^2\,dy
			\le
			C\left(
			\|\nabla_y\widetilde u(s)\|_{L^4_y}^4
			+
			\|\widetilde P(s)\|_{L^4_y}^4
			\right).
		\end{equation}
		By Lemma \ref{lem2.9} with \(p=4\), \(0\le\widetilde\rho\le2\), and
		\(\widetilde P\le2^\gamma\), we have
		\[
		\begin{aligned}
			&\|\nabla_y\widetilde u(s)\|_{L^4_y}^4
			+
			\|\widetilde P(s)\|_{L^4_y}^4\\
			\le{}&
			C
			\left(
			\|\nabla_y\widetilde u(s)\|_{L^2_y}^2
			+
			\|\widetilde P(s)\|_{L^2_y}^2
			\right)
			\int_{\mathbb R^2}\widetilde\rho|\dot{\widetilde u}|^2\,dy
			+
			C\|\widetilde P(s)\|_{L^3_y}^3.
		\end{aligned}
		\]
		Using \eqref{3.1.3}, we infer
		\begin{equation}\label{dotu-diff-2}
			\frac{d}{ds}
			\int_{\mathbb R^2}\widetilde\rho|\dot{\widetilde u}|^2\,dy
			+
			c\int_{\mathbb R^2}|\nabla_y\dot{\widetilde u}|^2\,dy
			\le
			C
			\int_{\mathbb R^2}\widetilde\rho|\dot{\widetilde u}|^2\,dy
			+
			C\|\widetilde P(s)\|_{L^3_y}^3.
		\end{equation}
		Multiplying \eqref{dotu-diff-2} by \(s\), we obtain
		\[
		\begin{aligned}
			&\frac{d}{ds}
			\left[
			s\int_{\mathbb R^2}\widetilde\rho|\dot{\widetilde u}|^2\,dy
			\right]
			+
			cs\int_{\mathbb R^2}|\nabla_y\dot{\widetilde u}|^2\,dy\\
			\le&
			\int_{\mathbb R^2}\widetilde\rho|\dot{\widetilde u}|^2\,dy
			+
			C s
			\int_{\mathbb R^2}\widetilde\rho|\dot{\widetilde u}|^2\,dy
			+
			Cs\|\widetilde P(s)\|_{L^3_y}^3.
		\end{aligned}
		\]
		Integrating over \([1,\tau]\), using Lemma \ref{lem3.5} in normalized variables and
		\eqref{large-time-H-decay-final}, we get
		\begin{equation}\label{dotu-decay-final}
			\begin{aligned}
				&\sup_{1\le s\le \tau}
				s\int_{\mathbb R^2}\widetilde\rho|\dot{\widetilde u}(s)|^2\,dy
				+
				\int_1^\tau
				s\int_{\mathbb R^2}|\nabla_y\dot{\widetilde u}|^2\,dyds\\
				\le&
				C\Bigg[
				\bar\rho E_0
				+
				\left(
				\bar\rho^\gamma M_0
				+
				(\bar\rho^\gamma M_0)^{-1}
				\right)
				\bar\rho^{1-\gamma}
				\left(
				\|\nabla u_0\|_{L^2}^2+\|P_0\|_{L^2}^2
				\right)
				\Bigg].
			\end{aligned}
		\end{equation}
		Combining \eqref{large-time-H-decay-final} and \eqref{dotu-decay-final} gives
		\eqref{large-time-L2-decay-normalized}. Scaling back through
		\eqref{scale-large-L2}--\eqref{scale-large-P3} yields
		\eqref{large-time-L2-decay-original}. The proof is complete.
	\end{proof}

	The a priori $L^{6}$ bound used in three dimensions is unavailable
	and we need new methods to estimate the term $\|\rho|u|^{r}\|_{L^1}$. We shall make full use of \eqref{2.6} but the crucial step to overcome this difficulty is to obtain a positive lower bound for the local mass of the density.
	
	\begin{lemma}\label{lem3.3}
		Let \((\rho,u)\) solve \eqref{1.1}--\eqref{1.4.1} smoothly on
		\(\mathbb R^2\times(0,T]\) satisfying \eqref{3.1.3} and the assumption of
		Theorem \ref{th1.1}. Then there exists a positive constant
		\begin{equation}\label{3.26.1-new}
			N_1
			:=
			4\left[
			1
			+
			\left(
			\frac{\displaystyle\int_{\mathbb R^2}\rho_0(x)|x|^\alpha\,dx}{M_0}
			\right)^{1/\alpha}
			+
			\left(
			\frac{2E_0}{M_0}
			\right)^{1/2}
			\right],
		\end{equation}
		such that
		\begin{equation}\label{3.26.2-new}
			\inf_{0\le s\le T}
			\int_{B_{N_1(1+s)}}\rho(x,s)\,dx
			\ge
			\frac{M_0}{2}.
		\end{equation}
	\end{lemma}
	
	\begin{proof}
		Multiplying the continuity equation \eqref{1.1}$_1$ by
		\((1+|x|^2)^{1/2}\)
		and integrating over \(\mathbb R^2\), we obtain, after integration by parts,
		\[
		\frac{d}{dt}
		\int_{\mathbb R^2}
		\rho(x,t)(1+|x|^2)^{1/2}\,dx
		=
		\int_{\mathbb R^2}
		\rho u\cdot\nabla(1+|x|^2)^{1/2}\,dx.
		\]
		Since
		\[
		\left|\nabla(1+|x|^2)^{1/2}\right|
		=
		\frac{|x|}{(1+|x|^2)^{1/2}}
		\le 1,
		\]
		we have
		\[
		\frac{d}{dt}
		\int_{\mathbb R^2}
		\rho(x,t)(1+|x|^2)^{1/2}\,dx
		\le
		\int_{\mathbb R^2}\rho |u|\,dx.
		\]
		By mass conservation,
		\[
		\int_{\mathbb R^2}\rho(x,t)\,dx=M_0,
		\qquad 0\le t\le T,
		\]
		and by  \eqref{3.1.1} we have
		\[
		\int_{\mathbb R^2}\rho|u|^2\,dx\le 2E_0.
		\]
		Thus Cauchy's inequality gives
		\[
		\int_{\mathbb R^2}\rho|u|\,dx
		\le
		\left(
		\int_{\mathbb R^2}\rho\,dx
		\right)^{1/2}
		\left(
		\int_{\mathbb R^2}\rho|u|^2\,dx
		\right)^{1/2}
		\le
		\sqrt{2M_0E_0}.
		\]
		Therefore,
		\begin{equation}\label{3.28.1-new}
			\frac{d}{dt}
			\int_{\mathbb R^2}
			\rho(x,t)(1+|x|^2)^{1/2}\,dx
			\le
			\sqrt{2M_0E_0}.
		\end{equation}
		
		Next we estimate the initial first moment. Since
		\[
		(1+|x|^2)^{1/2}\le 1+|x|,
		\]
		we have
		\[
		\int_{\mathbb R^2}
		\rho_0(x)(1+|x|^2)^{1/2}\,dx
		\le
		M_0+\int_{\mathbb R^2}\rho_0(x)|x|\,dx.
		\]
		Using H\"older's inequality and \(\alpha>1\), we obtain
		\[
		\begin{aligned}
			\int_{\mathbb R^2}\rho_0(x)|x|\,dx
			&=
			\int_{\mathbb R^2}
			\rho_0(x)^{1-\frac1\alpha}
			\left(\rho_0(x)|x|^\alpha\right)^{1/\alpha}\,dx\\
			&\le
			\left(
			\int_{\mathbb R^2}\rho_0(x)\,dx
			\right)^{1-\frac1\alpha}
			\left(
			\int_{\mathbb R^2}\rho_0(x)|x|^\alpha\,dx
			\right)^{1/\alpha}\\
			&=
			M_0^{1-\frac1\alpha}
			\left(
			\int_{\mathbb R^2}\rho_0(x)|x|^\alpha\,dx
			\right)^{1/\alpha}.
		\end{aligned}
		\]
		Consequently,
		\begin{equation}\label{3.28.2-new-initial}
			\int_{\mathbb R^2}
			\rho_0(x)(1+|x|^2)^{1/2}\,dx
			\le
			M_0
			+
			M_0^{1-\frac1\alpha}
			\left(
			\int_{\mathbb R^2}\rho_0(x)|x|^\alpha\,dx
			\right)^{1/\alpha}.
		\end{equation}
		
		Integrating \eqref{3.28.1-new} over \([0,t]\), and using
		\eqref{3.28.2-new-initial}, we get
		\begin{equation}\label{3.28.2-new}
			\begin{aligned}
				&\sup_{0\le \tau\le t}
				\int_{\mathbb R^2}
				\rho(x,\tau)(1+|x|^2)^{1/2}\,dx\\
				\le&
				M_0
				+
				M_0^{1-\frac1\alpha}
				\left(
				\int_{\mathbb R^2}\rho_0(x)|x|^\alpha\,dx
				\right)^{1/\alpha}
				+
				\sqrt{2M_0E_0}\,t .
			\end{aligned}
		\end{equation}
		
		Now fix any \(s\in[0,T]\). For \(|x|\ge N_1(1+s)\), one has
		\[
		(1+|x|^2)^{1/2}\ge |x|\ge N_1(1+s).
		\]
		Hence
		\[
		\begin{aligned}
			\int_{|x|\ge N_1(1+s)}\rho(x,s)\,dx
			&\le
			\frac{1}{N_1(1+s)}
			\int_{|x|\ge N_1(1+s)}
			\rho(x,s)(1+|x|^2)^{1/2}\,dx\\
			&\le
			\frac{1}{N_1(1+s)}
			\int_{\mathbb R^2}
			\rho(x,s)(1+|x|^2)^{1/2}\,dx.
		\end{aligned}
		\]
		Using \eqref{3.28.2-new}, we obtain
		\[
		\begin{aligned}
			\int_{|x|\ge N_1(1+s)}\rho(x,s)\,dx
			&\le
			\frac{
				M_0
				+
				M_0^{1-\frac1\alpha}
				\left(
				\int_{\mathbb R^2}\rho_0|x|^\alpha\,dx
				\right)^{1/\alpha}
				+
				\sqrt{2M_0E_0}\,s
			}
			{N_1(1+s)}\\
			&=
			\frac{M_0}{N_1(1+s)}
			\left[
			1
			+
			\left(
			\frac{\displaystyle\int_{\mathbb R^2}\rho_0|x|^\alpha\,dx}{M_0}
			\right)^{1/\alpha}
			+
			\left(
			\frac{2E_0}{M_0}
			\right)^{1/2}s
			\right].
		\end{aligned}
		\]
		Since \(s\le1+s\), one has
		\[
		\begin{aligned}
			&1
			+
			\left(
			\frac{\displaystyle\int_{\mathbb R^2}\rho_0|x|^\alpha\,dx}{M_0}
			\right)^{1/\alpha}
			+
			\left(
			\frac{2E_0}{M_0}
			\right)^{1/2}s\le
			\left[
			1
			+
			\left(
			\frac{\displaystyle\int_{\mathbb R^2}\rho_0|x|^\alpha\,dx}{M_0}
			\right)^{1/\alpha}
			+
			\left(
			\frac{2E_0}{M_0}
			\right)^{1/2}
			\right](1+s).
		\end{aligned}
		\]
		By the definition of \(N_1\), we get
		\[
		\int_{|x|\ge N_1(1+s)}\rho(x,s)\,dx
		\le
		\frac{M_0}{4}.
		\]
		Hence we have
		\[
		\begin{aligned}
			\int_{B_{N_1(1+s)}}\rho(x,s)\,dx
			&=
			\int_{\mathbb R^2}\rho(x,s)\,dx
			-
			\int_{|x|\ge N_1(1+s)}\rho(x,s)\,dx\\
			&\ge
			M_0-\frac{M_0}{4}
			=
			\frac{3M_0}{4}
			\ge
			\frac{M_0}{2}.
		\end{aligned}
		\]
		Taking the infimum over \(0\le s\le T\), we obtain \eqref{3.26.2-new}.
		This proves the lemma.
	\end{proof}

	\begin{lemma}\label{lem-rho-xbara-zlotnik}
		Let \((\rho,u)\) solve \eqref{1.1}--\eqref{1.4.1} smoothly on
		\(\mathbb R^2\times(0,T]\) satisfying \eqref{3.1.3} and the assumptions of
		Theorem \ref{th1.1}. Then, for any \(t\in[0,T]\),
		\begin{equation}\label{rho-xbara-zlotnik}
			\begin{aligned}
				&\sup_{0\le s\le t}
				\bar\rho^\gamma
				\|\rho(\cdot,s)\bar x^{1/2}\|_{L^1}\\
				\le{}&
				C(\alpha)
				\left(
				\bar\rho^\gamma M_0
				+
				(\bar\rho^\gamma M_0)^{-1}
				+
				\bar\rho^{\gamma+\frac{\alpha(\gamma+1)}2}
				\|\rho_0|x|^\alpha\|_{L^1}
				\right)
				(1+\bar\rho E_0)
				(1+\bar\rho^\gamma t).
			\end{aligned}
		\end{equation}
	\end{lemma}
	
	\begin{proof}
		We normalize the density upper bound.
		For \(0\le \tau\le \bar\rho^\gamma t\), define
		\[
		\widetilde Y(\tau)
		:=
		\int_{\mathbb R^2}\widetilde\rho(y,\tau)\,\overline y^{1/2}\,dy,
		\qquad
		\overline y
		:=
		(e+|y|^2)^{1/2}\log^2(e+|y|^2).
		\]
		Multiplying the normalized continuity equation \eqref{1.1}$_{1}$
		\[
		\widetilde\rho_\tau+\div_y(\widetilde\rho\widetilde u)=0
		\]
		by \(\overline y^{1/2}\), and integrating by parts over \(\mathbb R^2\), we get
		\begin{equation}\label{rho-xbara-z-p1}
			\widetilde Y'(\tau)
			=
			\int_{\mathbb R^2}
			\widetilde\rho\,\widetilde u\cdot\nabla \overline y^{1/2}\,dy.
		\end{equation}
		Since
		\[
		|\nabla \overline y^{1/2}|\le C,
		\]
		we obtain from \eqref{rho-xbara-z-p1} that
		\[
		\begin{aligned}
			\widetilde Y'(\tau)
			&\le
			C
			\int_{\mathbb R^2}
			\widetilde\rho|\widetilde u|\,dy\le
			C
			\left(
			\int_{\mathbb R^2}\widetilde\rho\,dy
			\right)^{1/2}
			\left(
			\int_{\mathbb R^2}\widetilde\rho|\widetilde u|^2\,dy
			\right)^{1/2}.
		\end{aligned}
		\]
		By mass conservation,
		\[
		\int_{\mathbb R^2}\widetilde\rho(y,\tau)\,dy
		=
		\int_{\mathbb R^2}\widetilde\rho_0(y)\,dy
		=
		\bar\rho^\gamma M_0,
		\]
		and \eqref{3.1.1} in the normalized variables,
		\[
		\sup_{0\le \tau\le \bar\rho^\gamma T}
		\int_{\mathbb R^2}\widetilde\rho|\widetilde u|^2\,dy
		\le
		C\bar\rho E_0,
		\]
		we have
		\begin{equation}\label{rho-xbara-z-diff}
			\widetilde Y'(\tau)
			\le
			C
			(\bar\rho^\gamma M_0)^{1/2}
			(\bar\rho E_0)^{1/2}.
		\end{equation}
		Integrating \eqref{rho-xbara-z-diff} over \([0,\tau]\), we get
		\begin{equation}\label{rho-xbara-z-Y-bound}
			\widetilde Y(\tau)
			\le
			\widetilde Y(0)
			+
			C
			(\bar\rho^\gamma M_0)^{1/2}
			(\bar\rho E_0)^{1/2}\tau .
		\end{equation}
		
		We estimate \(\widetilde Y(0)\). Since \(\alpha>1\), one has
		\[
		\overline y^{1/2}
		\le
		C(\alpha)(1+|y|^\alpha).
		\]
		Hence
		\begin{equation}\label{rho-xbara-z-Y0}
			\begin{aligned}
				\widetilde Y(0)
				&=
				\int_{\mathbb R^2}\widetilde\rho_0(y)\overline y^{1/2}\,dy\\
				&\le
				C(\alpha)
				\left(
				\int_{\mathbb R^2}\widetilde\rho_0(y)\,dy
				+
				\int_{\mathbb R^2}\widetilde\rho_0(y)|y|^\alpha\,dy
				\right)\\
				&=
				C(\alpha)
				\left[
				\bar\rho^\gamma M_0
				+
				\bar\rho^{\gamma+\frac{\alpha(\gamma+1)}2}
				\|\rho_0|x|^\alpha\|_{L^1}
				\right].
			\end{aligned}
		\end{equation}
		Combining \eqref{rho-xbara-z-Y-bound} and \eqref{rho-xbara-z-Y0}, we get
		\[
		\begin{aligned}
			\sup_{0\le r\le \tau}\widetilde Y(r)
			\le{}&
			C(\alpha)
			\left[
			\bar\rho^\gamma M_0
			+
			\bar\rho^{\gamma+\frac{\alpha(\gamma+1)}2}
			\|\rho_0|x|^\alpha\|_{L^1}
			\right]\\
			&+
			C
			(\bar\rho^\gamma M_0)^{1/2}
			(\bar\rho E_0)^{1/2}(1+\tau).
		\end{aligned}
		\]
		Using
		\[
		(\bar\rho^\gamma M_0)^{1/2}
		\le
		\bar\rho^\gamma M_0+(\bar\rho^\gamma M_0)^{-1},
		\qquad
		(\bar\rho E_0)^{1/2}
		\le
		1+\bar\rho E_0,
		\]
		we infer that
		\begin{equation}\label{rho-xbara-z-normal}
			\begin{aligned}
				\sup_{0\le r\le \tau}\widetilde Y(r)
				\le{}&
				C(\alpha)
				\left(
				\bar\rho^\gamma M_0
				+
				(\bar\rho^\gamma M_0)^{-1}
				+
				\bar\rho^{\gamma+\frac{\alpha(\gamma+1)}2}
				\|\rho_0|x|^\alpha\|_{L^1}
				\right)\\
				&\times
				(1+\bar\rho E_0)(1+\tau).
			\end{aligned}
		\end{equation}
		
		Finally, we return to the original variables. Since
		\[
		y=\bar\rho^{\frac{\gamma+1}{2}}x,
		\qquad
		\tau=\bar\rho^\gamma s,
		\]
		we have
		\[
		\begin{aligned}
			\widetilde Y(\tau)
			&=
			\int_{\mathbb R^2}
			\widetilde\rho(y,\tau)\overline y^{1/2}\,dy\\
			&=
			\bar\rho^\gamma
			\int_{\mathbb R^2}
			\rho(x,s)
			\left(e+\bar\rho^{\gamma+1}|x|^2\right)^{1/4}
			\log\!\left(e+\bar\rho^{\gamma+1}|x|^2\right)\,dx.
		\end{aligned}
		\]
		Since \(\bar\rho\ge1\), one has
		\[
		\left(e+\bar\rho^{\gamma+1}|x|^2\right)^{1/4}
		\log\!\left(e+\bar\rho^{\gamma+1}|x|^2\right)
		\ge
		\bar x^{1/2}.
		\]
		Hence
		\[
		\bar\rho^\gamma
		\|\rho(\cdot,s)\bar x^{1/2}\|_{L^1}
		\le
		\widetilde Y(\tau),
		\qquad
		\tau=\bar\rho^\gamma s.
		\]
		Taking \(\tau=\bar\rho^\gamma t\) in \eqref{rho-xbara-z-normal}, we obtain
		\[
		\begin{aligned}
			\sup_{0\le s\le t}
			\bar\rho^\gamma
			\|\rho(\cdot,s)\bar x^{1/2}\|_{L^1}
			\le{}&
			C(\alpha)
			\left(
			\bar\rho^\gamma M_0
			+
			(\bar\rho^\gamma M_0)^{-1}
			+
			\bar\rho^{\gamma+\frac{\alpha(\gamma+1)}2}
			\|\rho_0|x|^\alpha\|_{L^1}
			\right)\\
			&\times
			(1+\bar\rho E_0)(1+\bar\rho^\gamma t).
		\end{aligned}
		\]
		This proves \eqref{rho-xbara-zlotnik}. The proof is complete.
	\end{proof}
	
	The next lemma establishes the uniform-in-time upper bound for the density
	required in Proposition \ref{p3.1}.
	\begin{lemma}\label{lem:uniform-density-bound}
		Let \((\rho,u)\) solve \eqref{1.1}--\eqref{1.4.1} smoothly on
		\(\mathbb R^2\times(0,T]\) satisfying \eqref{3.1.3} and the assumptions of
		Theorem \ref{th1.1}.
		Then
		\begin{equation}\label{uniform-density-bound}
			0\le \rho(x,t)\le \frac32\bar\rho,
			\qquad
			(x,t)\in\mathbb R^2\times[0,T].
		\end{equation}
	\end{lemma}
	
	\begin{proof}
		We first normalize the density upper bound.
		By \eqref{3.1.3}, we have
		\begin{equation}\label{zlotnik-scale-smallness}
			\begin{aligned}
				&\Xi_0(1+\bar\rho E_0)^{3/2}
				\bar\rho^{1-\gamma}
				\sup_{0\le t\le T}
				\left(
				\|\nabla u(t)\|_{L^2}^2+\|P(t)\|_{L^2}^2
				\right)
				\le 2\varepsilon .
			\end{aligned}
		\end{equation}
		In particular, since \(\Xi_0\ge1\) and \(1+\bar\rho E_0\ge1\),
		\begin{equation}\label{zlotnik-ST-small}
			\bar\rho^{1-\gamma}
			\sup_{0\le t\le T}
			\left(
			\|\nabla u(t)\|_{L^2}^2+\|P(t)\|_{L^2}^2
			\right)
			\le 2\varepsilon.
		\end{equation}
		We choose \(\varepsilon>0\) sufficiently small such that
		\begin{equation}\label{zlotnik-eps-small}
			C\varepsilon^{1/2}\le \frac14,
			\qquad
			C\varepsilon^{6/7}\le \frac14,
		\end{equation}
		where \(C\) denotes the constants appearing in \eqref{short-F-L1-Linfty-small} and \eqref{large-F-Linfty-bound-new}.
		
		We first record the estimate for \(\|\nabla_y\widetilde F\|_{L^{40}_y}\), which will be
		used both in the short-time and large-time estimates. It follows from Lemma \ref{lem2.9} that
		\[
		\|\nabla_y\widetilde F(\tau)\|_{L^{40}_y}
		\le
		C\|\widetilde\rho\dot{\widetilde u}(\tau)\|_{L^{40}_y}.
		\]
		By Lemma \ref{lem3.3} in the normalized variables,
		\[
		\inf_{0\le s\le \tau}
		\int_{B_{\widetilde N_1(1+s)}}\widetilde\rho(y,s)\,dy
		\ge
		\frac{\bar\rho^\gamma M_0}{2},
		\]
		where
		\[
		\widetilde N_1
		=
		4\left[
		1
		+
		\left(
		\frac{
			\bar\rho^{\gamma+\frac{\alpha(\gamma+1)}2}
			\|\rho_0|x|^\alpha\|_{L^1}
		}{
			\bar\rho^\gamma M_0
		}
		\right)^{1/\alpha}
		+
		\left(
		\frac{2\bar\rho E_0}{\bar\rho^\gamma M_0}
		\right)^{1/2}
		\right].
		\]
		We claim that
		\begin{equation}\label{N1-control-new}
			\widetilde N_1
			\left(1+(\bar\rho^\gamma M_0)^{-1/2}\right)
			\le
			C(\alpha)\Xi_0^3(1+\bar\rho E_0)^{1/2}.
		\end{equation}
		Indeed, since \(\Xi_0\ge1\),
		\(\Xi_0\ge(\bar\rho^\gamma M_0)^{-1}\), and
		\(\Xi_0\ge
		\bar\rho^{\gamma+\frac{\alpha(\gamma+1)}2}
		\|\rho_0|x|^\alpha\|_{L^1}\),
		we have
		\[
		\left(
		\frac{
			\bar\rho^{\gamma+\frac{\alpha(\gamma+1)}2}
			\|\rho_0|x|^\alpha\|_{L^1}
		}{
			\bar\rho^\gamma M_0
		}
		\right)^{1/\alpha}
		\left(1+(\bar\rho^\gamma M_0)^{-1/2}\right)
		\le
		C(\alpha)\Xi_0^3.
		\]
		Similarly,
		\[
		\left(
		\frac{2\bar\rho E_0}{\bar\rho^\gamma M_0}
		\right)^{1/2}
		\left(1+(\bar\rho^\gamma M_0)^{-1/2}\right)
		\le
		C \Xi_0(1+\bar\rho E_0)^{1/2}
		\le
		C \Xi_0^3(1+\bar\rho E_0)^{1/2}.
		\]
		This proves \eqref{N1-control-new}.
		
		Now apply Lemma \ref{lem2.4.1} with
		\[
		a=\frac12,\qquad r=40,\qquad v=\dot{\widetilde u}.
		\]
		Since \(0\le\widetilde\rho\le2\), \(m_1=2\), \(m_2=\bar\rho^\gamma M_0/2\), and
		\(N_*=\widetilde N_1(1+\tau)\), we get
		\[
		\begin{aligned}
			\|\widetilde\rho\dot{\widetilde u}(\tau)\|_{L^{40}_y}
			\le{}&
			C(\alpha)\widetilde N_1(1+\tau)
			\left(1+(\bar\rho^\gamma M_0)^{-1/2}\right)\\
			&\times
			\left(1+\|\widetilde\rho(\tau)\overline y^{1/2}\|_{L^1_y}\right)^{1/45}\\
			&\times
			\left(
			\|\widetilde\rho^{1/2}\dot{\widetilde u}(\tau)\|_{L^2_y}
			+
			\|\nabla_y\dot{\widetilde u}(\tau)\|_{L^2_y}
			\right).
		\end{aligned}
		\]
		By Lemma \ref{lem-rho-xbara-zlotnik},
		\[
		1+\|\widetilde\rho(\tau)\overline y^{1/2}\|_{L^1_y}
		\le
		C(\alpha)\Xi_0(1+\bar\rho E_0)(1+\tau).
		\]
		Combining this with \eqref{N1-control-new}, we obtain
		\begin{equation}\label{gradF-L40-common}
			\begin{aligned}
				\|\nabla_y\widetilde F(\tau)\|_{L^{40}_y}
				\le{}&
				C(\alpha)
				\Xi_0^{3+\frac1{45}}
				(1+\bar\rho E_0)^{\frac12+\frac1{45}}
				(1+\tau)^{1+\frac1{45}}\\
				&\times
				\left(
				\|\widetilde\rho^{1/2}\dot{\widetilde u}(\tau)\|_{L^2_y}
				+
				\|\nabla_y\dot{\widetilde u}(\tau)\|_{L^2_y}
				\right).
			\end{aligned}
		\end{equation}
		
		We now establish the short-time \(L^1_\tau L^\infty_y\)-estimate of
		\(\widetilde F\). From
		\[
		\widetilde F=(2\mu+\lambda)\div_y\widetilde u-\widetilde P
		\]
		and \eqref{zlotnik-ST-small}, we have
		\begin{equation}\label{short-F-L2}
			\|\widetilde F(\tau)\|_{L^2_y}
			\le
			C
			\left[
			\bar\rho^{1-\gamma}
			\sup_{0\le t\le T}
			\left(
			\|\nabla u(t)\|_{L^2}^2+\|P(t)\|_{L^2}^2
			\right)
			\right]^{1/2}.
		\end{equation}
		Moreover,
		\[
		\|\nabla_y\widetilde F(\tau)\|_{L^2_y}
		\le
		C\|\widetilde\rho^{1/2}\dot{\widetilde u}(\tau)\|_{L^2_y}.
		\]
		By Lemma \ref{lem3.5} in the normalized variables,
		there holds
		\begin{equation}\label{short-sigma-dotu-used}
			\begin{aligned}
				&\sup_{0\le \tau\le \bar\rho^\gamma T}
				\widetilde\sigma(\tau)^2
				\int_{\mathbb R^2}
				\widetilde\rho|\dot{\widetilde u}|^2\,dy
				+
				\int_0^{\bar\rho^\gamma T}
				\widetilde\sigma(\tau)^2
				\int_{\mathbb R^2}
				|\nabla_y\dot{\widetilde u}|^2\,dy\,d\tau\\
				\le&
				C
				\bar\rho^{1-\gamma}
				\left(
				\|\nabla u_0\|_{L^2}^2+\|P_0\|_{L^2}^2
				\right)
				\le
				C
				\bar\rho^{1-\gamma}
				\sup_{0\le t\le T}
				\left(
				\|\nabla u(t)\|_{L^2}^2+\|P(t)\|_{L^2}^2
				\right).
			\end{aligned}
		\end{equation}
		For \(0<\tau\le1\), \(\widetilde\sigma(\tau)=\tau\). Therefore, by the Gagliardo--Nirenberg inequality \eqref{2.2.1},
		\[
		\begin{aligned}
			\|\widetilde F(\tau)\|_{L^{200}_y}
			&\le
			C
			\|\widetilde F(\tau)\|_{L^2_y}^{1/100}
			\|\nabla_y\widetilde F(\tau)\|_{L^2_y}^{99/100}\\
			&\le
			C
			\left[
			\bar\rho^{1-\gamma}
			\sup_{0\le t\le T}
			\left(
			\|\nabla u(t)\|_{L^2}^2+\|P(t)\|_{L^2}^2
			\right)
			\right]^{1/2}
			\tau^{-99/100}.
		\end{aligned}
		\]
		For \(0<\tau\le1\), it follows from \eqref{gradF-L40-common} and \(1+\tau\le2\) that
		\begin{equation}\label{short-gradF-L40}
			\begin{aligned}
				\|\nabla_y\widetilde F(\tau)\|_{L^{40}_y}
				\le&
				C(\alpha)
				\Xi_0^{3+\frac1{45}}
				(1+\bar\rho E_0)^{\frac12+\frac1{45}}\times
				\left(
				\|\widetilde\rho^{1/2}\dot{\widetilde u}(\tau)\|_{L^2_y}
				+
				\|\nabla_y\dot{\widetilde u}(\tau)\|_{L^2_y}
				\right).
			\end{aligned}
		\end{equation}
		Now the Gagliardo--Nirenberg inequality \eqref{2.3.1} with \(q=200\) and \(r=40\) gives
		\[
		\|\widetilde F\|_{L^\infty_y}
		\le
		C
		\|\widetilde F\|_{L^{200}_y}^{1-\theta}
		\|\nabla_y\widetilde F\|_{L^{40}_y}^{\theta},
		\qquad
		\theta=
		\frac{2r}{2r+q(r-2)}
		=
		\frac1{96}.
		\]
		Therefore,
		\begin{equation}\label{short-F-Linfty-pointwise}
			\begin{aligned}
				\|\widetilde F(\tau)\|_{L^\infty_y}
				\le{}&
				C(\alpha)
				\Xi_0^{\frac{17}{540}}
				(1+\bar\rho E_0)^{\frac{47}{8640}}\\
				&\times
				\left[
				\bar\rho^{1-\gamma}
				\sup_{0\le t\le T}
				\left(
				\|\nabla u(t)\|_{L^2}^2+\|P(t)\|_{L^2}^2
				\right)
				\right]^{\frac{95}{192}}\\
				&\times
				\tau^{-\frac{1881}{1920}}
				\left(
				\|\widetilde\rho^{1/2}\dot{\widetilde u}(\tau)\|_{L^2_y}
				+
				\|\nabla_y\dot{\widetilde u}(\tau)\|_{L^2_y}
				\right)^{\frac1{96}}.
			\end{aligned}
		\end{equation}
		Integrating \eqref{short-F-Linfty-pointwise} over \([0,1]\), and using
		\[
		\left(
		\|\widetilde\rho^{1/2}\dot{\widetilde u}\|_{L^2_y}
		+
		\|\nabla_y\dot{\widetilde u}\|_{L^2_y}
		\right)^{\frac1{96}}
		=
		\tau^{-\frac1{96}}
		\left[
		\tau^2
		\left(
		\|\widetilde\rho^{1/2}\dot{\widetilde u}\|_{L^2_y}
		+
		\|\nabla_y\dot{\widetilde u}\|_{L^2_y}
		\right)^2
		\right]^{\frac1{192}},
		\]
		we obtain by H\"older's inequality with exponents \(192/191\) and \(192\) that
		\[
		\begin{aligned}
			&\int_0^1
			\tau^{-\frac{1881}{1920}}
			\left(
			\|\widetilde\rho^{1/2}\dot{\widetilde u}\|_{L^2_y}
			+
			\|\nabla_y\dot{\widetilde u}\|_{L^2_y}
			\right)^{\frac1{96}}d\tau\\
			\le&
			\left(
			\int_0^1
			\tau^{-\frac{1901}{1920}\cdot\frac{192}{191}}
			d\tau
			\right)^{\frac{191}{192}}\times
			\left(
			\int_0^1
			\tau^2
			\left[
			\|\widetilde\rho^{1/2}\dot{\widetilde u}\|_{L^2_y}
			+
			\|\nabla_y\dot{\widetilde u}\|_{L^2_y}
			\right]^2d\tau
			\right)^{\frac1{192}}.
		\end{aligned}
		\]
		The first integral is finite since
		\[
		\frac{1901}{1920}\cdot\frac{192}{191}
		=
		\frac{1901}{1910}
		<1.
		\]
		Using \eqref{short-sigma-dotu-used}, we infer that
		\[
		\begin{aligned}
			&\int_0^1
			\|\widetilde F(\tau)\|_{L^\infty_y}\,d\tau
			\le
			C(\alpha)
			\Xi_0^{\frac{17}{540}}
			(1+\bar\rho E_0)^{\frac{47}{8640}}
			\left[
			\bar\rho^{1-\gamma}
			\sup_{0\le t\le T}
			\left(
			\|\nabla u(t)\|_{L^2}^2+\|P(t)\|_{L^2}^2
			\right)
			\right]^{\frac12}.
		\end{aligned}
		\]
		Since \(17/540<1/2\), \(47/8640<3/4\), \(\Xi_0\ge1\), and \(1+\bar\rho E_0\ge1\), it follows that
		\[
		\begin{aligned}
			&\int_0^1
			\|\widetilde F(\tau)\|_{L^\infty_y}\,d\tau
			\le
			C(\alpha)
			\left[
			\Xi_0(1+\bar\rho E_0)^{3/2}
			\bar\rho^{1-\gamma}
			\sup_{0\le t\le T}
			\left(
			\|\nabla u(t)\|_{L^2}^2+\|P(t)\|_{L^2}^2
			\right)
			\right]^{1/2}.
		\end{aligned}
		\]
		Combining this with \eqref{zlotnik-scale-smallness}, we get
		\begin{equation}\label{short-F-L1-Linfty-small}
			\int_0^1
			\|\widetilde F(\tau)\|_{L^\infty_y}\,d\tau
			\le
			C\varepsilon^{1/2}.
		\end{equation}
		
		Next we establish the large-time estimate of \(\widetilde F\). In normalized variables,
		Lemma \ref{lem:large-time-L2-decay} gives
		\begin{equation}\label{large-L2-decay-used}
			\begin{aligned}
				&\sup_{1\le \tau\le \bar\rho^\gamma T}
				\tau
				\left(
				\|\nabla_y\widetilde u(\tau)\|_{L^2_y}^2
				+
				\|\widetilde P(\tau)\|_{L^2_y}^2
				+
				\int_{\mathbb R^2}\widetilde\rho|\dot{\widetilde u}|^2\,dy
				\right)\\
				&+
				\int_1^{\bar\rho^\gamma T}
				\tau
				\int_{\mathbb R^2}|\nabla_y\dot{\widetilde u}|^2\,dy\,d\tau
				+
				\int_1^{\bar\rho^\gamma T}
				\tau
				\|\widetilde P(\tau)\|_{L^3_y}^3\,d\tau\\
				&\le
				C\left[
				\bar\rho E_0
				+
				\left(
				\bar\rho^\gamma M_0
				+
				(\bar\rho^\gamma M_0)^{-1}
				\right)
				\bar\rho^{1-\gamma}
				\left(
				\|\nabla u_0\|_{L^2}^2+\|P_0\|_{L^2}^2
				\right)
				\right].
			\end{aligned}
		\end{equation}
		Using \eqref{zlotnik-scale-smallness}, \(\Xi_0\ge
		\bar\rho^\gamma M_0+(\bar\rho^\gamma M_0)^{-1}\), and the fact that the initial value
		is bounded by the supremum in \eqref{zlotnik-ST-small}, we obtain
		\[
		\left(
		\bar\rho^\gamma M_0
		+
		(\bar\rho^\gamma M_0)^{-1}
		\right)
		\bar\rho^{1-\gamma}
		\left(
		\|\nabla u_0\|_{L^2}^2+\|P_0\|_{L^2}^2
		\right)
		\le C.
		\]
		Thus \eqref{large-L2-decay-used} implies
		\begin{equation}\label{large-L2-decay-simplified}
			\begin{aligned}
				&\sup_{1\le \tau\le \bar\rho^\gamma T}
				\tau
				\left(
				\|\nabla_y\widetilde u(\tau)\|_{L^2_y}^2
				+
				\|\widetilde P(\tau)\|_{L^2_y}^2
				+
				\int_{\mathbb R^2}\widetilde\rho|\dot{\widetilde u}|^2\,dy
				\right)\\
				&+
				\int_1^{\bar\rho^\gamma T}
				\tau
				\int_{\mathbb R^2}|\nabla_y\dot{\widetilde u}|^2\,dy\,d\tau
				+
				\int_1^{\bar\rho^\gamma T}
				\tau
				\|\widetilde P(\tau)\|_{L^3_y}^3\,d\tau
				\le
				C(1+\bar\rho E_0).
			\end{aligned}
		\end{equation}
		By Lemma \ref{lem2.9},
		\[
		\|\widetilde F(\tau)\|_{L^2_y}
		+
		\|\nabla_y\widetilde F(\tau)\|_{L^2_y}
		\le
		C\left(
		\|\nabla_y\widetilde u(\tau)\|_{L^2_y}
		+
		\|\widetilde P(\tau)\|_{L^2_y}
		+
		\|\widetilde\rho^{1/2}\dot{\widetilde u}(\tau)\|_{L^2_y}
		\right).
		\]
		By \eqref{large-L2-decay-simplified},
		\begin{equation}\label{F-L40-decay}
			\|\widetilde F(\tau)\|_{L^{40}_y}
			\le
			C(1+\bar\rho E_0)^{1/2}\tau^{-1/2},
			\qquad 1\le\tau\le\bar\rho^\gamma T.
		\end{equation}
		On the other hand, by \eqref{zlotnik-ST-small}, Corollary \ref{cor3.1.1}, and
		Lemma \ref{lem3.5} in normalized variables,
		\begin{equation}\label{F-L40-small}
			\|\widetilde F(\tau)\|_{L^{40}_y}
			\le
			C\left[
			\bar\rho^{1-\gamma}
			\sup_{0\le t\le T}
			\left(
			\|\nabla u(t)\|_{L^2}^2+\|P(t)\|_{L^2}^2
			\right)
			\right]^{1/2}.
		\end{equation}
		Interpolating \eqref{F-L40-decay} and \eqref{F-L40-small} with weights \(7/12\) and \(5/12\), we obtain
		\begin{equation}\label{F-L40-interpolated}
			\begin{aligned}
				\|\widetilde F(\tau)\|_{L^{40}_y}
				\le{}&
				C(1+\bar\rho E_0)^{7/24}\tau^{-7/24}\times
				\left[
				\bar\rho^{1-\gamma}
				\sup_{0\le t\le T}
				\left(
				\|\nabla u(t)\|_{L^2}^2+\|P(t)\|_{L^2}^2
				\right)
				\right]^{5/24}.
			\end{aligned}
		\end{equation}
		For \(\tau\ge1\), \eqref{gradF-L40-common} gives
		\begin{equation}\label{gradF-L40-large}
			\begin{aligned}
				\|\nabla_y\widetilde F(\tau)\|_{L^{40}_y}
				\le{}&
				C(\alpha)
				\Xi_0^{3+\frac1{45}}
				(1+\bar\rho E_0)^{\frac12+\frac1{45}}
				(1+\tau)^{1+\frac1{45}}\\
				&\times
				\left(
				\|\widetilde\rho^{1/2}\dot{\widetilde u}(\tau)\|_{L^2_y}
				+
				\|\nabla_y\dot{\widetilde u}(\tau)\|_{L^2_y}
				\right).
			\end{aligned}
		\end{equation}
		By the Gagliardo--Nirenberg inequality \eqref{2.3.1} with \(q=r=40\),
		\[
		\|\widetilde F\|_{L^\infty_y}^4
		\le
		C
		\|\widetilde F\|_{L^{40}_y}^{19/5}
		\|\nabla_y\widetilde F\|_{L^{40}_y}^{1/5}.
		\]
		Substituting \eqref{F-L40-interpolated} and \eqref{gradF-L40-large}, and using
		\((1+\tau)^{1+\frac1{45}}\le C\tau^{1+\frac1{45}}\) for \(\tau\ge1\), we obtain
		\begin{equation}\label{F-Linfty-pointwise-40}
			\begin{aligned}
				\|\widetilde F(\tau)\|_{L^\infty_y}^4
				\le{}&
				C(\alpha)
				\Xi_0^{\frac{136}{225}}
				(1+\bar\rho E_0)^{\frac{2183}{1800}}\\
				&\times
				\left[
				\bar\rho^{1-\gamma}
				\sup_{0\le t\le T}
				\left(
				\|\nabla u(t)\|_{L^2}^2+\|P(t)\|_{L^2}^2
				\right)
				\right]^{\frac{19}{24}}\\
				&\times
				\tau^{-\frac{1627}{1800}}
				\left(
				\|\widetilde\rho^{1/2}\dot{\widetilde u}(\tau)\|_{L^2_y}
				+
				\|\nabla_y\dot{\widetilde u}(\tau)\|_{L^2_y}
				\right)^{1/5}.
			\end{aligned}
		\end{equation}
		Indeed, the power of \(\tau\) is
		\[
		-\frac{7}{24}\cdot\frac{19}{5}
		+
		\left(1+\frac1{45}\right)\frac15
		=
		-\frac{1627}{1800}.
		\]
		Integrating \eqref{F-Linfty-pointwise-40} over \([1,\bar\rho^\gamma T]\), and using
		H\"older's inequality with exponents \(10/9\) and \(10\), we obtain
		\[
		\begin{aligned}
			&\int_1^{\bar\rho^\gamma T}
			\|\widetilde F(\tau)\|_{L^\infty_y}^4\,d\tau\\
			\le& C(\alpha)\Xi_0^{\frac{136}{225}}(1+\bar\rho E_0)^{\frac{2183}{1800}}
			\left[\bar\rho^{1-\gamma}\sup_{0\le t\le T}\big(\|\nabla u(t)\|_{L^2}^2+\|P(t)\|_{L^2}^2\big)\right]^{\frac{19}{24}}\\
			&\times\left(\int_1^\infty \tau^{-\frac{1627}{1620}}\,d\tau\right)^{\frac{9}{10}}
			\left(\int_1^{\bar\rho^\gamma T}\big[\|\widetilde\rho^{1/2}\dot{\widetilde u}\|_{L^2_y}+\|\nabla_y\dot{\widetilde u}\|_{L^2_y}\big]^2d\tau\right)^{\frac{1}{10}}.
		\end{aligned}
		\]
		Since
		\[
		-\frac{1627}{1800}\cdot\frac{10}{9}
		=
		-\frac{1627}{1620}
		<
		-1,
		\]
		the first integral is finite. Moreover, by Corollary \ref{cor3.1.1} and
		Lemma \ref{lem3.5} in normalized variables,
		\[
		\int_1^{\bar\rho^\gamma T}
		\left[
		\|\widetilde\rho^{1/2}\dot{\widetilde u}\|_{L^2_y}
		+
		\|\nabla_y\dot{\widetilde u}\|_{L^2_y}
		\right]^2d\tau
		\le
		C
		\bar\rho^{1-\gamma}
		\left(
		\|\nabla u_0\|_{L^2}^2+\|P_0\|_{L^2}^2
		\right).
		\]
		Thus, using the fact that the initial value is bounded by the supremum,
		\begin{equation}\label{large-F-Linfty-pre-final}
			\begin{aligned}
				&\int_1^{\bar\rho^\gamma T}
				\|\widetilde F(\tau)\|_{L^\infty_y}^4\,d\tau\le
				C(\alpha)
				\Xi_0^{\frac{136}{225}}
				(1+\bar\rho E_0)^{\frac{2183}{1800}}
				\left[
				\bar\rho^{1-\gamma}
				\sup_{0\le t\le T}
				\left(
				\|\nabla u(t)\|_{L^2}^2+\|P(t)\|_{L^2}^2
				\right)
				\right]^{\frac{107}{120}}.
			\end{aligned}
		\end{equation}
		By \eqref{zlotnik-ST-small}, the last bracket is not larger than \(1\). Since
		\[
		\frac{136}{225}<\frac{6}{7},
		\qquad
		\frac{2183}{1800}<\frac{9}{7},
		\qquad
		\frac{107}{120}>\frac67,
		\]
		and \(\Xi_0\ge1\), \(1+\bar\rho E_0\ge1\), we infer from
		\eqref{large-F-Linfty-pre-final} that
		\begin{equation}\label{large-F-Linfty-bound-new}
			\begin{aligned}
				&\int_1^{\bar\rho^\gamma T}
				\|\widetilde F(\tau)\|_{L^\infty_y}^4\,d\tau\\
				&\le
				C(\alpha)
				\left[
				\Xi_0
				(1+\bar\rho E_0)^{3/2}
				\bar\rho^{1-\gamma}
				\sup_{0\le t\le T}
				\left(
				\|\nabla u(t)\|_{L^2}^2+\|P(t)\|_{L^2}^2
				\right)
				\right]^{6/7}\\
				&\le C\varepsilon^{6/7}.
			\end{aligned}
		\end{equation}
		
		We now verify the condition for Zlotnik's inequality. Fix arbitrary
		\[
		0\le \tau\le \bar\rho^\gamma T,
		\qquad y\in\mathbb R^2.
		\]
		Let \(X(s;y,\tau)\) be the particle trajectory generated by \(\widetilde u\):
		\[
		\frac{d}{ds}X(s;y,\tau)=\widetilde u(X(s;y,\tau),s),
		\qquad
		X(\tau;y,\tau)=y.
		\]
		Along this trajectory, the continuity equation gives
		\[
		\frac{d}{ds}\widetilde\rho(X(s;y,\tau),s)
		=
		-\widetilde\rho\,\div_y\widetilde u(X(s;y,\tau),s).
		\]
		Since
		\[
		\div_y\widetilde u
		=
		\frac{\widetilde F+\widetilde P}{2\mu+\lambda},
		\qquad
		\widetilde P=\widetilde\rho^\gamma,
		\]
		we obtain
		\[
		\frac{d}{ds}\widetilde\rho(X(s;y,\tau),s)
		=
		-\frac{\widetilde\rho^{\gamma+1}}{2\mu+\lambda}(X(s;y,\tau),s)
		-
		\frac{\widetilde\rho\widetilde F}{2\mu+\lambda}(X(s;y,\tau),s).
		\]
		This is the Zlotnik form
		\[
		\frac{d}{ds}\widetilde\rho(X(s;y,\tau),s)
		=
		g\bigl(\widetilde\rho(X(s;y,\tau),s)\bigr)+b'(s),
		\]
		where
		\[
		g(\zeta):=-\frac{\zeta^{\gamma+1}}{2\mu+\lambda},
		\]
		and
		\[
		b(s):=
		-\frac1{2\mu+\lambda}
		\int_0^s
		\widetilde\rho(X(\ell;y,\tau),\ell)
		\widetilde F(X(\ell;y,\tau),\ell)\,d\ell .
		\]
		
		We estimate \(b(s_2)-b(s_1)\) on the short-time and large-time parts separately.
		
		If \(0\le s_1<s_2\le1\), then using \(0\le\widetilde\rho\le2\) and
		\eqref{short-F-L1-Linfty-small}, we have
		\[
		\begin{aligned}
			b(s_2)-b(s_1)
			&\le
			\frac2{2\mu+\lambda}
			\int_{s_1}^{s_2}
			\|\widetilde F(\ell)\|_{L^\infty_y}\,d\ell\\
			&\le
			C\varepsilon^{1/2}
			\le
			\frac14.
		\end{aligned}
		\]
		If \(1\le s_1<s_2\le \bar\rho^\gamma T\), then using \(0\le\widetilde\rho\le2\),
		H\"older's inequality, and \eqref{large-F-Linfty-bound-new}, we have
		\[
		\begin{aligned}
			b(s_2)-b(s_1)
			&\le
			\frac2{2\mu+\lambda}
			\int_{s_1}^{s_2}
			\|\widetilde F(\ell)\|_{L^\infty_y}\,d\ell\\
			&\le
			C
			\left(
			\int_1^{\bar\rho^\gamma T}
			\|\widetilde F(\ell)\|_{L^\infty_y}^4\,d\ell
			\right)^{1/4}
			(s_2-s_1)^{3/4}\\
			&\le
			C\varepsilon^{3/14}(s_2-s_1)^{3/4}.
		\end{aligned}
		\]
		By Young's inequality,
		\[
		C\varepsilon^{3/14}(s_2-s_1)^{3/4}
		\le
		C\varepsilon^{6/7}
		+
		\frac1{4(2\mu+\lambda)}(s_2-s_1),
		\]
		and hence, by \eqref{zlotnik-eps-small},
		\[
		b(s_2)-b(s_1)
		\le
		\frac14
		+
		\frac1{4(2\mu+\lambda)}(s_2-s_1).
		\]
		Therefore, for arbitrary \(0\le s_1<s_2\le\bar\rho^\gamma T\), splitting the interval
		\([s_1,s_2]\) into its intersection with \([0,1]\) and with
		\([1,\bar\rho^\gamma T]\), the preceding two estimates give
		\begin{equation}\label{global-b-increment}
			b(s_2)-b(s_1)
			\le
			\frac12
			+
			\frac1{4(2\mu+\lambda)}(s_2-s_1).
		\end{equation}
		
		For every \(\zeta\ge1\), one has
		\[
		g(\zeta)
		=
		-\frac{\zeta^{\gamma+1}}{2\mu+\lambda}
		\le
		-\frac1{2\mu+\lambda}
		\le
		-\frac1{4(2\mu+\lambda)}.
		\]
		Thus Zlotnik's inequality applies on \([0,\tau]\) with
		\[
		\bar\zeta=1,
		\qquad
		N_0=\frac12,
		\qquad
		N_1=\frac1{4(2\mu+\lambda)}.
		\]
		Since
		\[
		0\le\widetilde\rho_0\le1,
		\]
		Zlotnik's inequality yields
		\[
		\widetilde\rho(y,\tau)
		\le
		\max\{\|\widetilde\rho_0\|_{L^\infty},1\}
		+
		\frac12
		\le
		\frac32.
		\]
		Since \(y\in\mathbb R^2\) and \(\tau\in[0,\bar\rho^\gamma T]\) are arbitrary, we get
		\[
		0\le\widetilde\rho(y,\tau)\le\frac32,
		\qquad
		(y,\tau)\in\mathbb R^2\times[0,\bar\rho^\gamma T].
		\]
		Returning to the original variables,
		\[
		\rho(x,t)
		=
		\bar\rho\,
		\widetilde\rho\!\left(
		\bar\rho^{\frac{\gamma+1}{2}}x,\bar\rho^\gamma t
		\right),
		\]
		we obtain
		\[
		0\le\rho(x,t)\le\frac32\bar\rho,
		\qquad
		(x,t)\in\mathbb R^2\times[0,T].
		\]
		This proves \eqref{uniform-density-bound}. The proof is complete.
	\end{proof}

	\begin{corollary}\label{cor:scale-invariant-smallness-closed}
		Let \((\rho,u)\) solve \eqref{1.1}--\eqref{1.4.1} smoothly on
		\(\mathbb R^2\times(0,T]\) satisfying \eqref{3.1.3} and the assumptions of
		Theorem \ref{th1.1}. Then there exists a positive constant \(c_0\), depending only on
		\(\mu,\lambda,\gamma,\alpha\), such that if
		\begin{equation}\label{epsilon-epsilon0-relation}
			\varepsilon\ge c_0\varepsilon_0,
		\end{equation}
		then, for any \(t\in[0,T]\),
		\begin{equation}\label{scale-smallness-closed}
			\begin{aligned}
				&\Xi_0(1+\bar\rho E_0)^{3/2}
				\bar\rho^{1-\gamma}
				\sup_{0\le s\le t}
				\left(
				\|\nabla u(s)\|_{L^2}^2+\|P(s)\|_{L^2}^2
				\right)\le
				\frac32\varepsilon .
			\end{aligned}
		\end{equation}
	\end{corollary}
	
	\begin{proof}
		By Corollary \ref{cor3.1.1}, we have
		\begin{equation}\label{low-order-close-used}
			\begin{aligned}
				&\sup_{0\le s\le t}
				\left(
				\|\nabla u(s)\|_{L^2}^2+\|P(s)\|_{L^2}^2
				\right)
				+
				\int_0^t
				\int_{\mathbb R^2}\rho|\dot u|^2\,dxds
				+
				\int_0^t
				\|P(s)\|_{L^3}^3\,ds\\
				\le&
				C
				\left(
				\|\nabla u_0\|_{L^2}^2+\|P_0\|_{L^2}^2
				\right),
			\end{aligned}
		\end{equation}
		where \(C\) depends only on \(\mu,\lambda,\gamma\). Hence
		\[
		\begin{aligned}
			&\Xi_0(1+\bar\rho E_0)^{3/2}
			\bar\rho^{1-\gamma}
			\sup_{0\le s\le t}
			\left(
			\|\nabla u(s)\|_{L^2}^2+\|P(s)\|_{L^2}^2
			\right)\\
			\le&
			C
			\Xi_0(1+\bar\rho E_0)^{3/2}
			\bar\rho^{1-\gamma}
			\left(
			\|\nabla u_0\|_{L^2}^2+\|P_0\|_{L^2}^2
			\right).
		\end{aligned}
		\]
		Using the initial smallness condition \eqref{initial-scale-smallness}, we obtain
		\[
		\begin{aligned}
			\Xi_0(1+\bar\rho E_0)^{3/2}
			\bar\rho^{1-\gamma}
			\sup_{0\le s\le t}
			\left(
			\|\nabla u(s)\|_{L^2}^2+\|P(s)\|_{L^2}^2
			\right)
			\le
			C\varepsilon_0.
		\end{aligned}
		\]
		Now choose
		\[
		c_0:=\frac{2C}{3}.
		\]
		Then \eqref{epsilon-epsilon0-relation} implies
		\[
		C\varepsilon_0\le \frac32\varepsilon.
		\]
		Therefore,
		\[
		\begin{aligned}
			&\Xi_0(1+\bar\rho E_0)^{3/2}
			\bar\rho^{1-\gamma}
			\sup_{0\le s\le t}
			\left(
			\|\nabla u(s)\|_{L^2}^2+\|P(s)\|_{L^2}^2
			\right)
			\le
			\frac32\varepsilon.
		\end{aligned}
		\]
		This proves \eqref{scale-smallness-closed}. The proof is complete.
	\end{proof}

	\section{Higher Order A Priori Estimates}
	From now on, assume that \((\rho,u)\) is a smooth solution of
	\eqref{1.1}--\eqref{1.4.1} on \(\mathbb R^2\times(0,T]\) satisfying
	\eqref{3.1.3} with smooth initial data \((\rho_0,u_0)\) satisfying the assumptions \eqref{assumption}
	of Theorem \ref{th1.1}. Moreover, besides the parameters
	\(T,\ \mu,\lambda,\gamma,\ q,\ \alpha,\ \bar\rho,\ M_0^{-1}\),
	the generic positive constant \(C\) may also depend on \(\|\bar x^\alpha\rho_0\|_{L^1\cap H^1\cap W^{1,q}},
	\|u_0\|_{D^1},
	\|\rho_0^{1/2}u_0\|_{L^2}\) appearing in Theorem \ref{th1.1}.

	The following higher order estimates are needed only to extend the local classical
	solution to a global one. Their proofs are the same as those in
	\cite[Lemmas 4.1--4.3]{Li-xin-2019}, hence we omit the details.
	
	\begin{lemma}\label{lem4.1}
		Under the assumptions of Theorem \ref{th1.1}, the following estimates hold:
		\begin{equation}\label{high-est-1}
			\begin{aligned}
				&\sup_{0\le t\le T}
				\left(
				\|\rho(t)\|_{H^1\cap W^{1,q}}
				+
				\|\nabla u(t)\|_{L^2}
				+
				t\|\nabla^2u(t)\|_{L^2}^2
				\right)\\
				&+
				\int_0^T
				\left(
				\|\nabla^2u\|_{L^2}^2
				+
				\|\nabla^2u\|_{L^q}^{\frac{q+1}{q}}
				+
				t\|\nabla^2u\|_{L^q}^2
				\right)\,dt
				\le C,
			\end{aligned}
		\end{equation}
		\begin{equation}\label{high-est-weighted-rho}
			\sup_{0\le t\le T}
			\|\bar x^\alpha\rho(t)\|_{L^1\cap H^1\cap W^{1,q}}
			\le C,
		\end{equation}
		\begin{equation}\label{high-est-ut-time-weighted}
			\sup_{0\le t\le T}
			t
			\left(
			\|\rho^{1/2}u_t(t)\|_{L^2}^2
			+
			\|\nabla u(t)\|_{H^1}^2
			\right)
			+
			\int_0^T
			t\|\nabla u_t\|_{L^2}^2\,dt
			\le C.
		\end{equation}
	\end{lemma}
	
	Similarly,  assume that
	\((\rho,u)\) is a smooth solution of \eqref{1.1}--\eqref{1.4.1} on
	\(\mathbb R^2\times(0,T]\) satisfying \eqref{3.1.3} with smooth initial data
	\((\rho_0,u_0)\) satisfying  assumptions \eqref{assumption} of Theorem \ref{th1.1}. In addition,
	assume that the initial data satisfy the extra hypotheses \eqref{1.12.1} of Theorem \ref{th1.2}. In addition,
	besides the parameters already allowed in Theorem \ref{th1.1},
	the generic positive constant \(C\) may also depend on
	\(\delta_0,
	\|\nabla^2\rho_0\|_{L^2\cap L^q},
	\|\nabla^2P_0\|_{L^2\cap L^q},
	\|\bar x^{\delta_0}\nabla^2\rho_0\|_{L^2},
	\|\bar x^{\delta_0}\nabla^2P_0\|_{L^2},
	\|\nabla^2u_0\|_{L^2},
	\|g\|_{L^2}.\)
	
	Due to \cite{Liang2014 2D}, we also omit the details and give directly the following estimates.
	\begin{lemma}\label{lem4.2}
		Under the assumptions of Theorem \ref{th1.2}, the following estimates hold:
		\begin{equation}\label{high-thm2-ut-unweighted}
			\sup_{0\le t\le T}
			\left(
			\|\rho^{1/2}u_t(t)\|_{L^2}
			+
			\|\nabla u(t)\|_{H^1}
			\right)
			+
			\int_0^T
			\left(
			\|\nabla u_t\|_{L^2}^2
			+
			\|u_t\bar x^{-1}\|_{L^2}^2
			\right)\,dt
			\le C,
		\end{equation}
		\begin{equation}\label{high-thm2-weighted-rho-second}
			\sup_{0\le t\le T}
			\left(
			\|\bar x^{\delta_0}\nabla^2\rho(t)\|_{L^2}
			+
			\|\bar x^{\delta_0}\nabla^2P(t)\|_{L^2}
			\right)
			\le C,
		\end{equation}
		\begin{equation}\label{high-thm2-ut-time}
			\sup_{0\le t\le T}
			t\|\nabla u_t(t)\|_{L^2}^2
			+
			\int_0^T
			t
			\left(
			\|\rho^{1/2}u_{tt}\|_{L^2}^2
			+
			\|\nabla^2u_t\|_{L^2}^2
			\right)\,dt
			\le C,
		\end{equation}
		\begin{equation}\label{high-thm2-rho-second-Lq}
			\sup_{0\le t\le T}
			\left(
			\|\nabla^2\rho(t)\|_{L^q}
			+
			\|\nabla^2P(t)\|_{L^q}
			\right)
			\le C,
		\end{equation}
		\begin{equation}\label{high-thm2-top}
			\begin{aligned}
				&\sup_{0\le t\le T}
				t\Big(
				\|\rho^{1/2}u_{tt}(t)\|_{L^2}
				+
				\|\nabla^3u(t)\|_{L^2\cap L^q}
				+
				\|\nabla u_t(t)\|_{H^1}
				+
				\|\nabla^2(\rho u)(t)\|_{L^{\frac{q+2}{2}}}
				\Big)\\
				&+
				\int_0^T
				t^2
				\left(
				\|\nabla u_{tt}\|_{L^2}^2
				+
				\|u_{tt}\bar x^{-1}\|_{L^2}^2
				\right)\,dt
				\le C.
			\end{aligned}
		\end{equation}
	\end{lemma}
	
	\section{Proofs of Theorems \ref{th1.1} and \ref{th1.2}}

	\noindent With the a priori estimates established in Sections 3 and 4, we proceed to prove the main theorems of this paper in this section.
	
	\noindent\textbf{The proof of Theorem \ref{th1.1}. }	
	By Lemma \ref{lem2.1}, there exists a \(T_*>0\) such that the Cauchy problem
	\eqref{1.1}--\eqref{1.4.1} has a unique strong solution \((\rho,u)\) on
	\(\mathbb R^2\times(0,T_*]\). We will use the a priori estimates,
	Proposition \ref{p3.1} to extend the local strong solution \((\rho,u)\)
	to all time.
	
	Recall that
	\[
	\Xi_0:=
	\bar\rho^\gamma M_0
	+
	(\bar\rho^\gamma M_0)^{-1}
	+
	\bar\rho^{\gamma+\frac{\alpha(\gamma+1)}2}
	\|\rho_0|x|^\alpha\|_{L^1},
	\]
	and
	\[
	I_0
	:=
	\Xi_0(1+\bar\rho E_0)^{3/2}
	\bar\rho^{1-\gamma}
	\left(
	\|\nabla u_0\|_{L^2}^2+\|\rho_0^\gamma\|_{L^2}^2
	\right).
	\]
	Let \(\varepsilon_1\in(0,1]\) be the small constant in Proposition
	\ref{p3.1}, and let \(c_0>0\) be the constant in Corollary
	\ref{cor:scale-invariant-smallness-closed}. We choose the small constant
	\(\varepsilon_0\) in Theorem \ref{th1.1} sufficiently small such that
	\[
	2c_0\varepsilon_0\le \varepsilon_1.
	\]
	Set \(\varepsilon:=2c_0\varepsilon_0.\) Then
	\[
	I_0\le \varepsilon_0<2\varepsilon,
	\qquad
	\varepsilon\ge c_0\varepsilon_0,
	\qquad
	\varepsilon\le \varepsilon_1.
	\]
	
	For any \(T>0\) on which the strong solution exists, define
	\begin{equation}\label{def-IT}
		I_T
		:=
		\Xi_0(1+\bar\rho E_0)^{3/2}
		\bar\rho^{1-\gamma}
		\sup_{0\le t\le T}
		\left(
		\|\nabla u(t)\|_{L^2}^2+\|P(t)\|_{L^2}^2
		\right).
	\end{equation}
	First, since
	\[
	0\le \rho_0\le \bar\rho,
	\qquad
	I_0\le \varepsilon_0<2\varepsilon,
	\]
	there exists a \(T_1\in(0,T_*]\) such that
	\begin{equation}\label{bootstrap-holds-T1}
		\sup_{\mathbb R^2\times[0,T_1]}\rho\le2\bar\rho,
		\qquad
		I_{T_1}\le2\varepsilon.
	\end{equation}
	
	Next, set
	\begin{equation}\label{def-T-star}
		T^*
		:=
		\sup\left\{
		T>0\ \middle|\
		\text{the strong solution exists on }[0,T]
		\text{ and }
		\sup_{\mathbb R^2\times[0,T]}\rho\le2\bar\rho,\
		I_T\le2\varepsilon
		\right\}.
	\end{equation}
	Then \(T^*\ge T_1>0\). As in the proof of Theorem 1.1 in \cite{Li-xin-2019}, we get $T^*=+\infty$. This completes the
	proof of Theorem \ref{th1.1}. $\hfill\square $
	
	\noindent\textbf{The proof of Theorem \ref{th1.2}. } By the same argument as Theorem \ref{th1.1}, Theorem \ref{th1.2} can be established using Lemma \ref{lem2.1}, Proposition \ref{p3.1}, and Lemmas \ref{lem4.1}--\ref{lem4.2}. $\hfill\square $
	
	\section{Appendix: Scaling-invariant quantities}
	
	In this appendix, we verify the scaling invariance of the initial quantity
	\(I_0\) and explain why the power of \(\Xi_0\) used in \eqref{definition-I0}
	does not lead to a data-independent positive lower bound.
	The following lemmas show that scaling invariance alone is not enough: a
	candidate invariant quantity must also be capable of becoming genuinely small.
	
	\begin{lemma}\label{lem6.1}
		The quantity
		\[
		I_0
		:=
		\Xi_0(1+\bar\rho E_0)^{3/2}
		\bar\rho^{1-\gamma}
		\left(
		\|\nabla u_0\|_{L^2}^2+\|\rho_0^\gamma\|_{L^2}^2
		\right)
		\]
		is invariant under the natural scaling of the two-dimensional isentropic
		compressible Navier--Stokes equations.
	\end{lemma}
	
	\begin{proof}
		The system \eqref{1.1} is invariant under the scaling
		\begin{equation}\label{eq:scaling-2d}
			\rho_\ell(x,t)
			=
			\ell^{1/\gamma}
			\rho\!\left(
			\ell^{\frac{\gamma+1}{2\gamma}}x,\ell t
			\right),
			\qquad
			u_\ell(x,t)
			=
			\ell^{\frac{\gamma-1}{2\gamma}}
			u\!\left(
			\ell^{\frac{\gamma+1}{2\gamma}}x,\ell t
			\right),
		\end{equation}
		for any \(\ell>0\). The corresponding initial data satisfy
		\[
		\rho_{0,\ell}(x)
		=
		\ell^{1/\gamma}
		\rho_0\!\left(
		\ell^{\frac{\gamma+1}{2\gamma}}x
		\right),
		\qquad
		u_{0,\ell}(x)
		=
		\ell^{\frac{\gamma-1}{2\gamma}}
		u_0\!\left(
		\ell^{\frac{\gamma+1}{2\gamma}}x
		\right).
		\]
		If \(0\le \rho_0\le \bar\rho\), then
		\[
		0\le \rho_{0,\ell}(x)\le \bar\rho_\ell,
		\qquad
		\bar\rho_\ell:=\ell^{1/\gamma}\bar\rho.
		\]
		Let
		\[
		M_{0,\ell}:=\int_{\mathbb R^2}\rho_{0,\ell}(x)\,dx,
		\qquad
		E_{0,\ell}
		:=
		\int_{\mathbb R^2}
		\left(
		\frac12\rho_{0,\ell}|u_{0,\ell}|^2
		+
		\frac{1}{\gamma-1}\rho_{0,\ell}^\gamma
		\right)\,dx .
		\]
		Changing variables
		\(z=\ell^{\frac{\gamma+1}{2\gamma}}x\), we have
		\[
		M_{0,\ell}
		=
		\ell^{1/\gamma}
		\ell^{-\frac{\gamma+1}{\gamma}}
		M_0
		=
		\ell^{-1}M_0.
		\]
		Moreover,
		\[
		\begin{aligned}
			E_{0,\ell}
			&=
			\int_{\mathbb R^2}
			\left(
			\frac12
			\ell^{1/\gamma}\rho_0(z)
			\ell^{\frac{\gamma-1}{\gamma}}|u_0(z)|^2
			+
			\frac1{\gamma-1}
			\ell\rho_0^\gamma(z)
			\right)
			\ell^{-\frac{\gamma+1}{\gamma}}\,dz\\
			&=
			\ell^{-1/\gamma}
			\int_{\mathbb R^2}
			\left(
			\frac12\rho_0|u_0|^2
			+
			\frac1{\gamma-1}\rho_0^\gamma
			\right)\,dz\\
			&=
			\ell^{-1/\gamma}E_0.
		\end{aligned}
		\]
		Hence
		\[
		\bar\rho_\ell^\gamma M_{0,\ell}
		=
		(\ell^{1/\gamma}\bar\rho)^\gamma\ell^{-1}M_0
		=
		\bar\rho^\gamma M_0,
		\]
		\[
		(\bar\rho_\ell^\gamma M_{0,\ell})^{-1}
		=
		(\bar\rho^\gamma M_0)^{-1},
		\qquad
		\bar\rho_\ell E_{0,\ell}
		=
		\bar\rho E_0.
		\]
		For the weighted density term, we have
		\[
		\begin{aligned}
			\|\rho_{0,\ell}|x|^\alpha\|_{L^1}
			&=
			\int_{\mathbb R^2}
			\ell^{1/\gamma}
			\rho_0\!\left(
			\ell^{\frac{\gamma+1}{2\gamma}}x
			\right)
			|x|^\alpha\,dx\\
			&=
			\ell^{1/\gamma}
			\ell^{-\frac{\alpha(\gamma+1)}{2\gamma}}
			\ell^{-\frac{\gamma+1}{\gamma}}
			\int_{\mathbb R^2}\rho_0(z)|z|^\alpha\,dz\\
			&=
			\ell^{-1-\frac{\alpha(\gamma+1)}{2\gamma}}
			\|\rho_0|x|^\alpha\|_{L^1}.
		\end{aligned}
		\]
		Consequently,
		\[
		\begin{aligned}
			&\bar\rho_\ell^{\gamma+\frac{\alpha(\gamma+1)}2}
			\|\rho_{0,\ell}|x|^\alpha\|_{L^1}\\
			&=
			(\ell^{1/\gamma}\bar\rho)^{\gamma+\frac{\alpha(\gamma+1)}2}
			\ell^{-1-\frac{\alpha(\gamma+1)}{2\gamma}}
			\|\rho_0|x|^\alpha\|_{L^1}\\
			&=
			\bar\rho^{\gamma+\frac{\alpha(\gamma+1)}2}
			\|\rho_0|x|^\alpha\|_{L^1}.
		\end{aligned}
		\]
		Similarly, we define
		\[
		\Xi_{0,\ell}
		:=
		\bar\rho_\ell^\gamma M_{0,\ell}
		+
		(\bar\rho_\ell^\gamma M_{0,\ell})^{-1}
		+
		\bar\rho_\ell^{\gamma+\frac{\alpha(\gamma+1)}2}
		\|\rho_{0,\ell}|x|^\alpha\|_{L^1},
		\]
		then we obtain
		\[
		\Xi_{0,\ell}=\Xi_0.
		\]
		It remains to check the last term in \(I_0\). Since
		\[
		\nabla u_{0,\ell}(x)
		=
		\ell^{\frac{\gamma-1}{2\gamma}}
		\ell^{\frac{\gamma+1}{2\gamma}}
		\nabla u_0\!\left(
		\ell^{\frac{\gamma+1}{2\gamma}}x
		\right)
		=
		\ell
		\nabla u_0\!\left(
		\ell^{\frac{\gamma+1}{2\gamma}}x
		\right),
		\]
		we get
		\[
		\|\nabla u_{0,\ell}\|_{L^2}^2
		=
		\ell^2
		\ell^{-\frac{\gamma+1}{\gamma}}
		\|\nabla u_0\|_{L^2}^2
		=
		\ell^{\frac{\gamma-1}{\gamma}}
		\|\nabla u_0\|_{L^2}^2.
		\]
		Similarly,
		\[
		\|\rho_{0,\ell}^\gamma\|_{L^2}^2
		=
		\ell^{\frac{\gamma-1}{\gamma}}
		\|\rho_0^\gamma\|_{L^2}^2 ,\quad
		\bar\rho_\ell^{1-\gamma}
		=
		(\ell^{1/\gamma}\bar\rho)^{1-\gamma}
		=
		\ell^{\frac{1-\gamma}{\gamma}}
		\bar\rho^{1-\gamma},
		\]
		implies that
		\[
		\begin{aligned}
			&\bar\rho_\ell^{1-\gamma}
			\left(
			\|\nabla u_{0,\ell}\|_{L^2}^2
			+
			\|\rho_{0,\ell}^\gamma\|_{L^2}^2
			\right)=
			\bar\rho^{1-\gamma}
			\left(
			\|\nabla u_0\|_{L^2}^2
			+
			\|\rho_0^\gamma\|_{L^2}^2
			\right).
		\end{aligned}
		\]
		Combining the above identities yields
		\[
		\begin{aligned}
			I_0[\rho_{0,\ell},u_{0,\ell}]
			&=
			\Xi_{0,\ell}
			(1+\bar\rho_\ell E_{0,\ell})^{3/2}
			\bar\rho_\ell^{1-\gamma}
			\left(
			\|\nabla u_{0,\ell}\|_{L^2}^2
			+
			\|\rho_{0,\ell}^\gamma\|_{L^2}^2
			\right)\\
			&=
			\Xi_0
			(1+\bar\rho E_0)^{3/2}
			\bar\rho^{1-\gamma}
			\left(
			\|\nabla u_0\|_{L^2}^2
			+
			\|\rho_0^\gamma\|_{L^2}^2
			\right)\\
			&=
			I_0[\rho_0,u_0].
		\end{aligned}
		\]
		The proof is complete.
	\end{proof}
	
	The next lemma explains the obstruction caused by taking a too high power of
	\(\Xi_0\). It shows that, although the resulting quantity is still invariant
	under the natural scaling, it cannot be used as a smallness condition with an
	arbitrarily small universal threshold.
	
	\begin{lemma}\label{lem6.2}
		Under the assumption of Theorem \ref{th1.1}, we set
		\[
		J_\theta
		:=
		\Xi_0^\theta(1+\bar\rho E_0)^{3/2}
		\bar\rho^{1-\gamma}
		\left(
		\|\nabla u_0\|_{L^2}^2+\|\rho_0^\gamma\|_{L^2}^2
		\right),
		\]
		for
		\(\theta\ge2\gamma+\frac{4(2\gamma-1)}{\alpha}.\)
		Then there exists a positive constant \(c=c(\alpha,\gamma)>0\), independent of
		the initial data and \(\bar\rho\), such that \(J_\theta\ge c\).
		Consequently, a smallness condition based on such a high power of \(\Xi_0\)
		cannot have an arbitrarily small threshold independent of the initial data.
	\end{lemma}
	
	\begin{proof}
		Let
		\[
		R:=\left(\frac{2\int_{\mathbb R^2}\rho_0|x|^\alpha\,dx}{M_0}\right)^{1/\alpha}.
		\]
		Then
		\[
		\int_{\mathbb R^2\setminus B_R}\rho_0\,dx
		\le
		R^{-\alpha}\int_{\mathbb R^2}\rho_0|x|^\alpha\,dx
		=\frac{M_0}{2},
		\qquad
		\int_{B_R}\rho_0\,dx\ge \frac{M_0}{2}.
		\]
		By H\"older's inequality,
		\[
		\frac{M_0}{2}
		\le
		|B_R|^{1-\frac1{2\gamma}}
		\left(\int_{B_R}\rho_0^{2\gamma}\,dx\right)^{\frac1{2\gamma}}.
		\]
		Then we have
		\[
		\int_{\mathbb R^2}\rho_0^{2\gamma}\,dx
		\ge c(\alpha,\gamma)M_0^{2\gamma}|B_R|^{-(2\gamma-1)}.
		\]
		Since \(|B_R|=CR^2\), the definition of \(R\) gives
		\begin{equation}\label{lem6.2-1}
			\int_{\mathbb R^2}\rho_0^{2\gamma}\,dx
			\ge
			c(\alpha,\gamma)
			M_0^{2\gamma+\frac{2(2\gamma-1)}{\alpha}}
			\left(\int_{\mathbb R^2}\rho_0|x|^\alpha\,dx\right)^{-\frac{2(2\gamma-1)}{\alpha}}.
		\end{equation}
		Multiplying \eqref{lem6.2-1} by \(\bar\rho^{1-\gamma}\), we obtain
		\begin{equation}\label{lem6.2-2}
			\begin{aligned}
				\bar\rho^{1-\gamma}\int_{\mathbb R^2}\rho_0^{2\gamma}\,dx
				&\ge c(\alpha,\gamma)
				\left(\bar\rho^\gamma M_0\right)^{2\gamma+\frac{2(2\gamma-1)}{\alpha}}\times
				\left(\bar\rho^{\gamma+\frac{\alpha(\gamma+1)}2}
				\int_{\mathbb R^2}\rho_0|x|^\alpha\,dx\right)^{-\frac{2(2\gamma-1)}{\alpha}}.
			\end{aligned}
		\end{equation}
		Indeed, the power of \(\bar\rho\) cancels because
		\[
		1-\gamma
		-\gamma\left(2\gamma+\frac{2(2\gamma-1)}{\alpha}\right)
		+\left(\gamma+\frac{\alpha(\gamma+1)}2\right)
		\frac{2(2\gamma-1)}{\alpha}
		=0.
		\]
		By the definition of \(\Xi_0\), one has
		\begin{equation}\label{lem6.2-3}
			\Xi_0\ge
			\left(\bar\rho^\gamma M_0\right)^{-1},
		\end{equation}
		and
		\begin{equation}\label{lem6.2-4}
			\Xi_0\ge
			\bar\rho^{\gamma+\frac{\alpha(\gamma+1)}2}
			\int_{\mathbb R^2}\rho_0|x|^\alpha\,dx.
		\end{equation}
		Now multiply \eqref{lem6.2-2} by
		\(
		\Xi_0^{2\gamma+\frac{4(2\gamma-1)}{\alpha}},
		\)
		we get
		\begin{equation}\label{lem6.2-5}
			\begin{aligned}
				&\Xi_0^{2\gamma+\frac{4(2\gamma-1)}{\alpha}}
				\bar\rho^{1-\gamma}
				\int_{\mathbb R^2}\rho_0^{2\gamma}\,dx \ge
				c(\alpha,\gamma)
				\left(
				\Xi_0\bar\rho^\gamma M_0
				\right)^{
					2\gamma+\frac{2(2\gamma-1)}{\alpha}}
				\left(
				\frac{\Xi_0}{
					\bar\rho^{\gamma+\frac{\alpha(\gamma+1)}2}
					\int_{\mathbb R^2}\rho_0|x|^\alpha\,dx}
				\right)^{\frac{2(2\gamma-1)}{\alpha}}.
			\end{aligned}
		\end{equation}
		By \eqref{lem6.2-3} and \eqref{lem6.2-4}, both factors on the right-hand side of
		\eqref{lem6.2-5} are bounded from below by \(1\). Hence
		it follows that, for \(\theta\ge 2\gamma+\frac{4(2\gamma-1)}{\alpha},\)
		\[
		\Xi_0^\theta
		\bar\rho^{1-\gamma}
		\int_{\mathbb R^2}\rho_0^{2\gamma}\,dx
		\ge
		c(\alpha,\gamma).
		\]
		Finally, we obtain
		\[
		\begin{aligned}
			J_\theta
			&=
			\Xi_0^\theta(1+\bar\rho E_0)^{3/2}
			\bar\rho^{1-\gamma}
			\left(
			\|\nabla u_0\|_{L^2}^2+\|\rho_0^\gamma\|_{L^2}^2
			\right) \ge
			\Xi_0^\theta
			\bar\rho^{1-\gamma}
			\int_{\mathbb R^2}\rho_0^{2\gamma}\,dx
			\ge
			c(\alpha,\gamma).
		\end{aligned}
		\]
		This proves the lemma.
	\end{proof}
	
	\begin{lemma}\label{lem6.3}
		The scaling-invariant quantity \(I_0\) defined in \eqref{definition-I0} has no
		positive lower bound depending only on \(\alpha\) and \(\gamma\). More precisely,
		there exists a sequence of smooth initial data
		\(\{(\rho_0^{(n)},u_0^{(n)})\}_{n=1}^\infty\) with
		\(0\le\rho_0^{(n)}\le1\), \(0<M_0^{(n)}<\infty\), and
		\(u_0^{(n)}=0\), such that
		\[
		I_0[\rho_0^{(n)},u_0^{(n)}]\to0
		\qquad\text{as } n\to\infty.
		\]
	\end{lemma}
	
	\begin{proof}
		Let \(0\le\varphi\in C_c^\infty(B_2\setminus B_1)\)  and
		\(0<\|\varphi\|_{L^\infty}\le1\). For \(n\ge1\), set
		\[
		\rho_0^{(n)}(x):=n^{-(\alpha+2)}\varphi(x/n),
		\qquad
		u_0^{(n)}(x):=0,
		\qquad
		\bar\rho=1.
		\]
		Then \(0\le\rho_0^{(n)}\le1\). Direct changes of variables give
		\[
		M_0^{(n)}=\int_{\mathbb R^2}\rho_0^{(n)}\,dx=C(\varphi) n^{-\alpha},
		\qquad
		\int_{\mathbb R^2}\rho_0^{(n)}|x|^\alpha\,dx=C(\varphi,\alpha),
		\]
		\[
		E_0^{(n)}=C(\varphi,\gamma)n^{2-\gamma(\alpha+2)},
		\qquad
		\|(\rho_0^{(n)})^\gamma\|_{L^2}^2
		=C(\varphi,\gamma)n^{2-2\gamma(\alpha+2)}.
		\]
		Consequently,
		\[
		\Xi_0^{(n)}
		=
		M_0^{(n)}+(M_0^{(n)})^{-1}
		+
		\int_{\mathbb R^2}\rho_0^{(n)}|x|^\alpha\,dx
		\le
		C n^\alpha,
		\qquad
		1+E_0^{(n)}\le C,
		\]
		for all \(n\ge1\). Since \(u_0^{(n)}=0\), we obtain
		\[
		\begin{aligned}
			I_0[\rho_0^{(n)},u_0^{(n)}]
			&=
			\Xi_0^{(n)}(1+E_0^{(n)})^{3/2}
			\|(\rho_0^{(n)})^\gamma\|_{L^2}^2\\
			&\le
			C n^{\alpha}
			n^{2-2\gamma(\alpha+2)}
			=
			C n^{-(2\gamma-1)(\alpha+2)}.
		\end{aligned}
		\]
		Because \(\gamma>1\) and \(\alpha>1\), the exponent
		\((2\gamma-1)(\alpha+2)>0\). Hence
		\(I_0[\rho_0^{(n)},u_0^{(n)}]\to0\) as \(n\to\infty\). This proves that the present
		choice of the first power of \(\Xi_0\) does not produce the obstruction
		described in Lemma \ref{lem6.2}.
	\end{proof}
	
	To compare the two smallness assumptions, fix in Li--Xin \cite{Li-xin-2019} the upper density bound \(\bar\rho=1\).
	Denote the corresponding energy threshold by \(\varepsilon_*>0\), thus the condition used for comparison is
	\(E_0\le\varepsilon_*\). The comparison below concerns only these scalar
	conditions; the regularity and compatibility assumptions are kept separate.
	
	\begin{lemma}[Comparison of scalar smallness conditions]\label{lem6.4}
		The conditions \(E_0\le\varepsilon_*\) and \(I_0\le\varepsilon_0\) are not
		equivalent. More precisely:
		\begin{enumerate}
			\item[(i)] For fixed \(\rho_*\), \(M_0=\int\rho_*dx\), and \(\bar\rho=1\),
			the condition \(E_0\le\varepsilon_*\) does not imply
			\(I_0\le\varepsilon_0\).
			\item[(ii)] There exists a sequence \(\{(\rho_n,u_n)\}_{n\ge1}\) with
			\(\bar\rho_n=1\) and \(M_n:=\int\rho_n dx\) such that
			\[
			I_0[\rho_n,u_n]\longrightarrow0,
			\qquad E_0[\rho_n,u_n]\longrightarrow+\infty.
			\]
			Consequently, \(I_0[\rho_n,u_n]\le\varepsilon_0\) does not imply
			\(E_0[\rho_n,u_n]\le\varepsilon_*\).
		\end{enumerate}
	\end{lemma}
	
	\begin{proof}
		For (i), let \(\rho_*(x)=a e^{-|x|^2}\), where \(a\in(0,1]\) is chosen so that
		\[
		\frac{a^\gamma}{\gamma-1}\int_{\mathbb R^2}e^{-\gamma|x|^2}\,dx
		<\frac{\varepsilon_*}{2}.
		\]
		Fix \(\psi\in C_c^\infty(B_1;\mathbb R^2)\) with
		\(\|\nabla\psi\|_{L^2}>0\), set \(y_n=(4\log n,0)\), and define
		\(u_n(x)=n\psi(x-y_n)\). For \(x\in\operatorname{supp}u_n\),
		\(|x|\ge4\log n-1\). Hence, for all sufficiently large \(n\),
		\[
		\rho_*(x)\le a e^{-(4\log n-1)^2}\le a n^{-6}.
		\]
		Therefore
		\[
		\frac12\int\rho_*|u_n|^2dx
		\le \frac{a n^2}{2n^6}\|\psi\|_{L^2}^2
		=\frac{a}{2}\|\psi\|_{L^2}^2n^{-4},
		\]
		and consequently \(E_0[\rho_*,u_n]<\varepsilon_*\) for all sufficiently large \(n\).
		On the other hand, translation invariance gives
		\[
		\|\nabla u_n\|_{L^2}^2
		=n^2\|\nabla\psi\|_{L^2}^2.
		\]
		Since \(\rho_*\), \(M_0\), and \(\Xi_0\) are fixed and
		\(\|\rho_*^\gamma\|_{L^2}>0\),
		\[
		I_0[\rho_*,u_n]
		\ge \Xi_*(1+E_0[\rho_*,u_n])^{3/2}
		\left(\|\nabla u_n\|_{L^2}^2+\|\rho_*^\gamma\|_{L^2}^2\right)
		\ge  n^2\|\nabla\psi\|_{L^2}^2\longrightarrow+\infty.
		\]
		
		For (ii), let \(0\le\varphi\in C_c^\infty(B_2\setminus B_1)\),
		\(m:=\int\varphi dx\), \(R_n:=n^{2/\alpha}\), and
		\[
		\rho_n(x):=\frac{n^{-1}}{mR_n^2}\varphi(x/R_n),
		\qquad
		s:=2\gamma+\frac{8\gamma-4}{\alpha}>2.
		\]
		Define
		\[
		c_\alpha:=m^{-1}\int\varphi(y)|y|^\alpha dy,\quad
		A_\gamma:=m^{-\gamma}\int\varphi(y)^\gamma dy,\quad
		B_\gamma:=m^{-2\gamma}\int\varphi(y)^{2\gamma}dy.
		\]
		\[
		\int\rho_n dx=n^{-1},\quad
		\int\rho_n|x|^\alpha dx=c_\alpha n,\quad
		\Xi_n\le Cn,
		\quad \|\rho_n^\gamma\|_2^2=B_\gamma n^{-s},\quad
		\int\rho_n^\gamma dx=A_\gamma n^{-\kappa_0},
		\]
		where \(\kappa_0=\gamma+4(\gamma-1)/\alpha>0\). Choose \(1<r<s/2\).
		Fix \(\eta\in C^\infty(\mathbb R)\) such that
		\(0\le\eta\le1\), \(\eta(t)=1\) for \(t\le0\), and \(\eta(t)=0\) for \(t\ge1\). Set
		\(\ell_n=\log(e^{n^{4r}}/(2R_n))\) and define the radial function
		\[
		\chi_n(x)=\begin{cases}
			1, & |x|\le2R_n,\\
			\eta\!\left(\dfrac{\log(|x|/(2R_n))}{\ell_n}\right), & 2R_n<|x|<e^{n^{4r}},\\
			0, & |x|\ge e^{n^{4r}}.
		\end{cases}
		\]
		Then \(\chi_n\in C_c^\infty(\mathbb R^2)\), \(0\le\chi_n\le1\), and
		\[
		\int_{\mathbb R^2}|\nabla\chi_n|^2dx\le\frac{2\pi\|\eta'\|_{L^\infty(\mathbb R)}^2}{\ell_n}
		\le Cn^{-4r},
		\]
		because \(\ell_n=n^{4r}-\log(2R_n)\ge\tfrac12n^{4r}\) for all large \(n\). Set \(u_n(x):=(n^{r/2}\chi_n(x),0)\).
		Since \(\operatorname{supp}\rho_n\subset B_{2R_n}\),
		\[
		E_0[\rho_n,u_n]=\frac12n^{r-1}+A_\gamma(\gamma-1)^{-1}n^{-\kappa_0},
		\qquad \|\nabla u_n\|_2^2\le Cn^{-3r}.
		\]
		For large \(n\), \((1+E_0)^{3/2}\le Cn^{\frac32(r-1)}\), and hence
		\[
		I_0[\rho_n,u_n]\le Cn^{\frac{3r-1}{2}}(n^{-3r}+n^{-s})
		\le C(n^{-r}+n^{2r-s})\to0,
		\qquad E_0[\rho_n,u_n]\ge\tfrac12n^{r-1}\to+\infty.
		\]
	\end{proof}

	\vspace{\baselineskip}	
	{\footnotesize
		
		\section*{Acknowledgements}
		The authors would like to thank Professor Huanyao Wen for some kindly and valuable comments.
		B. L\"{u} is supported by the National Natural Science Foundation of China (No.12371219),  and the Double-Thousand Plan of Jiangxi Province (No. jxsq2023201115).

		\section*{Declarations}
		
		\textbf{Conflict of interest }
		The authors declare that they have no conflict of interest.

	}
\end{document}